%% file: 0paper.tex
\documentclass[12pt,letterpaper]{article}

\usepackage{graphicx}

\usepackage{verbatim}
\usepackage{amsmath}
\numberwithin{equation}{section} 
\numberwithin{figure}{section} 

\usepackage{color}

\usepackage{amssymb}

\usepackage{hyperref}

\usepackage[margin=1in]{geometry}

\title{Billiards in Nearly Isosceles Triangles}
\author{W. Patrick Hooper \thanks{Supported by N.S.F. Postdoctoral
Fellowship DMS-0803013} \hskip 5 pt and
Richard Evan Schwartz \thanks{
Supported by 
N.S.F. Research Grant DMS-0305047}}
\date{April, 23 2013\footnote{This improves the published paper \cite{HS}; it corrects some errors in \S 9. See remark \ref{rem:corrections}.}}

\newtheorem{theorem}{Theorem}[section]
\newtheorem{proposition}[theorem]{Proposition}
\newtheorem{lemma}[theorem]{Lemma}

\newtheorem{remarks}[theorem]{Remarks}
\newtheorem{corollary}[theorem]{Corollary}
\newtheorem{conjecture}[theorem]{Conjecture}

\def\startproof{{\bf {\medskip}{\noindent}Proof: }}

\def\endproof{$\spadesuit$  \newline}

\def\C{\mbox{\boldmath{$C$}}}%
\def\D{\mbox{\boldmath{$D$}}}%
\def\N{\mbox{\boldmath{$N$}}}%
\def\Q{\mbox{\boldmath{$Q$}}}%
\def\R{\mbox{\boldmath{$R$}}}%
\def\Z{\mbox{\boldmath{$Z$}}}%

\begin{document}
\maketitle

\begin{abstract}
We prove that any sufficiently small perturbation
of an isosceles triangle has a periodic
billiard path.  Our proof involves the analysis
of certain infinite families of Fourier series
that arise in connection with triangular
billiards, and reveals some self-similarity
phenomena in irrational triangular
billiards.  Our analysis illustrates the surprising
fact that billiards on a triangle {\it near\/}
a Veech triangle is extremely complicated even
though billiards {\it on\/} a Veech triangle
is very well understood.
\end{abstract}

\input{1intro}

\input{2back}
\input{3pat}
\input{4abc}

\input{5word}
\input{6qrt}
\input{7pivot}
\input{8tile}

\input{9stability}

\bibliographystyle{amsalpha}
\bibliography{bibliography}

\end{document}

%% file: 1intro.tex
\section{Introduction}

This paper concerns periodic billiard paths in
triangles.  In some cases, quite a bit is known about
periodic billiard paths in triangles, and in some cases
surprisingly little is known.  For example,
it is still unknown if every triangle has a periodic
billiard path.

\begin{center}
\resizebox{!}{1.4in}{\includegraphics{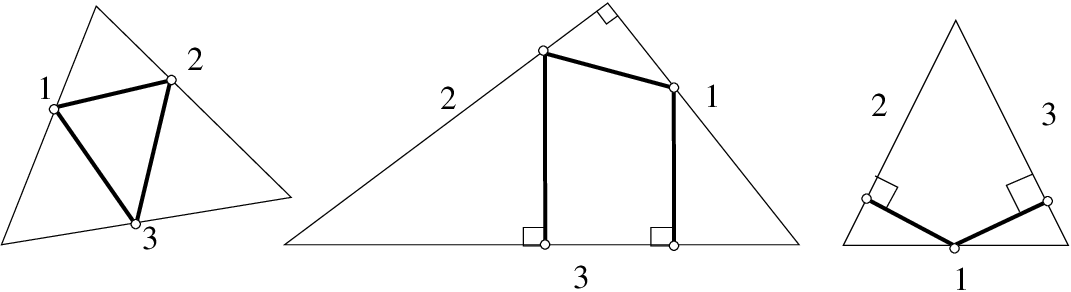}}
\vspace{-1em}
\end{center}
\begin{center}
Figure 1.1: Periodic paths having combinatorial type $123$ and $312321$
and $2131$.
\end{center}

One can show in elementary geometric ways
that acute, right, and isosceles triangles always
have periodic billiard paths.
The famous Fagnano path, which goes back to $1775$,
exists on a triangle if and only if the triangle is acute.
The Fagnano path has combinatorial type
$123$, meaning that it hits side $1$ of the triangle,
then side $2$, then side $3$, and then closes up.
Assuming the sides are appropriately labelled,
a periodic billiard path of
combinatorial type $312321$ exists on a triangle if and
only if the triangle is right and a periodic 
billiard path of combinatorial type
$2131$ exists on a triangle if and only if the triangle is
isosceles.  Examples of these paths are shown in Figure 1.1.

Periodic billiard paths in right triangles have
been fairly extensively studied.   The path
$312321$, mentioned above, is part of a
broader class of periodic billiard paths
in right triangles,
introduced in \cite{GSV}, which start out perpendicular to a
side of the triangle. One sometimes calls these
``perpendicular billiard paths''.
In \cite{CHK} it is shown that almost every perpendicular
billiard path is periodic.  In \cite{T05} this result
is refined, for right triangles with small angle
lying in $(\pi/6,\pi/4)$:  For such triangles,
all but one perpendicular billiard path is periodic,
and the union of these periodic perpendicular paths is
dense in the phase space.  In \cite{GZ} and
\cite{Vo05}, some classes of periodic billiard
paths in triangles are introduced and shown to be
{\it unstable\/} -- meaning that paths of the
same combinatorial type
do not exist on nearby acute or obtuse triangles.
Finally, in \cite{H2} it is shown
that all periodic billiard paths in right triangles
are unstable.

Every rational triangle -- i.e., a triangle whose
angles are rational multiples of $\pi$ -- has a periodic billiard path.
Indeed, any given rational polygon has
a dense set of periodic billiard
paths \cite{BGKT98}.  See also \cite{M} and \cite{Vo05}.
The subject of billiards on rational polygons
is a deep and extensive one.
For instance, see \cite{V} for connections between rational
billiards and Teichm\"uller Theory.  See
\cite{MT} and \cite{T} for surveys on
rational billiards.

Much less is known about periodic billiard paths in
obtuse, irrational triangles.  The paper
\cite{GSV}, the first to make serious inroads
into this question, produces some infinite families of
stable periodic billiard paths in obtuse
triangles. A periodic billiard path on a triangle
is called {\it stable\/} if all sufficiently 
nearby triangles have
a combinatorially identical periodic billiard path.
The paper \cite{HH} continues the program started
in \cite{GSV}, producing additional families of
stable periodic billiard paths on obtuse triangles.
The papers
\cite{H}, \cite{S1}, \cite{S2} exhibit additional 
infinite families of periodic billiard paths
for some obtuse triangles.
In particular, in \cite{S1},
\cite{S2} it is
shown that a triangle has a periodic billiard path
provided that all its angles are at most $100$
degrees.  

We already mentioned above that the
periodic billiard path $2131$ exists on
any isosceles triangle.   Unfortunately, this
path is unstable: It disappears as soon as we perturb the
triangle so that it is no longer isosceles.  One
might wonder if there is more to the story for
such perturbations.
Here is the main result of this paper.

\begin{theorem}
\label{isos}
Any sufficiently small perturbation of an isosceles triangle
has a periodic billiard path.
\end{theorem}

\subsection{Overview of the Proof}

The parameter space of (obtuse) triangles is an open
triangle $\Delta \subset \R^2$, where the
point $(x,y)$ corresponds to the obtuse triangle
whose small angles are $x$ and $y$ radians.
To each infinite periodic word $W$, with digits
in the set $\{1,2,3\}$, we assign the region
$O(W) \subset \Delta$ as follows:  A point
belongs to $O(W)$ if $W$ describes the
combinatorics of a periodic billiard path in
the corresponding triangle.  By this we mean that
we label the sides of the triangle $1$, $2$, and $3$,
and then read off $W$ as the sequence of successive
edges encountered by the billiard path.  We call
$O(W)$ an {\it orbit tile\/} and $W$
a {\it combinatorial type\/}.

The open line segment $\{x=y\} \subset \Delta$ parametrizes the obtuse
isosceles triangles.
We prove
Theorem \ref{isos} by covering a neighborhood of this segment
by orbit tiles.  This innocent-sounding
approach gives rise to an extremely intricate covering
involving infinitely many combinatorial types.
As we go along, we will try to explain why the complexity
seems necessary.

We define the
{\it Veech points\/} 
\begin{equation}
V_n=(\pi/2n,\pi/2n); \hskip 30 pt n=3,4,5...
\end{equation}  
These points are
special because they correspond to triangles
which have Veech's famous lattice property \cite{V}.
First of all, we will prove the following result.

\begin{theorem}
\label{nonveech}
A point on the obtuse isosceles line lies in the
interior of an orbit tile provided it is not of
the form $V_n$.
\end{theorem}
Theorem \ref{nonveech} involves a doubly-infinite family
of tiles -- with infinitely many tiles existing between
consecutive Veech points. See Figure 1.3 below. Experimentally, the family
in Theorem \ref{nonveech} seems to be the most efficient one
by far.
Theorem \ref{nonveech} focuses
our attention on the Veech points. 

We find it useful
to treat the points $V_n$ in a uniform way.  We 
decompose neighborhoods of the Veech points into quadrants.
Let $B(\epsilon)$ denote
the ball of radius $\epsilon$, centered at the origin.
Let $B_{\pm,\pm}(\epsilon)$ denote the
intersection of $B_{\epsilon}$ with the
open $(\pm,\pm)$ quadrant.
Now define
\begin{equation}
N_{\pm,\pm}(n,\epsilon)=V_n+B_{\pm,\pm}(\epsilon).
\end{equation}
Finally, let $N(n,\epsilon)$ denote the $\epsilon$
neighborhood of $V_n$.

\begin{theorem}
\label{easy}
For each $n\geq 4$ there are words
$A_n$, $B_n$, and $C_n$, and
some $\epsilon_n>0$ such that
\begin{itemize}
\item $N_{--}(n,\epsilon_n) \subset O(A_n)$.
\item $N_{-+}(n,\epsilon_n) \subset O(B_n)$
\item $N_{+-}(n,\epsilon_n) \subset O(C_n)$.
\item $B(\epsilon_n)-\overline N_{++} \subset
O(A_n) \cup O(B_n) \cup O(C_n)$
\end{itemize}
\end{theorem}
The last statement is present to take care of the boundaries of
the quadrants.
See Figures 1.2 and 1.4 below. 
It is worth remarking that the words $A_n$ are part of a larger family discovered in \cite{GSV}.  See also \cite{HH}.

Theorem \ref{easy} focuses our attention on the
regions $N_{++}(n,\epsilon)$.
Theorem \ref{easy} and Conjecture \ref{patveech} together
imply that these regions do not have
finite covers by orbit tiles, at least  when $n$ is a power of $2$.
We deal with all values of $n$ at once, by introducing
a doubly infinite
family $\{W_{nk}\}$ of words.  Here $n=4,5,6...$ and
$k=0,1,2...$.  
Figure 1.4 below shows some of the corresponding
orbit tiles for $n=4$.

\begin{theorem}
\label{cover1}
For each $n \geq 3$, there is some $\epsilon=\epsilon_n$
such that 
$$\overline N_{++}(n,\epsilon)-\{V_n\}\subset \bigcup_{k=0}^{\infty} O(W_{nk}).$$
\end{theorem}

Theorems \ref{nonveech}, \ref{easy}, and \ref{cover1} take care
of neighborhoods of all points  except $V_2$ and $V_3$.
The example we work out in the next chapter shows that
$V_3 \in O(W)$ for a certain word $W$ of length $22$.
See Corollary \ref{V3}.  Alternatively,
Theorem \ref{thm:stability} handles a neighborhood of $V_3$.  Finally,
in \cite{S2} we proved that a neighborhood of $V_2$, the
point corresponding to the right-angled isosceles triangle,
is contained in the union of $9$ orbit tiles. 
This completes the proof of Theorem \ref{isos}.

It is worth remarking that some of the complexity in
our proof seems necessary.  We will prove the following
result.

\begin{theorem}
\label{thm:stability}\hspace{1em}\vspace{-1em}\\
\begin{enumerate}
\item For $k=3,4,5...$ the triangle $V_{2^k}$  does not lie in the interior of an orbit tile.
\item For $n \geq 3$ and not a power of two, $V_n$ does lie in the interior of an orbit tile.
\end{enumerate}
\end{theorem}

In a separate and independent way, 
Statement 2 of Theorem \ref{thm:stability} handles
all the Veech points except $V_{2^k}$ for $k=3,4,5...$
However, this result does not really save us any time
in our analysis, because we still need to use
the analysis above to cover the remaining Veech points.
Statement 1 suggests
that perhaps the remaining points will be trouble.
Indeed, computer evidence strongly supports the following conjecture.

\begin{conjecture}
\label{patveech}
For $k=3,4,5...$ no neighborhood of $V_{2^k}$ has a
finite covering by orbit tiles. 
\end{conjecture}

\noindent
{\bf Remark:\/} 
In \cite{S1} it is proved that
no neighborhood of $(\pi/3,\pi/6)$ has a finite covering by
orbit tiles.
Thus, the covering constructed in \cite{S1} and \cite{S2}
for the ``100 degree result''mentioned above is necessarily infinite,
partly because of the ``trouble spot'' at the point
corresponding to the $(30,60,90)$-triangle.  In the
same way, Conjecture \ref{patveech} states that there are
an infinite number of ``trouble spots'' at various Veech points.
\newline

\subsection{Some Pictures of the Tiles}

We discovered all the results in this paper using our computer
program, McBilliards, a well-documented Java-based
program which is publicly available. \footnote{\url{http://mcbilliards.sourceforge.net}}  
 The reader
can see great pictures of our tiles using
McBilliards.  Here we reproduce a few of these pictures.

\begin{center}
\resizebox{!}{2in}{\includegraphics{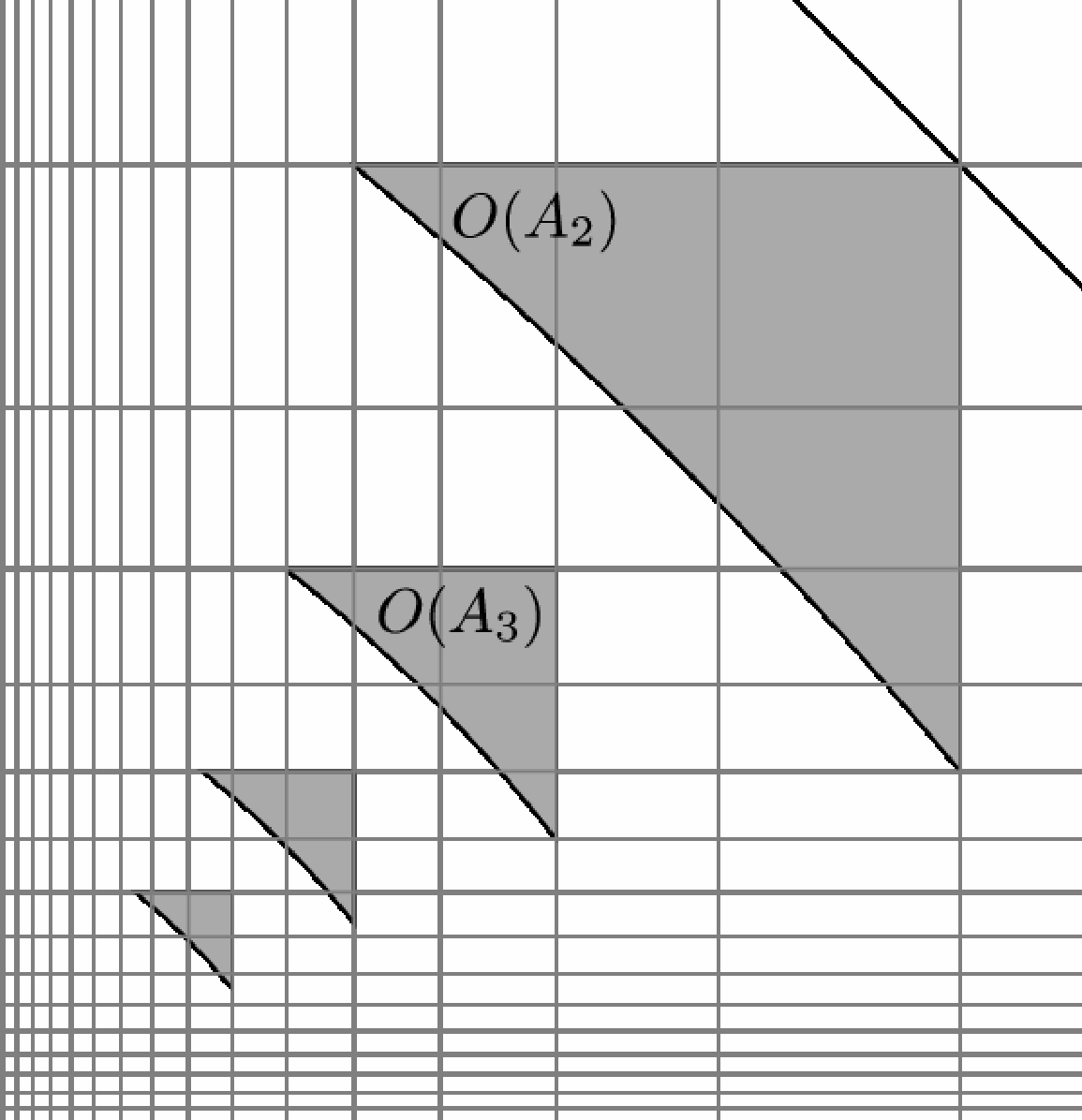}}
\vspace{-1em}
\end{center}
\begin{center}
Figure 1.2: The orbit tiles $O(A_{n})$ for $n=2,3,4,5$.
\end{center}

Figure 1.2 shows a picture 
the first few orbit tiles in the $A$ series.
The horizontal
grid lines have the form $x_2=\pi/n$ for $n=4,5,6...$
and the vertical grid lines have the form $x_1=\pi/n$ for $n=4,5,6$.
The right-angled tips of these tiles are the
Veech points.

The tiles in the $A$ series are also part of the tiles of we use
to prove Theorem \ref{nonveech}.   Theorem
\ref{nonveech} uses a double-infinite family $\{Y_{n,m}\}$ of
words with $m \in \N$ and $n \in \{2,3,4...\}$.  We have
$A_n=Y_{n,1}$.   For $n$ fixed, the tiles $\{Y_{n,m}\}$ live
between the two consecutive Veech points $V_n$ and $V_{n+1}$.
Figure 1.3 shows some of these tiles.

\begin{center}
\includegraphics[width=6.5in]{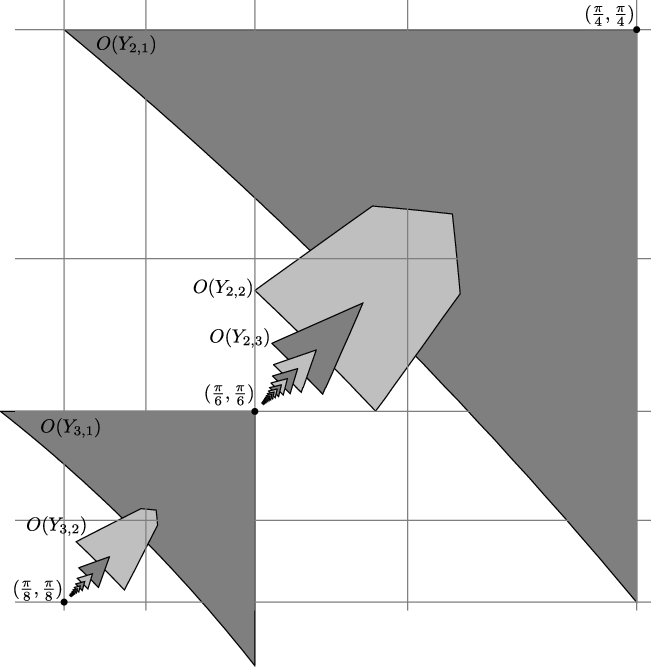}
\newline
Figure 1.3: Some of the tiles $O(Y_{n,m})$ with $n \in \{2,3\}$ and $m \in \{1,2,3,\ldots\}$. 
\end{center}

Figure 1.4 shows a neighborhood
of the point $V_4=(\pi/8,\pi/8)$, a point
in the bottom left corner of Figure 1.3.
The tiles $O(B_4)$ and $O(C_4)$ are tiny in comparison to
the size of $O(A_4)$, so much of $O(A_4)$ is off the
screen.  The union covers a neighborhood of $V_4$.

\begin{center}
\includegraphics[height=5in]{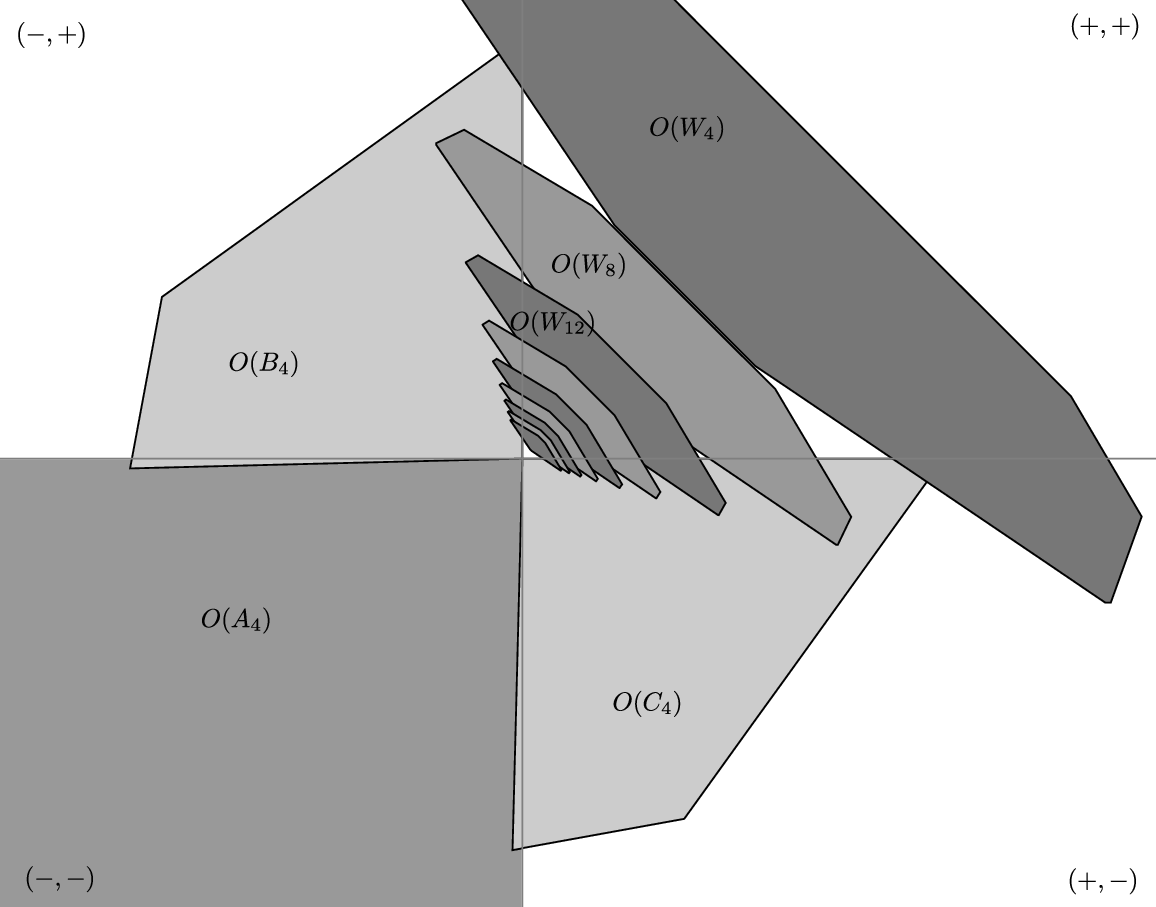}
\newline
Figure 1.4: The orbit tiles $O(A_4)$ and $O(B_4)$ and $O(C_4)$
and $O(W_{4k})$ for $k=1,...,9$.
\end{center}

\subsection{Asymptotic Self-Similarity and Fourier Series}

Let $O_{nk}=O(W_{nk})$.
The following self-similarity result, which is
the central technical result in the paper, implies
Theorem \ref{cover1}.

\begin{theorem}[Central Lemma]
\label{limit1}
Let $S_{nk}$ be the dilation which maps $V_n$ to $0$ and
expands distances by
$$
\zeta_n k^2; \hskip 30 pt \zeta_n:=2(n-1)\cot(\pi/2n) \approx 4 n^2/\pi
$$
If $n$ is held fixed and $k \to \infty$ then
the closure of $S_{nk}(O_{nk})$ Hausdorff-converges 
to the convex quadrilateral $Q_n$ with vertices
$$v_1=(-\frac{1}{n},1-\frac{1}{n}); \hskip 15 pt
v_2=(1-\frac{1}{n},-\frac{1}{n}); \hskip 15 pt
v_3=(a_n,a_n); \hskip 15 pt v_4=(\mu_n a_n,\mu_n a_n);$$
where
$$a_n=\frac{1}{2}-\frac{1}{2n}; \hskip 15 pt
\mu_n=\frac{1}{2}-\frac{\tan^2(\pi/2n)}{2}.$$
The convergence is such that
any compact subset $Q'_n \subset Q_n$ is contained
in $S_n(O_{nk})$ for $k$ sufficiently large in comparison
to $n$.
\end{theorem}

An example of the limiting quadrilateral $Q_n$ is shown in figure 1.5, below. 

\noindent
{\bf Proof of Theorem \ref{cover1}:\/}
Fixing $n$, Theorem \ref{limit1}
implies that there are constants $0<\epsilon_1<\epsilon_2$ such that,
for $k$ sufficiently large,
$O(W_{nk})$ contains the set $V_n+\Lambda_k$, where
$\Lambda_k$ is the convex hull of
$$
(\epsilon_1/k^2,0); \hskip 15 pt
(\epsilon_2/k^2,0); \hskip 15 pt
(0,\epsilon_1/k^2); \hskip 15 pt
(0,\epsilon_2/k^2)$$
for $k$ large.
But the union of the sets $V_n+\Lambda_k$ covers $N(n,\epsilon)$ for some
$\epsilon>0$.  These sets ``bunch up''
as $k \to \infty$.  Compare Figure 1.4.  So, Theorem \ref{limit1} 
proves Theorem \ref{cover1}.
\endproof 

\begin{center}
\includegraphics[height=2in]{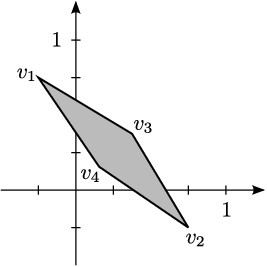}
\vspace{-1em}
\end{center}
\begin{center}
Figure 1.5: The limiting quadrilateral $Q_4$.
\end{center}

Our technique for proving Theorem \ref{limit1}
involves looking at the Fourier transforms of the analytic
functions which define the edges of the orbit tiles
of interest to us.  The Fourier transforms of these
functions are functions defined on $\Z^2$.
It turns out that the supports of these Fourier
transforms grow ``linearly'' with the parameter $k$
in a way we make precise in \S \ref{sect:rescaling}.  To deal with
with the situation as $k \to \infty$
we prove the Quadratic Rescaling Theorem, a
result which describes the asymptotic limits
of a family of functions which vary with the
prescribed growth.   One of the main technical
innovations in the paper is a combinatorial
method for understanding such growing families of
Fourier series.
Our technique seems to be more general than the
application we give here, but so far this is the
main application.

\subsection{Paper Outline}

In \S \ref{sect:background} we will give background information about
triangular billiards, and in particular discuss how
one computes the orbit tile $O(W)$ based on the
combinatorics of the word $W$.    All of the constructions in
\S \ref{sect:background} are programmed into McBilliards.
The interested reader can see these constructions
in action when using the program.

In section \ref{sect:non-veech}, we prove theorem \ref{nonveech}.

In \S \ref{sect:easy} we will prove Theorem \ref{easy}.

\S \ref{s3}-\ref{s6} are devoted to the proof of
the Central Lemma, namely Theorem \ref{limit1}. 
In \S \ref{s3} we will introduce our words $W_{nk}$,
and prove some preliminary results about the
orbit tiles $O_{nk}:=O(W_{nk}).$ In particular,
we will isolate a region $R_{nk} \subset \Delta$
with the property that (independent of $n$ and $k$)
a certain $16$ functions define $O_{nk} \cap R_{nk}$.

In \S \ref{sect:rescaling} we will prove the Quadratic Rescaling Theorem,
a result that is designed to analyze infinite families of
defining functions, such as the $16$ families we isolate in \S  \ref{s3}.
In \S \ref{sect:pivot} and \S \ref{s6} we will use the
Quadratic Rescaling Theorem to finish the proof of
Theorem \ref{limit1}.

In \S \ref{sect:stability}, which is logically independent
from the rest of the paper,
we prove theorem \ref{thm:stability}. This classifies the Veech triangles $V_n$ which lie in the interior of an orbit tile.

\newpage

%% file: 2back.tex
\section{Billiard Paths and Defining Functions}
\label{sect:background}

\subsection{Unfoldings}

The {\it unfolding\/} of a word $W$ with respect to a triangle $T$,
which we denote by $U(W,T)$, is the union of triangles
obtained by reflecting $T$ out according to the digits of
$W$.  This construction is discussed in detail in
\cite{S1} and in e.g. \cite{T}.
We will persistently abuse our notation in the following sense:
A point $X$ in parameter space represents a triangle $T=T_X$.
We will often write $U(W,X)$ in place of $U(W,T)$.

There is a
sequence of vertices which runs across the top of
$U(W,T)$.  We call these the {\it top vertices\/} and
label them $a_1,a_2,...$ from left to right.  There
is a sequence of vertices which runs across the bottom
of $U(W,T)$ and we label these $b_1,b_2,...$ from
left to right.  Figure 2.1 shows an example of
an unfolding, with respect to the Veech point $V_3$.

\begin{center}
\includegraphics[height=1.3in]{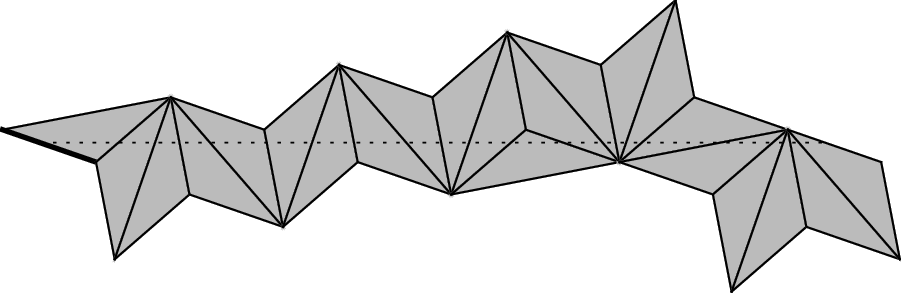}
\newline
Figure 2.1: $U(2323132313123232313131,V_3)$ with a centerline
\end{center}

It is worth pointing out that one of the apparent edges
of the unfolding in Figure 2.1 is not actually an
edge of reflection.  Figure 2.2 shows the unfolding
of the same word with respect to a different point.

\begin{center}
\includegraphics[height=1.3in]{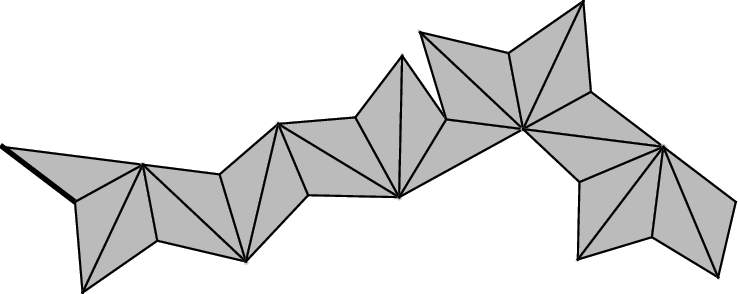}
\newline
Figure 2.2 $U(2323132313123232313131,(\pi/5,\pi/6))$
\end{center}

\noindent
{\bf Remark:\/}
We mainly care about unfoldings of isosceles triangles
or triangles that are very nearly isosceles.  From
here on in, most of our pictures show unfoldings of
isosceles triangles.  However, as with Figure 2.2, 
we sometimes show an unfolding of a non-isosceles
triangle so as to illustrate a point.
\newline

The first side of $U$ has been highlighted in both examples.
$W$ represents
a periodic billiard path in $T$ iff
the first and last sides of $U(W,T)$
are parallel and the interior of
$U(W,T)$ contains a line segment $L$,
called a {\it centerline\/},
such that $L$ intersects the first and
last sides at corresponding points.
In both examples above, the first and last
sides are parallel.  However, the
centerline only exists for Figure 2.1.
In particular, Figure 2.1 shows that the given word
describes a periodic billiard path for the
triangle corresponding to $V_3$.
As in Figures 2.1 and 2.2
we always rotate the picture so that
the first and last sides are related by
a horizontal translation.  We call this
horizontal translation the {\it holonomy\/}.

\subsection{Stability and Hexpaths}

A word $W$ is called {\it stable\/} if the first
and last sides of $U(W,T)$ are parallel for any
triangle $T$.  This implies that $O(W)$ is an
open set. 
In this section we 
will explain a combinatorial
criterion for stability.  The proof is well known,
and we omit it.  See \cite{S1} for details.

\begin{center}
\resizebox{!}{2.9in}{\includegraphics{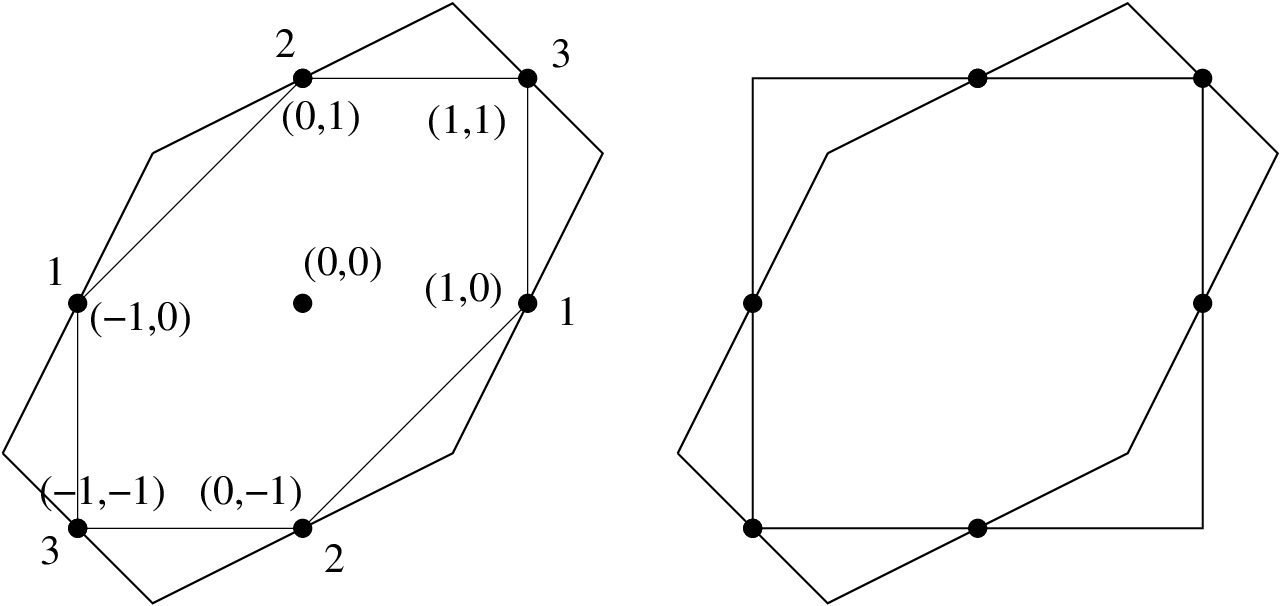}}
\newline
Figure 2.3: The fundamental hexagon
\label{fig2.3}
\end{center}

Let $H_0$ be the outer hexagon shown in Figure 2.3.
The shape if $H_0$ is a bit strange, but the
inscribed hexagon has vertices on the integer
lattice $\Z^2$ as shown.   Also, $H$ is well related
to a square of side-length $2$, as shown on the
right hand side of Figure 2.3.
The sides of $H_0$
are divided into $3$ types, according to their
label. Let $\cal H$ denote the tiling of
$\R^2$ by translates of $H_0$.  By $\cal H$
we really mean the union of edges of the tiling.
By construction, the midpoints of edges in $\cal H$
lie in $\Z^2$.

Given the word $W$,
we can draw a path in $\cal H$ by following the edges
as determined by the word: we move along the
$d$th family when we encounter the digit $d$.
Figure 2.4 shows the path corresponding
to the examples given in Figures 2.1 and 2.2.
The dot
in the picture 
indicates the start of the path.  We call this path
the {\it hexpath\/} and denote it by $H(W)$.

\begin{lemma}[Hexpath]
\label{closedhexpath}
The word $W$ is stable iff $H(W)$ is a closed path.
\end{lemma}

This condition in the Hexpath Lemma
 is equivalent to the better known condition, which appears
as lemma 3.3.1 in \cite{T}. We have restricted this lemma to our context.

\begin{lemma}
\label{lem:stability}
A word $W$ is stable iff the number of times each letter $\ell=1,2,3$ appears
in an odd position in $W$ equals the number of times $\ell$ appears in an even position.
\end{lemma}

This condition can easily be verified for the example we have been considering.
In addition it happens for squares of words of odd length.

\begin{corollary}[Odd squares are stable]
\label{cor:odd_stable}
If $W$ is a word of odd length, then $W^2$ is stable.
\end{corollary}

Now we can take care of the loose end from the introduction.

\begin{corollary}
\label{V3}
$V_3$ is contained in the interior of an orbit tile.
\end{corollary}

\startproof
Figure 2.1 shows that the given word describes
a periodic billiard path for the triangle
corresponding to $V_3$.  Figure 2.4 below shows
that the hexpath corresponding to this word is closed.
Hence, the corresponding billiard path is stable.
\endproof

\begin{center}
\includegraphics{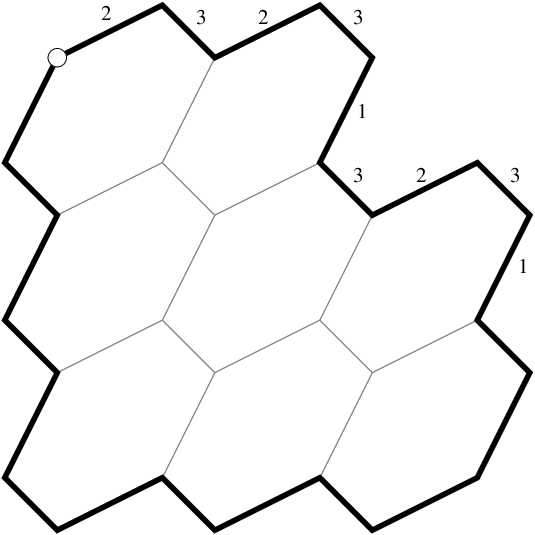}
\newline
Figure 2.4: The hexpath for $W=2323132313123232313131$.
\end{center}

\subsection{The Squarepath}

It turns out that the hexpath $H(W)$ contains precisely
the same information as a certain rectilinear path,
which we call the {\it squarepath\/}.  Each vertex of the
hexpath has a unique type $3$ edge emanating from it.
The squarepath is obtained by connecting the
midpoints of these type-3 edges together, in order.
We denote the squarepath by $\widehat Q(W)$.
We can also define similar paths based on the
edges of type $1$ or $2$.  These paths are
somewhat more complicated, though they will be
of theoretical importance for us.  In practice,
however, we will always try to work with the
type $3$ edges.

If we mark off points on the squarepath at
integer steps (starting with a vertex) the resulting
points are naturally in bijection with the type
$3$ edges of the unfolding.  In the next section
we will elaborate on this bijection.
Figure 2.5 shows the squarepath for the examples we have
been considering
The hexpath is drawn underneath in grey.

\begin{center}
\includegraphics{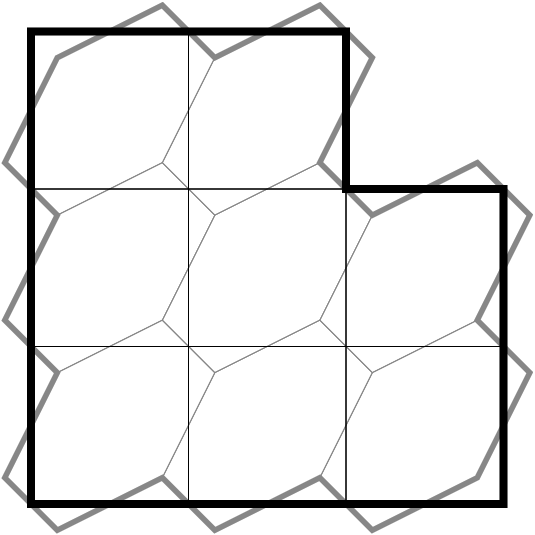}
\newline
Figure 2.5: $\widehat Q_3(W)$ in black and $H(W)$ in grey.
\end{center}

It is possible to reconstruct $H(W)$ from $\widehat Q(Q)$,
when $\widehat Q(W)$ is a closed loop.
When $\widehat Q(W)$ is embedded, this loop bounds a finite
union of squares. We simply replace each square by the
associated hexagon. Then $H(W)$ is the boundary of the
union of hexagons.  In general, $H(W)$ is the union of
all the edges of $\cal H$ which intersect $\widehat Q(W)$.
There is a natural ordering to these edges, and so
the union of all these edges naturally has the structure
of a closed loop.

It turns out that there is a simple algorithm for deducing the
combinatorics of the unfolding from the squarepath.   Say that a
$k$-{\it dart\/} is a union of $k$ triangles,
arranged around a common vertex, in the pattern shown in
Figure 2.6 for $k=2,3,4.$.  A $k$-dart is just an unfolding
with respect to either the word $(13)^{k-1}1$ or the word
$(23)^{k-1}2$. 
\newline
\newline
{\bf Remark:\/}  We shall almost always consider darts
made from isosceles triangles. Indeed, the idea of grouping
the unfolding into darts is mainly a combinatorial trick,
and in our applications we might as well perform the
trick with respect to unfoldings of isosceles triangles.
However, some of our pictures
show darts made from triangles that are not quite isosceles.

\begin{center}
\includegraphics{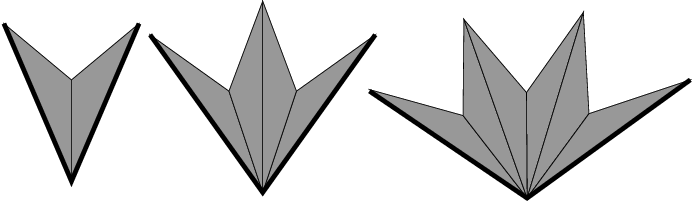}
\newline
Figure 2.6: $k$-darts for $k=1,2,3$
\end{center}

We say that the {\it $3$-spine\/} of the dart is the union of
the two outermost long edges.  We have highlighted
the spines of our darts in Figure 2.6.

The relation of $U(W,*)$ to $\widehat Q(W)$ is as follows:
\begin{itemize}
\item The maximal darts of the unfolding are in bijection
with the edges of the square path.  (The maximal
$k$-darts correspond to edges of length $2k$.)
The maximal darts are glued together along their $3$-spines. 
\item Two consecutive maximal darts lie on opposite
sides of their common $3$-edge iff $\widehat Q(W)$ makes a
northwest or southeast turn at the vertex
corresponding to this $3$-edge.
\end{itemize}
To make this work precisely, we need to take the
infinite periodic continuation of $U$, or else
identify the first and last sides of $U$ to make
an annulus.  As it is, the reader needs to take
special care in figuring out how the rightmost
maximal dart fits together with the leftmost one.
We have included a copy of Figure 2.2, except with
the spines of the maximal darts drawn in black. See Figure 2.7.

\begin{center}
\resizebox{!}{1.5in}{\includegraphics{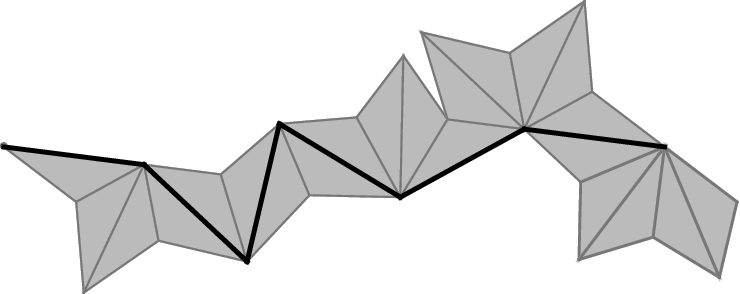}}
\newline
Figure 2.7: Dividing the unfolding into maximal darts.
\end{center}

We have taken a lot of trouble to describe the squarepath and its
relation to the hexpath and the unfolding because we plan to
specify all our words in terms of their squarepaths. Using
the square path gives a very simple description of the
word, and lets the reader best see the patterns which arise
in our families.

\subsection{Edge Labellings}
\label{s2.4}

We label each edge of $\cal H$ by the coordinates of its midpoint.
This labelling is canonical, once we decide which point
of $\Z^2$ gets labeled $(0,0)$.  
The {\it McBilliards convention\/} is to assign the
label $(0,0)$ to the edge of $H(W)$ corresponding to
the last digit of $W$.  This edge is the leftmost edge
of the unfolding $U(W,*)$.  In Figure 2.8 we have
labeled the origin and several nearby points.

\begin{center}
\includegraphics{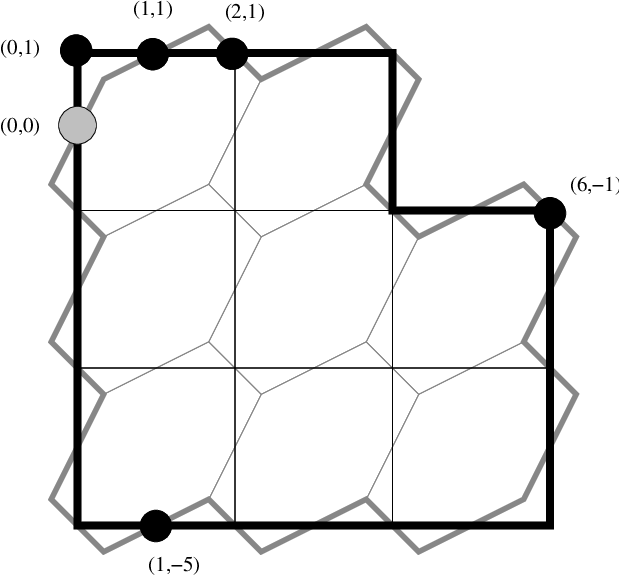}
\newline
Figure 2.8: Some labellings.
\end{center}

We identify $\widehat e$ with its label.
Our labelling has a geometric significance.
Let $X=(x_1,x_2)$ be a parameter point and
let $T_X$ be the corresponding triangle.
Let $e_1$ and $e_2$ be two edges of the unfolding $U(W,T_X)$.
Let $\theta(e_1,e_2)$ be the counterclockwise angle
through which we must rotate $e_1$ so as 
to produce an edge parallel to $e_2$.
We take $\theta$ mod $\pi$, so that the orientations
of $e_1$ and $e_2$ are irrelevant.
Then, as is easily established by induction:
\begin{equation}
\theta(e_1,e_2)=X \cdot (\widehat e_2-\widehat e_1).
\end{equation}

\subsection{Defining Functions}  

We frequently write
\begin{equation}
E(x)=\exp(i(x))
\end{equation}
for notational convenience.

Given two points $p,q \in \R^2$ we write
$$p \uparrow q; \hskip 20 pt 
p \updownarrow q; \hskip 20 pt
p \downarrow q$$
iff the $y$ coordinate respectively is
greater than, equal, or less than the $y$ coordinate of $q$.
Suppose that $p$ and $q$ are two vertices of our unfolding.
In this section we will give the formula
for a function $F=F_{p,q}$ which has the property
that $F=0$ iff $p \updownarrow q$. 
These {\it defining functions\/} are computed purely
from the word $W$.  The orbit tile $O(W)$ can be
described as the region where the defining functions
corresponding to the $(a_i,b_j)$ pairs are all positive.
The edges of $O(W)$ is defined in terms of the $0$-level
sets of the defining functions.

For any $d \in \{1,2,3\}$ there is an infinite,
periodic polygonal path made from type-$d$ edges
of the infinite periodic continuation of
$U(W,T)$.  The image of 
this path in $U(W,T)$ is what we call the
$d$-{\it spine\/}.   We have already encountered
the $3$-spine:  It is the union of the
$3$-spines of the maximal darts of $U(W,*)$.

\begin{center}
\includegraphics{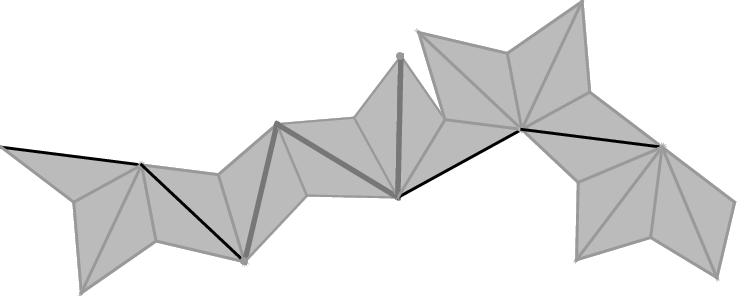}
\newline
Figure 2.9: A $3$-path and the $3$-spine.
\end{center}

Let $e_0,...,e_m$ denote the list of
edges, ordered from left to right, which
appear in the $d$-spine.
We say that the vertices $p$ and $q$ are
$d$-{\it connected\/} if
there is a polygonal path of type-$d$ edges
connecting $p$ to $q$.  In this case,
let $f_0,...,f_n$ denote these edges,
ordered from left to right.  We order $p$ and $q$ so that
$p$ is the left endpoint of the $d$-path and
$q$ is the right endpoint.  The $3$ thick grey edges
in Figure 2.9 show the $3$-path connecting
$p=b_4$ to $q=a_6$.  

We define
\begin{equation}
\label{alter}
P(X)=\pm \sum_{i=0}^n (-1)^i E(X \cdot \widehat f_i);
\hskip 30 pt
Q(X)=\sum_{i=0}^m (-1)^i E(X \cdot \widehat e_i).
\end{equation}
We will explain the global sign in front of $P$ below.
The reason for
the general alternation of the signs is explained in \cite{S1}.
Our functions have the following
geometric interpretation: If we normalize
so that the $d$ edges have length $1$ and rotate
$U(W,T)$ so that the first edge is horizontal,
then $\pm P(X)$ is the vector pointing from $p$ to $q$ and
$Q(X)$ is the translation vector. 
Therefore, 
\begin{equation}
\label{form1}
F:= {\rm Im\/}( \pm P \overline Q)=0 \hskip 15 pt
\Longleftrightarrow \hskip 15 pt p \updownarrow q.
\end{equation}

For the above example the sign in front of $P$ turns out to be a $(+)$.
(See below.)
We therefore have
$$P(X)=E(4x_1-x_2)-E(6x_1-x_2)+E(6x_1-3x_2).$$
$$Q(X)=E(x_2)\!-\!E(4x_1+x_2)+E(4x_1\!-\!x_2)\!-\!E(6x_1\!-\!x_2)+E(6x_1\!-\!5x_2)+E(\!-\!5x_2).$$
Here is what we call the {\it function tableau\/} for $P$.
$$
\begin{matrix}
(+)&4&-1 \cr
&6&-1 \cr
&6 & -3\end{matrix}
$$
When we reconstruct the function from its tableau, we
use the convention that the signs of the terms alternate. 
The $(+)$ of $(-)$ indicates the global sign in front of $P$.

\begin{center}
\resizebox{!}{4in}{\includegraphics{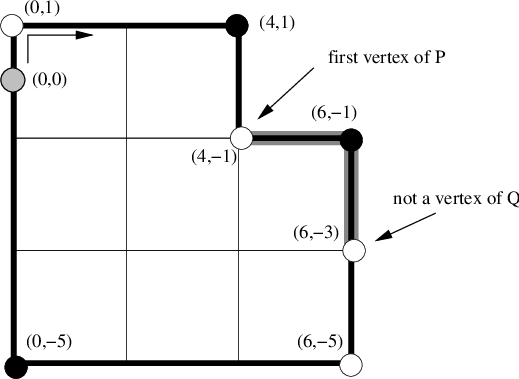}}
\newline
Figure 2.10 The paths $\widehat P$ (grey) and $\widehat Q$ (black)
\label{fig2.9}
\end{center}

When $d=3$, we can
represent both $P$ and $Q$ in terms of the squarepath.
First of all, the list of vertices of $\widehat Q$ is
precisely the function tableau for $Q$.  The situation
for $P$ is more involved:  The edges of
$U(W,*)$ are in canonical bijection with the
edges which emanate from the vertices of the hexpath.
Say that a $3$-edge of
$\cal Q$ is a {\it starter\/} if it corresponds to an
edge of $U(W,*)$ which is incident to $p$.
Say that a $3$-edge of 
$\cal H$ is a {\it finisher\/} if it corresponds to an
edge of $U(W,*)$ which is incident to $q$. 
Let $\widehat P$ denote the shortest
sub-path of $\widehat Q_d(W)$
whose initial endpoint is the midpoint of a
starter and whose final
endpoint is the midpoint of a finisher. 
Then the function tableau for $P$ is just the
list of coordinates of the vertices of $\widehat P$.
Figure 2.10 shows the paths corresponding
to $\widehat P$ and $\widehat Q$.  The origin is
marked with a grey dot.

We can interpret the path $\widehat Q$ as a function
from $\Z^2$ to $\Z$, as follows.  We alternately
color the vertices encountered by $Q$ black and white,
starting with white.
$\widehat Q$ assigns the value
$x_1-x_2$ to $X \in \Z^2$ if $x_1$ white
vertices of $\widehat Q$ coincide with $X$ and
if $x_2$ black vertices of $\widehat Q$ coincide with $X$.
We make the same definition for $\widehat P$, except
that we have to take care whether or not to color
the first vertex encountered by $\widehat P$ black or
white.  (See below.)
With this interpretation,
$\widehat Q$ is the Fourier series of $Q$.
\begin{equation}
Q(X)=\sum_{V \in \Z^2} \widehat Q(X) E(X \cdot V).
\end{equation}
The same goes for $P$ and $\widehat P$.
\newline
\newline
{\bf The Global Sign:\/}
This discussion supposes that $F>0$ if $q \uparrow p$.
(As above, $p$ is on the left.)
We also suppose that the
initial vertex of $\widehat P$ is also a vertex of
$\widehat Q$.  In this case, the sign in front of $P$ is $(-1)^u$, where $u$
is the number of vertices of $\widehat Q$ (starting with
the first one) which lie before the first vertex of $\widehat P$.
That is, the initial vertex of $\widehat P$ should get the
same color whether it is considered a vertex of $\widehat P$
or a vertex of $\widehat Q$.
For example, we can see from Figure 2.10 that $u=2$ and so
the sign is a $(+)$.
This rule has a simple geometric proof:  When $p$ and $q$ are the
first and last vertices of the $3$-spine of $U$, then $P=Q$
and so the sign definitely should be a $(+)$.  If we 
move $q$ along the $3$-spine, the sign does not change,
by ``continuity'':  Moving either vertex by $1$ ``click'' should
produce a nearby value for $P$.  However, moving the
$p$ vertex changes the global sign, given the form
of Equation \ref{alter}.
In general the first vertex of
$\widehat P$ need not be a vertex of $\widehat Q$.
This irritating situation does not
arise in this paper. McBilliards has a general algorithm
which correctly determines the sign in every possible case.

\subsection{The Dart Lemma}
\label{sect:dart}

Figure 2.11 shows a typical dart.
We say that the {\it inferior\/}
vertices of $D$ are the ones which are not adjacent
to the $3$-spine and not on the $3$-spine.  The
inferior vertices are marked with white dots in
Figure 2.11.  We call the other vertices of the
dart {\it superior\/}. In Figure 2.11 the superior
vertices are in black or grey and the inferior
vertices are in white.

\begin{center}
\includegraphics{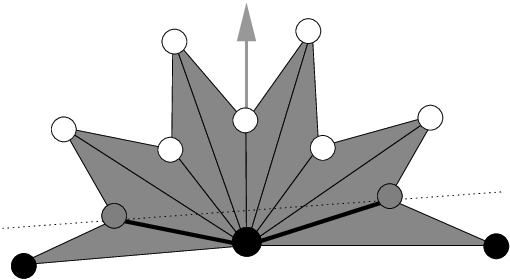}
\newline
Figure 2.11: An Acute Dart
\end{center}

Recall that the unfolding $U$ can be written as a union
of maximal darts.  We say that a vertex of $U$ is
{\it inferior\/} if it is an inferior vertex of
one of the maximal darts.
Let $\delta(W)$ denote the largest $k$ such that
$U(W,*)$ contains a $k$-dart.  Here is the main
result of this section:

\begin{lemma}[Dart]
\label{lem:dart}
Let $X=(x_1,x_2)$.  Suppose that
$$\max(x_1,x_2) \leq \frac{2\pi}{\delta(W)}.$$
Suppose also that all the
top superior vertices of $U(W,X)$ lie above all the bottom
superior vertices of $U(W,X)$.  Then $X \in O(W)$.
\end{lemma}

\noindent
{\bf Remark:\/} There is a more restrictive angle condition that almost immediately
guarantees that the maximal darts are acute. This condition is given by
$$\max(x_1,x_2) \leq \frac{\pi}{2 \delta(W)-2}.$$
This is precisely the condition we will use in the proof of
Theorem \ref{nonveech}.   Our condition is weaker than this
and does not, in itself, guarantee that the maximal darts are
acute.  However, our weaker condition combines with the second
hypothesis of the Dart Lemma to establish the acuteness.  See
the very end of our proof.  We mention this because the
discrepancy is likely to otherwise cause confusion.
\newline

We will prove the Dart Lemma in several stages.
We say that the {\it base\/} of a dart is the vertex which
is common to all the triangles.  The base is denoted
by a big black dot in Figure 2.11.
We say that the
{\it centerline\/} of the dart is the ray of bilateral
symmetry, emanating from the base.  The centerline is
indicated by a ray in Figure 2.11.  We say that the
dart {\it points up\/} if the ray points upward, and
{\it points down\/} if the ray points downward.  Let
$D_V$ denote the union
of outermost edges of $D$ which are not the longest edges.
This set is highlighted
in Figure 2.11.  We say that $D$
is {\it acute\/} if $D_V$ makes an acute
angle towards the centerline of the dart.
Figure 2.11 shows an acute dart.

\begin{lemma}
\label{dart}
If $D$ is an up-pointing acute dart,
then each inferior vertex of $D$ lies above some superior vertex.
Likewise, if $D$ is down-pointing and acute, then each
inferior vertex of $D$ lies above some superior vertex
of $D$.
\end{lemma}

\startproof
The short edges of $D$ have the same length.
Hence the line joining the
two superior vertices separates all the inferior vertices
from the base.
\endproof

We call $U$ {\it controlled\/} if the
following holds for all maximal darts $D$
of $U$:
\begin{itemize}
\item $D$ is acute.
\item If the base of $D$ is a bottom vertex of $U$ then $U$ points up
\item If the base of $D$ is a top vertex of $U$ then $U$ points down.
\end{itemize}

\begin{lemma}
\label{predart}
Suppose $U$ is a controlled unfolding.
Then the lowest top vertex of $U$ and
the highest bottom vertex of $U$ are both superior vertices.
\end{lemma}

\startproof
We will prove this statement for the top vertices.  The proof
for the bottom vertices is the same.  Let
$v$ be an inferior top vertex of $U$. Then there is some
maximal dart $D$ such that $v$ is an inferior vertex of $D$.
Each edge of $D$, except possibly the edges on the $3$-spine,
is an edge of reflection of $U$.  Thus, the inferior vertices
of $D$ all have the opposite type (top or bottom) from the
base.  Likewise for the superior vertices of $D$.  Hence,
the inferior vertices and the superior vertices of
$D$ all have the same type.  Since
one inferior vertex of $D$ is a top vertex, the base
of $D$ is a bottom vertex.  Since $U$ is controlled,
$D$ points up.  Lemma \ref{dart}
now implies that $v$ is higher than
one of the superior vertices $v'$ of $D$. As we already
mentioned, $v'$ is also a top vertex of $U$.  Hence,
we have found another top vertex, $v'$, which is lower
than $v$.
\endproof

To finish the proof of the Dart Lemma, we just have
to establish that $U=U(W,X)$ is a controlled unfolding.
Let $D$ be a maximal dart of $U$.  Assume without loss of
generality that the basepoint of $D$ is a bottom
vertex.  Each edge of
$D_V$ is an edge of reflection of $U$.  Hence
the endpoints of $D_V$ are top vertices.
All these vertices are superior vertices.
Hence, the endpoints of $D_V$ lie above
the basepoint of $D_V$.  Our restriction
on $X$ guarantees that the angle of $D_V$
is at most $2 \pi$.  Hence $D_V$ must actually make
an acute angle, since both its edges point up.
The centerline lies between these two
up-pointing edges.  Hence $D$ itself is
up-pointing.   Since $D$ is arbitrary,
we see that $U$ is controlled.  This completes
the proof of the Dart Lemma.

\subsection{Pseudo-Parallel Families}
\label{pp}

Suppose that $e$ is an edge of the unfolding
$U(W,*)$.
When $U(W,X)$ is rotated so that
it has horizontal holonomy, the line containing
$e$ is parallel to the complex number
\begin{equation}
E(\widehat e \cdot X) U_Q(X).
\end{equation}

We say that the edges $\{e_0,...,e_n\}$ form a
{\it pseudo-parallel family\/} relative to
the point $X_0$ if the dot product
$e_j \cdot X_0$ is independent of $j$. 
In this case, the edges $e_0,...,e_n$ are all
parallel in $U(W,X_0)$.  We assume that these
edges have negative slope in $U(W,X_0)$.
The points $\widehat e_0,...,\widehat e_n$ must lie
on a line segment in $\R^2$.  In our examples in
this paper, the line in question always has slope $-1$
because $X_0$ lies on the isosceles line.
We order our edges so that $\widehat e_0,...,\widehat e_n$
appear in order on the line.

Let $R'(e_j)$ denote the region in parameter space such
that $e_j$ has negative slope in the unfolding.  Let
$R(e_j)$ denote the path connected component of
$R'(e_j)$ which contains $X_0$. 

\begin{lemma}[Convex Hull]
$R(e_0) \cap R(e_n) \subset R(e_j)$
for all $j$.
\end{lemma}

\startproof
We think of $\{X_t\}$ as a path in
$R(e_0) \cap R(e_n)$ which connects
$X_0$ to some other point $X_1$.
Let $S^1$ denote the unit complex numbers and
let $E: \R \to S^1$ be the universal covering map.
For each object $z \in S^1$ we let
$\widetilde z$ denote the lift to $\R$, so that
$E(\widetilde z)=z$.  In particular, we define
$$U(t)=U_Q(X_t); \hskip 30 pt 
\widetilde E(t)=\widehat e_j \cdot X_t;
\hskip 30 pt
E_j(t)=E(\widehat x_j \cdot X_t).$$
Let $\widetilde I_t \subset \R$ be the interval whose
endpoints are $E_0(t)$ and $E_n(t)$.
By convexity $\widetilde E_j(t) \subset \widetilde I_t$ for all $t$.
The edges $e_0(t)$ and $e_n(t)$ have negative slope
for all $t$.  Hence
$\widetilde I_t$ has length less than $\pi/2$ for any
$t \in [0,1]$.   Hence
$E_t(j)$ lies in the arc $I_t$, which has length less than
$\pi/2$.
If we rotate $S^1$ so that $U(t)=1$ then 
the endpoints of $I_t$, namely $E_0(t)$ and $E_j(t)$,
are both contained in the same negative quadrant of
$\R^2$. (Either $(-+)$ or $(+-)$.)  Hence
$I_t$ is contained in one of the negative quadrants.
Hence $E_j(t)$ is also contained in one of
these quadrants.  That is, $e_j$ has negative
slope for any parameter value $t$.
\endproof

\newpage

%% file: 3pat.tex
\section{Proof of Theorem \ref{nonveech}}
\label{sect:non-veech}

\def\u{W} 

Our proof of Theorem
\ref{nonveech} 
is based on the $2$ parameter family $Y_{m,n}$ of
odd-length words.
\begin{equation}
\label{eq:Y}
Y_{n,m}=1 (\u_n)^m 3 2; \hskip 50 pt
\u_n = (31)^{n-1}  (32)^{n-1}
\end{equation} 
Figure 1.4 of the introduction shows the orbit tiles for some of the words $Y_{n,m}$.
The family $\u_n$, which we call the {\it unstable family\/},
describes unstable periodic billiard paths in certain
isosceles triangles of interest.   The square words
$Y_{m,n}^2$ are stable by  Corollary
\ref{cor:odd_stable}.
Theorem \ref{nonveech} is therefore a consequence of the
following result.

\begin{theorem}
\label{thm:Ycovers}
For every integer $n \ge 2$ and real number $x$ so that $\frac{\pi}{2n+2} < x < \frac{\pi}{2n}$,
there is a periodic billiard path in $T_x$ with combinatorial type $Y_{n,m}$ for some $m \in \N$.
\end{theorem}

Here $T_{x}$ denote the obtuse isosceles triangle corresponding to the point $(x,x)$ in the plane.  The two small angles of $T_x$ have measure $x$-radians.

\subsection{The Unstable Family}

\begin{proposition}
\label{prop:unstable}
$\u_n$ describes a periodic billiard path in $T_x$ for all
$x<\frac{\pi}{2n-2}$.
\end{proposition}

\startproof
The unfolding for the word 
$W_n$ consists of two maximal $n-1$ darts.   Given our bounds on $x$,
we satisfy the hypotheses given in the remark immediately following
the statement of the Dart Lemma.  Thus, it suffices to consider the
superior vertices of the unfolding.
There are $4$ superior top vertices, labelled $A,B,C,D$.  Likewise,
there are $4$ superior bottom vertices, labelled $E,F,G,H$.  See
Figure 3.1.  Thus, by the
Dart Lemma, it suffices to show that $X \uparrow Y$ for
each $X \in \{A,B,C,D\}$ and $Y \in \{E,F,G,H\}$.
We normalize coordinates so that $A=(0,0)$, and the long side has length one. Then, we can compute the coordinates for the $3$-spine.\
$$E=(\sin (n-1)x, -\cos (n-1)x). \quad D=(2 \sin (n-1)x,0). \quad H=(3\sin (n-1)x, -\cos (n-1)x).$$
With this choice, the unfolding is horizontal as desired. (That is, $A \updownarrow D$.)

\begin{center}
\includegraphics[width=4.2in]{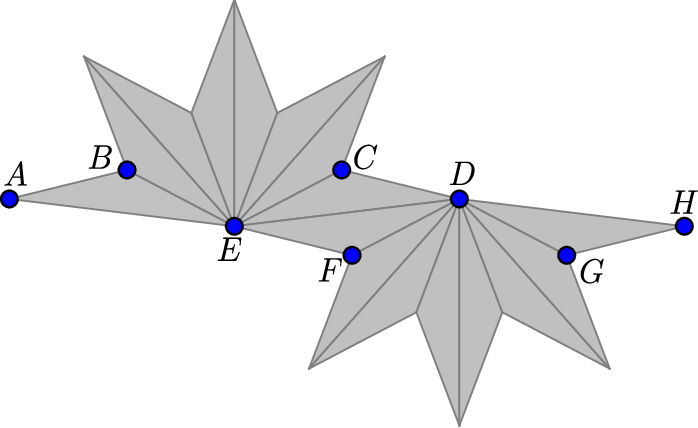}
\newline
Figure 3.1: An unfolding for the word $W_5$. One period is shown, which begins at edge 
$\overline{AE}$ and ends at the parallel edge $\overline{DH}$.
\end{center}

By symmetry, each point has a partner-point at the same height:
$$A \updownarrow D \quad B \updownarrow C \quad E \updownarrow H \quad F \updownarrow G$$
Thus, it is sufficient to concentrate on the central rhombus, $CDFE$. Given a vector ${\bf v} \in \R^2$, we use
$d({\bf v})$ to denote the angle of the vector made with the horizontal. It is sufficient to check that the four vectors
$\overrightarrow{EC}$, $\overrightarrow{ED}$, $\overrightarrow{FC}$, and $\overrightarrow{FD}$ point upward. That is, for ${\bf v}$ equal
each of those four vectors, we must have $0<d({\bf v})<\pi$.
We compute
$$d(\overrightarrow{FC})=d(\overrightarrow{EC})=\frac{\pi}{2}-(n-2)x \quad d(\overrightarrow{ED})=\frac{\pi}{2}-(n-1)x
\quad  d(\overrightarrow{FC})=\pi-(n-1)x$$
In all cases, we have $0<d({\bf v})<\pi$ for $0<x<\frac{\pi}{2n-2}$.
\endproof

\subsection{The stable family $Y_{n,m}$}

The word $Y_{n,m}$ has an additional special symmetry. 
If you write $Y_{n,m}$ in reverse and swap the letters $1$ and $2$,
you get $Y_{n,m}$ back. Given a word $W$, let $\widehat{W}$ denote $W$ written in reverse with $1$ swapped with $2$.   There is some word $W=W_{m,n}$ such that
\begin{equation}
\label{specialsymm}
Y_{m,n}=1W3\widehat W 2,
\end{equation}
\noindent
{\bf Remark:\/}
It is a consequence of work in \cite{H2} that every stable periodic billiard path in an isosceles triangle has an combinatorial type
$W$ with the symmetry $\widehat{W}=W$.   This fact, however, is not necessary for
our proof here.
\newline

We now record some special properties of words having the form given by the
right hand side of Equation \ref{specialsymm}.

\begin{proposition}
\label{prop:special_word}
Let $Y$ be a word of the form $Y=1W3 \widehat{W}2$, and $T$ be an obtuse isosceles triangle. Consider the unfolding
$U(Y^2,T)$ chosen so that the translation bringing the first edge to the last is horizontal. Then the long edge (edge 3)
of the first triangle in the unfolding is horizontal.
\end{proposition}
\startproof
Consider the bi-infinite repeating word $\overline{Y}$. 
This word has some symmetry, which is revealed by expanding the word out.
$$\overline{Y}=\ldots 1W3\widehat{W}~21W3\widehat{W}2~\Big|~1W3\widehat{W}~21W3\widehat{W}$$
Reflection in the vertical line above swaps the letters $1$ and $2$ while preserving $3$. This is precisely how the reflective symmetry
of the isosceles triangles permutes the labeling of the sides. Thus, this symmetry extends to the bi-infinite unfolding 
$U(\overline{Y},T)$. The direction of the holonomy of $U(Y^2,T)$ must be the eigenvector
corresponding to eigenvalue $-1$ of the reflective symmetry of $U(\overline{Y},T)$. But this reflection is just the 
reflective symmetry of the first triangle in the unfolding. So these two directions are parallel.
\endproof

We will use the following principle for detecting our billiard path.
Recall that side $3$ denotes the long side of an isosceles triangle.

\begin{proposition}
\label{prop:symmetry_implication}
Suppose that a billiard path in an obtuse isosceles triangle
starts out parallel to side $3$, and has initial combinatorial
type $1W3$, where the final $3$ corresponds to an edge which
the path hits at the midpoint, $M$.   Then the billiard path
is closed and has combinatorial type $1W3\widehat W 2$.
\end{proposition}

See Figure 3.2 for a case when this proposition applies. 

\begin{center}
\includegraphics[width=6in]{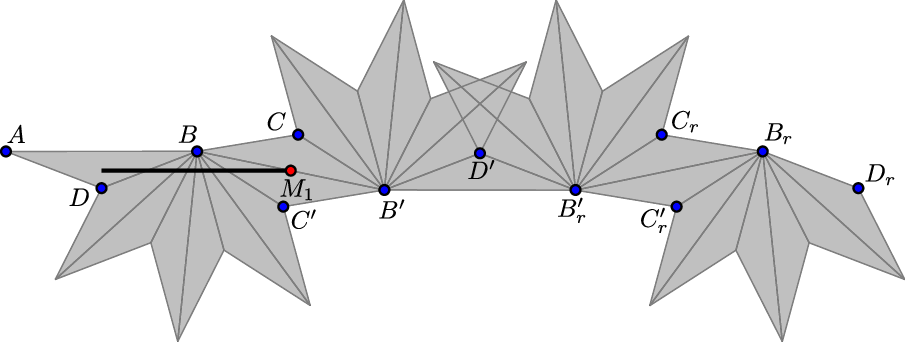}
\newline
Figure 3.2: An unfolding for the square of the word $Y_{4,1}$.
\end{center}

\startproof
The trajectory $t_1$ described in the proposition lies within the unfolding of the initial word $1W3$ and then hits $M$.
The unfolding for the word $1W3\widehat W 2$ has 180 degree rotational symmetry $\phi$ around the point $M$. 
Thus, the longer trajectory $t_2=t_1 \cup \phi(t_1)$ lies within the unfolding of $1W3\widehat W 2$.
Now consider the unfolding of the even length word $(1W3\widehat W 2)^2$. 
This unfolding has vertical reflective symmetry $\rho$ which swaps the two halves of the word. (It is vertical assuming the trajectory is horizontal.) The trajectory $t_3=t_2 \cup \rho(t_2)$ lies within the unfolding of $1W3\widehat W 2$.
\endproof

We will break the proof of Theorem \ref{thm:Ycovers} into two cases. 
The first case is easiest.

\begin{lemma}
For each $x$ satisfying $\frac{\pi}{2n+1}\leq x<\frac{\pi}{2n}$ there is a periodic billiard path in $T_x$ with combinatorial type $Y_{n,1}$.
\end{lemma}

\startproof
Given the triangle $T_x$, unfold the triangle according to the square of the word $Y_{n,1}$ as in Figure 3.2.
Let $M_1$ be the midpoint we must hit. This is the first midpoint of a long side which is the fixed point of 
a 180 degree rotational symmetry of the the bi-infinite unfolding, $U(\overline{Y_{n,1}},T_x)$.

We coordinatize the unfolding so that $M_1$ is given coordinates $(0,0)$. We will show that all the top vertices have
positive $y$-coordinate, and all the bottom vertices have negative $y$-coordinate. Regardless of $n$, 
the Dart Lemma tells us that most of the vertices are irrelevant. It is enough to prove this statement for those
vertices, who are given names in Figure 3.2. We have named four vertices $A$, $B$, $C$ and $D$. 
The other vertices are either images of these under the rotational symmetry about $M_1$ (denoted by $\ast'$), 
images under reflection in the vertical line through $D'$ (denoted by $\ast_r$), or images under the composition.
So, it is enough to show that the statement is true for the vertices $A$, $B$, $C$ and $D$. 

The points $A$ and $B$ have the same $y$-coordinate by Proposition \ref{prop:special_word}.
$M_1$ lies below them, because angle $\angle ABM_1=2n x< \pi$. Also angle $\angle ABC=(2n+1)x \geq \pi$, so $C$ has $y$-coordinate 
greater than or equal to the $y$-coordinates of $A$ and $B$. Finally, $M_1$ lies closer to $B$ then $D$. Furthermore, 
the vector $\overrightarrow{BM_1}$ is closer
to horizontal than the vector $\overrightarrow{DB}$. ($\angle DBA=x=\angle M_1BC$, but the horizontal
direction lies strictly between the directions of $\overrightarrow{BM_1}$ and $\overrightarrow{BC}$.)
Thus the $y$-coordinate of $D$ must be negative.
\endproof

\noindent
{\bf Remark:\/}
The words $Y_{n,1}$ are the same as the words $A_{n-1}$, which appear in Figure 1.1 
and play a prominent role in Theorem \ref{easy}.
\newline

The second case is more complicated. While we could give a constructive proof, as above, we
find that a non-constructive proof clarifies the situation.
To illustrate this case, we consider the word $Y_{4,2}$ and the triangle $T$ in Figure 3.3. The unfolding 
$U(Y_{4,2},T)$ depicted in this Figure contains a horizontal segment joining the first triangle to the midpoint of the long side.
This segment hits the sequence of sides $1(31)^4(32)^43$. Let $W=(31)^4(32)^4$. By 
Proposition \ref{prop:symmetry_implication}, there is a periodic billiard path in $T$ with combinatorial type 
$Y_{4,2}=1W3\widehat{W}2$.
The significant point is that by Proposition \ref{prop:symmetry_implication},
we only need to consider the unfolding for an initial subword.

\begin{center}
\includegraphics[width=6in]{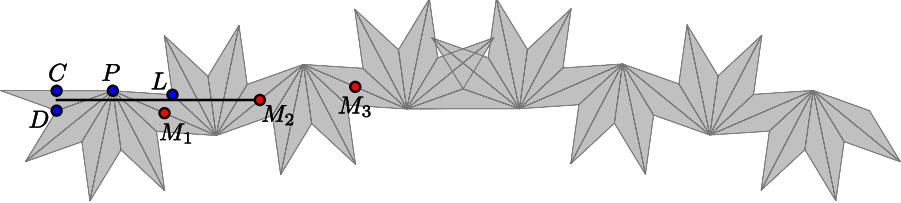}
\newline
Figure 3.3: An unfolding of the square of the word $Y_{4,2}$.
\end{center}

\begin{lemma}
\label{lem:rest}
For each $x$ with $\frac{\pi}{2n+2}<x<\frac{\pi}{2n+1}$, 
there is a periodic billiard path in $T_x$ with combinatorial type $Y_{n,m}$
for some $m \in \N$.
\end{lemma}
\startproof
Consider the unfolding of $T_x$ according to the infinite word $1 (\u_n)^{\infty}$. See Figure 
3.4. We normalize the unfolding so that the initial
long side of $T_x$ is horizontal. 

We will show that there is some index $m$ such that $M_m$ lies below all preceding
top vertices and above all preceding bottom vertices.   Here the points
$M_0,M_1,...$ are the midpoints of some of the long segments.  See Figure 3.4.

\begin{center}
\includegraphics[width=6in]{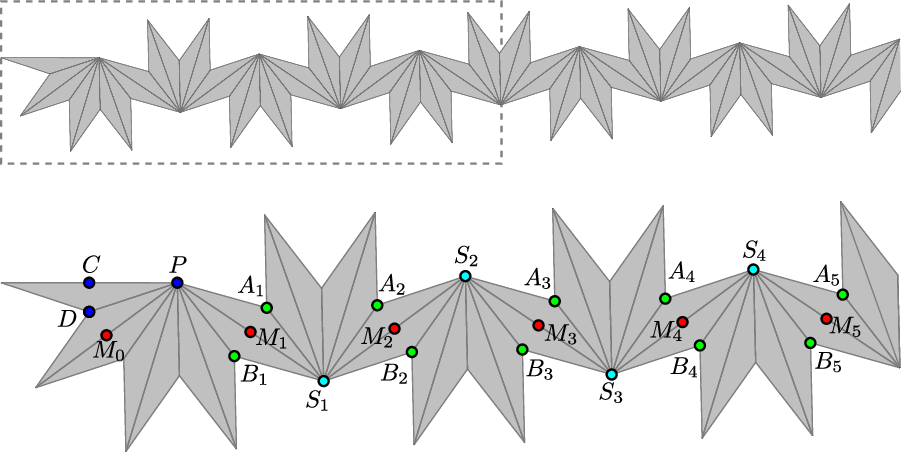}
\newline
Figure 3.4: The top figure shows an unfolding for the word $1 (\u_4)^{\infty}$. The bottom figure shows a blow up of the unfolding with names of important vertices labeled.
\end{center}

Many of the vertices in Figure 3.4 are given names.  The unlabelled vertices are
inferior, and may be ignored by the Dart Lemma.
Understanding this unfolding is made much easier by the fact that
$\u_n$ is the combinatorial type of a periodic billiard path in $T_x$.
See Proposition  \ref{prop:unstable}.
We consider the vector ${\bf v}=\overrightarrow{M_1M_2}$.
This vector points in the direction $(n+1)x-\frac{\pi}{2}$ (measured relative to the horizontal in polar coordinates). 
In particular, ${\bf v}$ has positive $y$-coordinate, since $x>\frac{\pi}{2n+2}$.
To compute this direction, note that $\overrightarrow{S_1M_1}$ points in the direction $2nx$ and $\overrightarrow{S_1M_2}$
points in direction $2x$. We always normalize so that the long side has length $1$. This gives us the formula 
$${\bf v}=\frac{1}{2} (\cos 2x, \sin 2x) - \frac{1}{2} (\cos 2nx, \sin 2nx).$$

Now we will eliminate some top vertices. The vector $\overrightarrow{P A_1}$ points in direction $\pi+2nx<2 \pi$. So
$A_1 \downarrow P$. Also $\overrightarrow{A_1 A_k}$ is parallel to ${\bf v}$ for $k>1$, so 
$A_1 \downarrow A_k$. And $\overrightarrow{A_{2k} S_{2k}}$ points in direction $x$, so $A_{2k} \downarrow S_{2k}$.
It follows that $A_1$ has the least $y$-coordinate of any vertex in the unfolding $1 (\u_n)^{\infty}$.

Now we will eliminate some bottom vertices. $\overrightarrow{B_k B_{k+1}}$ points in the same direction as ${\bf v}$, so 
$B_{k+1} \uparrow B_k$. The vector $\overrightarrow{B_{2k+1} S_{2k+1}}$ points in direction
$(2n+1)x-\pi<0$, so $B_{2k+1} \uparrow S_{2k+1}$. It follows that the bottom vertex with greatest $y$-coordinate in the unfolding up to 
the appearance of $M_m$ is either $D$ or $B_m$. (We ignore the fact that when $m$ is even $B_m$ appears later in the unfolding
than $M_m$.) 

Letting $P_y$ denote the $y$-coordinate of an arbitrary point $P$, we want to
show that
$$(M_m)_y>(B_m)_y; \hskip 50 pt
D_y<(M_m)_y<(A_1)_y,$$
for some $m$.
When $m$ is even, $\overrightarrow{B_m M_m}$ points in the direction $2x+\frac{\pi}{2}$, so 
$M_m \uparrow B_m$. When $m$ is odd, $\overrightarrow{B_m M_m}$ points in the direction $2nx-\frac{\pi}{2}$. So, for any choice of $m$, the first inequality holds.
The first endpoint
$M_0$ lies below $D$ since $\overrightarrow{M_0 D}$ points in the direction $2x+\frac{\pi}{2}$.
Let \begin{equation}
f(x)=(A_1)_y-D_y. \hskip 30 pt
g(x)=(M_{i+1})_y-(M_i)_y.
\end{equation} 
Here $g(x)$ is independent of $i$.
We compute
\begin{equation}
g(x)=\frac{1}{2} (\sin 2x-\sin 2nx); \hskip 30 pt
f(x)=\frac{1}{2 \cos x}(\sin x- \sin (2n+1)x).
\end{equation}
The formula for $f$ follows from the fact that the length
of the short side is $\frac{1}{2 \cos x}$, and $\overrightarrow{P D}$ and
$\overrightarrow{P A_1}$ point in directions $\pi+x$ and $\pi+(2n+1)x$ respectively.

A sufficient criterion for the second equation is that $g(x)<f(x)$.
That this is true follows from some trigonometry. First we reduce $f(x)$ and $g(x)$ to more convenient forms. 
\begin{equation}
\cos (x) f(x) = -\sin nx \cos (n+1)x 
\quad \textrm{and} \quad
g(x)=-\sin (n-1)x \cos (n+1)x.
\end{equation}
Thus, $$\frac{f(x)}{g(x)}=\frac{\sin nx }{\cos x \sin (n-1)x}>\frac{1}{\cos x}>1.$$
This completes the proof.
\endproof


\newpage

%% file: 4abc.tex
\section{Proof of Theorem \ref{easy}}
\label{sect:easy}

\subsection{The A Family}

Here we introduce the words $\{A_n\}$ for $n \geq 2$.
These words
already appear in \cite{HH}, and their analysis
is quite easy.
$A_n$ is the square of a word of odd
length.  Listing out the first few words explicitly and
then writing the general pattern, we have:
\begin{equation}
A_2=(2323131)^2; \hskip 15 pt
A_3=(23232313131)^2; \hskip 15 pt
A_n=((23)^n(13)^{n-1}1)^2
\end{equation}
The squarepath $\widehat Q(A_n)$ is a square of sidelength
$2n$.   Hence $U(A_n,*)$ is the
union of $4$ maximal $n$-darts.
The $3$-spine for $U(A_n,V_n)$ is contained in a
straight line.  There are two top vertices on
this straight line and two bottom vertices.
The top vertices are $a_1$ and $a_2$.  The
bottom vertices are $b_{2n}$ and $b_{2n+1}$.
Figures 4.1, 4.2, and 4.3 show the first few examples.

\begin{center}
\includegraphics[height=1.2in]{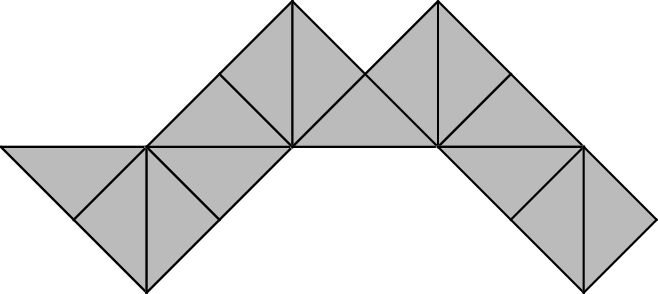}
\newline
Figure 4.1 $U(A_2,V_2)$.
\end{center}

\begin{center}
\includegraphics[height=1.2in]{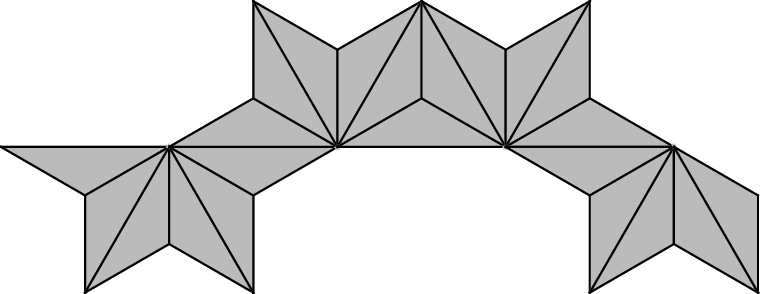}
\newline
Figure 4.2 $U(A_3,V_3)$.
\end{center}

\begin{center}
\includegraphics[height=1.2in]{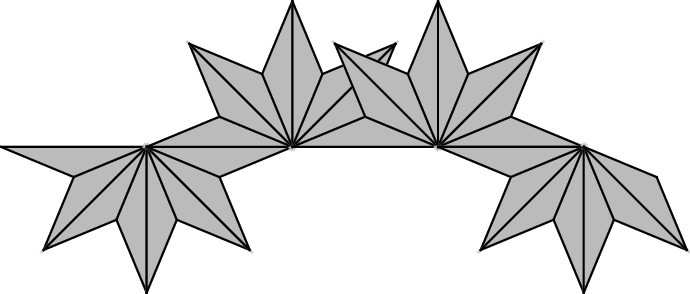}
\newline
Figure 4.3 $U(A_4,V_4)$.
\end{center}

If $X$ is a parameter point sufficiently close to $V_n$
then the lowest top vertices of $U(A_n,X)$ remain $a_1$ and $a_2$
and the highest bottom vertices remain $b_{2n}$ and
$b_{2n+1}$.   
When $X \in N_{--}(n,\epsilon)$ the $3$-spine for
$U(A_n,X)$ is no longer a straight line segment,
but rather makes a zig-zag.  Both obtuse angles
in the unfolding are slightly smaller, and this
causes the $3$ spine to make an acute angle in
the directions of the centerlines of the maximal
darts, as shown in Figures 4.4 and 4.5.

\begin{center}
\includegraphics[height=1.3in]{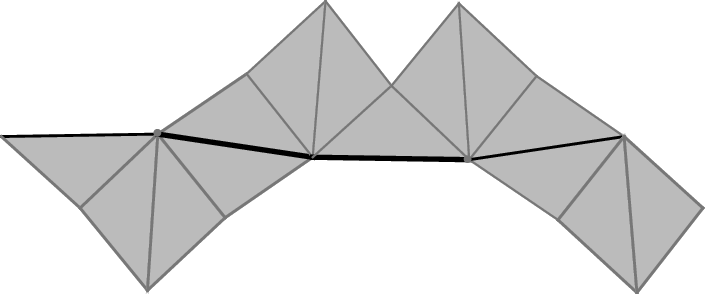}
\newline
Figure 4.4: $U(A_2,X)$ for $X \in N_{--}$
\end{center}

\begin{center}
\includegraphics[height=1.3in]{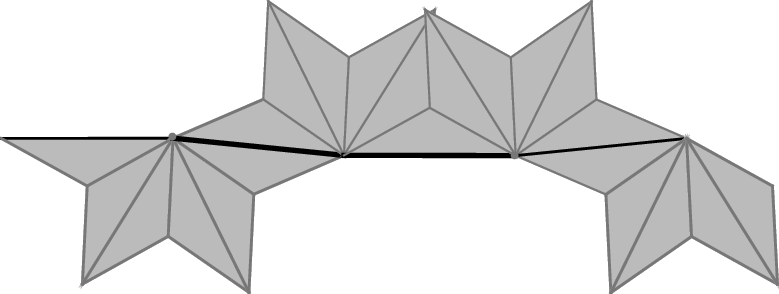}
\newline
Figure 4.5 $U(A_3,X)$ for $X \in N_{--}$.
\end{center}

From this geometric picture we see easily that the $a_1$ and $a_2$
lie above $b_{2n}$ and $b_{2n+1}$ for points in 
$N_{--}(n,\epsilon_n)$.

As an alternate argument we note that, since $A_n$ is an odd
square, the unfolding $U(A_n,*)$ has glide-reflection symmetry.
Thus, if $a_i$ is the lowest top vertex then $b_{2n+i-1}$ is the
highest bottom vertex.  Thus, it suffices to show that
$a_1 \uparrow b_{2n}$ and $a_2 \uparrow b_{2n+1}$.  We
compute the defining
function $F$ for $(a_1,b_{2n})$ and find that
\begin{equation}
\label{FF}
F(x_1,x_2)= -4 \sin^2(nx)\sin(2nx).
\end{equation} 
For $(x_1,x_2)$ near $V_n=(\pi/2n,\pi/2n)$,
the above expression is negative iff
$x_2<\pi/2n$.
That is, $a_1 \uparrow b_{2n}$ iff $x_2<\pi/2n$
and $x_2$ is sufficiently close to $\pi/2n$.
The calculation for the
pair $(a_2,b_{2n+1})$ yields the same result, but
with $x_1$ and $x_2$ interchanged. 

This takes care of the first statement of Theorem \ref{easy}.

\subsection{The B Family}

We will show that
\begin{equation}
\label{veech1}
N_{-+}(n,\epsilon_n) \subset O(B_n); \hskip 30 pt n=4,5,6...
\end{equation}
By symmetry,
\begin{equation}
\label{veech11}
N_{+-}(n,\epsilon_n) \subset O(C_n); \hskip 30 pt n=4,5,6...
\end{equation}
Our argument will show that the two segments bounding
$N_{--}(n,\epsilon_n)$ (except for $V_n$ itself)
 are respectively contained in
$O(B_n)$ and $O(C_n)$. This takes care of the fourth
statement of Theorem \ref{easy}.

The word $B_n$ has length $40n-60$.  This
word is determined by its
squarepath $\widehat Q_n:=\widehat Q(B_n)$, which we now describe.
We will draw $\widehat Q_4$ and $\widehat Q_5$, with
the understanding that $\widehat Q_{n+1}$ is obtained from
$\widehat Q_n$ by lengthening each edge by $2$ units. 
The small grey squares in Figure 4.6 have 
edgelength $2$.
We have drawn some of the edges in grey to help the
reader parse the loops.  These loops are homeomorphic
to figure $8$ curves.
The grey dot indicates the origin.

\begin{center}
\includegraphics[height=3.5in]{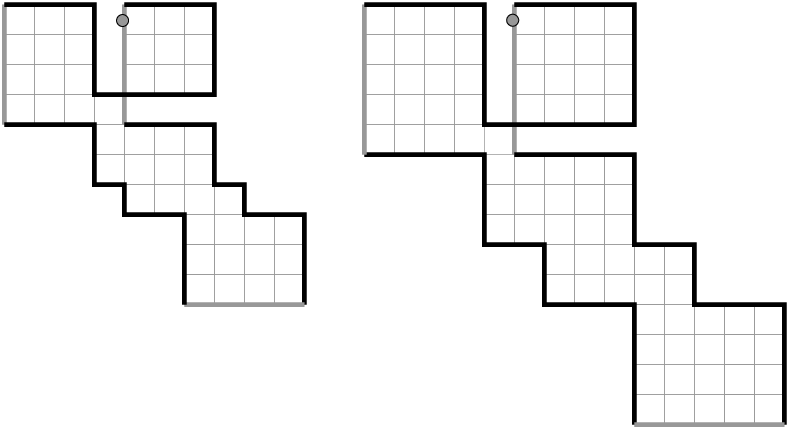}
\newline
Figure 4.6: $\widehat Q_n$ for $n=4,5$.
\end{center}

The shortest unfolding $U(B_4,V_4)$ has $100$ triangles in it. Here
it is.  We have highlighted the $3$-spine. 

\begin{center}
\includegraphics[width=6.25in]{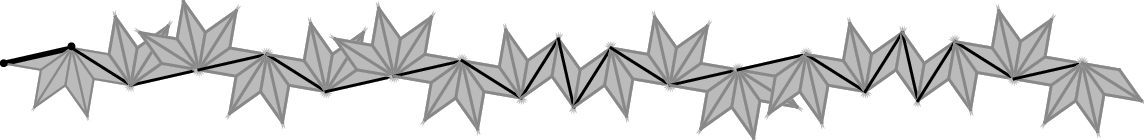}
\newline
Figure 4.7: $U(B_4,V_4)$.
\end{center}

Equation \ref{master} shows three lists.  The list $L_j$ is the list
of $j$th coordinates of the successive vertices of the squarepath.
The list $L$ computes either $L_1+L_2$ mod $4n$ or $L_1+L_2+2n$ mod $4n$
depending on the parity of the vertex.   

\begin{equation}
\label{master}
\begin{matrix}{\bf L:\/} \cr 
1\cr-1\cr1\cr1\cr-1\cr1\cr1\cr-1\cr3\cr-3\cr3\cr-1\cr1\cr1\cr-1\cr3\cr-3\cr3\cr-1\cr 1 \end{matrix}
\hskip 30 pt
\begin{matrix}
{\bf L_1:\/} \cr 
    0 n +0\cr
    2 n -2\cr
    2 n -2\cr
    0 n -2\cr
    0 n -2\cr
   -2 n +0\cr
   -2 n +0\cr
    0 n -2\cr
    0 n -2\cr
    2 n -8\cr
    2 n -8\cr
    4 n -12\cr
    4 n -12\cr
    6 n -12\cr
    6 n -12\cr
    4 n -8\cr
    4 n -8\cr
    2 n -2\cr
    2 n -2\cr
    0 n +0 \end{matrix}    
\hskip 30 pt
\begin{matrix}{\bf L_2:\/}
\cr 
0 n +1\cr
    0 n +1\cr
   -2 n +3\cr
   -2 n +3\cr
    0 n +1\cr
    0 n +1\cr
   -2 n +1\cr
   -2 n +1\cr
   -4 n +5\cr
   -4 n +5\cr
   -6 n +11\cr
   -6 n +11\cr
   -8 n +13\cr
   -8 n +13\cr
   -6 n +11\cr
   -6 n +11\cr
   -4 n +5\cr
   -4 n +5\cr
   -2 n +1\cr
   -2 n +1 \end{matrix}
\end{equation}

Let
\begin{equation}
\omega_n=E(\frac{\pi}{2n}).
;\hskip 30 pt
E(x)=\exp(ix).
\end{equation} 
We will often write $\omega=\omega_n$ when the dependence on
$n$ is clear.
The holonomy of $U(B_n,V_n)$ is obtained as the alternating sum
of the vertices of $\widehat Q$.
Since we are evaluating this sum at $V_n$, each vertex
contributes some power of $\omega$ to the sum.
The list $L$ above tells us which power.
In deriving this list, we used the relation
$\omega^{a+2n}=-\omega^a$.  Note that $L$ is independent of $n$.
We have:
\begin{equation}
Q(V_n)=\sum_{i=1}^{20} \omega^{L(i)}=
8\omega+4\omega^3+6\omega^{-1}+2 \omega^{-3}.
\end{equation}
This agrees with the McBilliards Calculations.
An easy calculus argument shows that
$Q(V_n)$ lies between $1$ and $\omega$ on the unit
circle. 

\subsection{Reducing to Six Vertices}

Any point $X$ sufficiently near $V_n$ satisfies the hypothesis
of the Dart Lemma with respect to $B_n$.  Hence, we just have
to show that all the top superior vertices of $U(B_n,X)$ lie
above all the bottom superior vertices of $U(B_n,X)$ for
$X \in N_{-+}(n,\epsilon)$ when $\epsilon$ is
sufficiently small.  In this section we will reduce this
problem to checking $6$ superior vertices.

Since $U_n$ decomposes
into $20$ maximal darts, there are at most $80$
superior vertices, independent of $n$. 
Each maximal dart has $2$ (O)uter superior vertices and
$2$ (I)nner superior vertices. The outer superior vertices
lie on the $3$-spine and the inner superior vertices
do not.  The superior vertices in each maximal dart
are naturally ordered from left to right.  There are
two superior vertices on the (L)eft and two on the (R)ight.
We denote the $4$ superior vertices of the $K$th maximal dart
(perhaps redundantly) by
\begin{equation}
(K,L,O); \hskip 15 pt
(K,L,I); \hskip 15 pt
(K,R,I); \hskip 15 pt
(K,R,O).
\end{equation}
Sometimes we will decorate our notation with an asterisk
to indicate whether it is a top vertex $(\ )^*$ or a
bottom vertex $(\ )_*$.

A {\it leader\/} is either a lowest top vertex
or a highest bottom vertex.  Here is the main result of
this section:

\begin{lemma}
\label{elim}
For each $n$, the leaders of $U(B_n,V_n)$ are
$$(1,L,I)_*; \hskip 10 pt
(4,L,O)_*; \hskip 10 pt
(16,L,I)_* \hskip 10 pt
(5,R,I)^*; \hskip 10 pt
(10,R,I)^*; \hskip 10 pt
(12,R,O)^*.
$$
Moreover, these points all have the same height.
\end{lemma}

Any top superior vertex $a_2$ which is not on our list should lie above
$b_1=(1,L,I)_*$.
We will symbolically 
compute 
$$F(V_n)={\rm Im\/}(P(V_n) \overline Q(V_n))$$ for such
pairs and show that the imaginary part of this function
is positive. Hence $a_2 \uparrow b_1$.  
Likewise, if $b_2$ is a bottom superior vertex not on the list, 
we will show that $F(V_n)<0$ where $F$ is the function
corresponding to $(b_1,b_2)$.
Finally we will show that the $P(V_n)$ is a real multiple
of $Q(V_n)$ when $P$ is defined relative to the pair $(b_1,c)$
and $c$ is on the list given in Lemma \ref{elim}.
The key to our calculations is a slick procedure
for computing these points.

In computing our points we will slightly modify the
method described in \S \ref{sect:background}, so as to use the
$3$-spine as much as possible.  Unfortunately,
it is not possible to directly connect all the
points of interest to us by a $3$-path. 
The work-around we explain below works for every
point except for
$a_1=(1,L,O)^*$, which we easily observe to
lie above the points listed in Lemma \ref{elim}.
(The edge connecting $(1,L,O)$ to $(1,L,I)$ has
negative slope.) For the rest of the
superior vertices we do the following:

\begin{itemize}
\item Connect $(1,L,I)_*$ to $(1,L,O)^*$ using the
common edge $e_0$.
\item Connect $(I,L,O)_*$ to a point $p'$ on the
$3$-spine which is adjacent to $p$, using the
fewest number of edges $e_1,...,e_s$ from the $3$-spine.
\item Connect $p'$ to $p$ by the edge $e_{s+1}$
which is incident to both vertices.
\end{itemize}

Once $p'$ is determined, there are either $2$ or
$3$ choices for $p$.
To use an analogy, the $3$-spine
is like the highway and the other edges we use
are like the off and on ramps.  Our first step is
to get onto the highway using $e_0$. Then
we drive along the highway using
$e_1,...,e_s$.  At this point we can either
take one of the off-ramps and stop the car or
else go one more mile and stop the car, depending
on our final destination.  Figure 4.8 shows a
fairly accurate picture for $s=2$.  
The dotted lines indicated some
of the triangles in the unfolding.

\begin{center}
\includegraphics{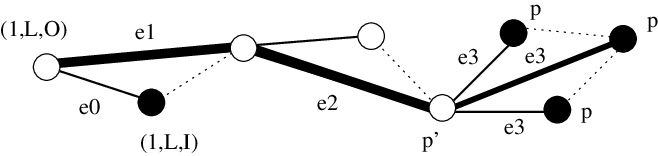}
\newline
Figure 4.8: The Connecting Path
\end{center}

We normalize so that the type $3$ (long) edges of
our triangles have length $1$.  It follows from the
Law of Sines that the short sides have length
\begin{equation}
\lambda=\frac{1}{\omega+\omega^{-1}}.
\end{equation}
The vector that points from
$(1,L,I)_*$ to $p$ is 
\begin{equation}
\label{elimX}
P(s,\delta)=-\lambda+\sum_{i=1}^{s} \omega^{L(i)}
+ \lambda^{|\delta|} \omega^{L(s+1)+\delta}.
\end{equation}
Referring to Figure 4.8,, the number
$\delta$ is $-1$ (bottom) or $0$ (middle) or $1$ (top) depending
on which of the three choices we make for $e_{s+1}$.
Here $L$ refers to the labeling in Equation \ref{master}.

To demonstrate Equation \ref{elimX}, we note that the
three paths suggested by Figure 4.8 lead to the three sums
$$-\lambda+\omega+\omega^{-1}+\lambda \omega^0; \hskip 30 pt
-\lambda + \omega + \omega^{-1} + \omega; \hskip 30 pt
-\lambda + \omega + \omega^{-1} +\lambda \omega^2.
$$

Let's concentrate on the right sum.  We will
use the notation $x \to y$ to denote that ${\rm Im\/}(x)$ and ${\rm Im\/}(y)$
are positive multiples of each other.
Our middle expression simplifies to
$$P=\frac{1+\omega^2+2\omega^4}{\omega(1+\omega^2)}.$$
Recalling our formula for the holonomy, we compute that
$$P \overline Q=
\left[\begin{matrix}0\cr 4\cr 8\cr 14\cr 10\cr 4\cr 0 \end{matrix} \right] \cdot 
\left[\begin{matrix}\omega^{-6} \cr \omega^{-4} \cr \omega^{-2} \cr \omega^0 \cr 
\omega^2 \cr \omega^4 \cr \omega^6\end{matrix}\right]\to 
\left[\begin{matrix}10\!-\!8 \cr 4\!-\!4 \cr 0 \end{matrix} \right] \cdot
\left[\begin{matrix}\omega^2 \cr \omega^4 \cr \omega^6\end{matrix}\right]=
\left[\begin{matrix}2 \cr 0 \cr 0\end{matrix}\right] \cdot
\left[\begin{matrix}\omega^2 \cr \omega^4 \cr \omega^6\end{matrix}\right] \to
\left[\begin{matrix}2 \cr 0 \cr 0\end{matrix}\right] \cdot
\left[\begin{matrix}\sin(1\pi/n) \cr \sin(2 \pi/n) \cr \sin(3 \pi/n)\end{matrix}\right].$$
Hence ${\rm Im\/}(F(V_n))>0$.  
The point of $U_n$ corresponding to the third sum above is
$(3,R,I)^*$, and this shows that
$(3,R,I)^* \uparrow (1,L,I)_*$.
It turns out that our sums always lead to the general expression
$$\left[\begin{matrix}c_{-6} \cr c_{-4} \cr c_{-2} \cr c_0 \cr c_2 \cr c_4 \cr c_6 \end{matrix} \right] \cdot 
\left[\begin{matrix}\omega^{-6} \cr \omega^{-4} \cr \omega^{-2} \cr \omega^0 \cr 
\omega^2 \cr \omega^4 \cr \omega^6\end{matrix} \right]\to
\left[\begin{matrix}c_2-c_{-2} \cr c_4-c_{-4} \cr c_6-c_{-6} \end{matrix}\right] \cdot
\left[\begin{matrix}\omega^2 \cr \omega^4 \cr \omega^6 \end{matrix} \right] \to
\left[\begin{matrix}a_1 \cr a_2 \cr a_3 \end{matrix} \right] \cdot
\left[\begin{matrix}\sin(1\pi/n) \cr \sin(2\pi/n) \cr \sin(3 \pi/n) \end{matrix} \right]
 \hskip 20 pt
a_j,c_j \in \Z.$$
In listing the results
of our calculations it suffices to list vector $(a_1,a_2,a_3)$.

For each $\delta \in \{-1,0,1\}$ and each
$\beta \in \{1,...,20\}$ we compute $P(\delta,s) \overline Q$ and
extract the coefficient vector $(a_1,a_2,a_3)$.

Here is the table for $\delta=-1$ and $\beta=1,...,10$.
\begin{equation}
\begin{matrix}
(0)&(-)&(-)&(-)&(-)&(-)&(-)&(-)&(-)&(-) \cr
0&0&-4&-2&-2&-6&-4&-4&-2&0 \cr
0&-2&-2&0&-2&-2&0&-2&0&0 \cr
0&0&0&0&0&0&0&0&0&-2\end{matrix}
\end{equation}

Here is the table for $\delta=-1$ and $\beta=1,...,10$.
\begin{equation}
\begin{matrix}
(-)&(-)&(-)&(-)&(-)&(0)&(-)&(-)&(-)&(-) \cr
0&0&-4&-2&-2&0&2&2&2&-2 \cr
-2&-2&-2&0&-2&0&0&-2&-2&-2 \cr
-2&0&0&0&0&0&-2&-2&0&0\end{matrix}
\end{equation}

Here is the table for $\delta=0$ and $\beta=1,...,10$.
\begin{equation}
\begin{matrix}
(+)&(-)&(0)&(+)&(-)&(-)&(+)&(-)&(+)&(-) \cr
4&-2&0&2&-4&-2&0&-6&2&-4 \cr
2&-2&0&2&-2&0&2&-2&4&-4 \cr
0&0&0&0&0&0&0&0&2&-2\end{matrix}
\end{equation}

Here is the table for $\delta=0$ and $\beta=11,...,20$.
\begin{equation}
\begin{matrix}
(+)&(-)&(0)&(+)&(-)&(+)&(-)&(+)&(-)&(+) \cr
4&-2&0&2&-4&4&-2&6&0&2 \cr
2&-2&0&2&-2&4&-4&2&-2&0 \cr
0&0&0&0&0&2&-2&0&0&0\end{matrix}
\end{equation}

Here is the table for $\delta=0$ and $\beta=1,...,10$.
\begin{equation}
\begin{matrix}
(+)&(+)&(+)&(+)&(+)&(0)&(+)&(+)&(+)&(+) \cr
6&2&2&4&0&0&2&-2&-2&-2 \cr
2&2&0&2&2&0&2&2&2&0 \cr
0&0&0&0&0&0&0&0&2&2\end{matrix}
\end{equation}

Here is the table for $\delta=0$ and $\beta=11,...,20$.
\begin{equation}
\begin{matrix}
(0)&(+)&(+)&(+)&(+)&(+)&(+)&(+)&(+)&(+) \cr
0&2&2&4&0&0&0&2&4&4 \cr
0&2&0&2&2&2&0&0&2&0 \cr
0&0&0&0&0&2&2&0&0&0\end{matrix}
\end{equation}

All the expressions on the first two tables correspond to bottom vertices.
An inspection of Figure 4.7 shows that the expressions
with $(+)$ signs correspond to top vertices and the expressions with
$(-)$ signs correspond to bottom vertices. (We copy the figure above.)  All the
expressions on the bottom two tables correspond to top vertices.
Finally, there are $6$ expressions which are real, independent of $n$,
and these correspond to the $6$ vertices in Lemma \ref{elim}.
This completes the proof of Lemma \ref{elim}.

\subsection{The End of the Proof}

Lemma \ref{elim} implies that $O(B_n)$ is determined by 
positions of the $6$ vertices
$$
\alpha_1=(12,R,O)^*; \hskip 15 pt
\alpha_2=(5,R,I)^*; \hskip 15 pt
\alpha_3=(10,R,I)^*.
$$
\begin{equation}
\beta_1=(4,L,O)_*; \hskip 15 pt
\beta_2=(1,L,I)_*; \hskip 15 pt
\beta_3=(16,L,I)_*;
\end{equation}

Let $H_{ij}$ be the defining function which measures
the height of $\alpha_i$ minus the height of
$\beta_j$, when the unfolding is normalized to
that the long edges have length $1$.   Below we will prove

\begin{lemma}
\label{final}
For every relevant pair of indices, we have
$\partial_1 H_{ij}(V_n)<0$ and
$\partial_2 H_{ij}(V_n) \geq 0.$
There is equality in the second equation iff
$(i,j)=(1,1)$.
\end{lemma}

Note that $H_{ij}(V_n)=0$.  Lemma \ref{final} therefore
says that $H_{ij}(X)>0$ provided that
$(i,j) \not = (1,1)$ and $X \in N_{-+}(n,\epsilon_n)$
for sufficiently small $\epsilon_n$.  If we knew that
$H_{11}$ vanished identically on the line
$x_1=\pi/2n$ then we could conclude the same result
for $(i,j)=(1,1)$.   Before proving Lemma \ref{final}
we will take care of the exceptional pair of indices.

\begin{lemma}
\label{yline}
$H_{11}$ vanishes identically
on the line $x=\pi/2n$.
\end{lemma}

\startproof
We will draw the picture for the case $n=4$, but the phenomenon we
describe is completely general.  Let $P$ and $Q$ and $F$ be the
defining functions associated to our two points. Figure 4.9 shows
the paths $\widehat Q$ and $\widehat P$.

\begin{center}
\includegraphics{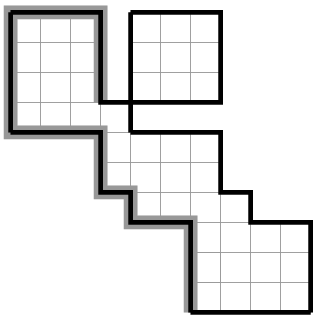}
\newline
Figure 4.9: $\widehat Q$ in black and $\widehat P$ in grey.
\end{center}

Note that $\widehat P$ covers half the vertices of $\widehat Q$ and
the vertices in the complement $\widehat Q-\widehat P$ are, as
a subset of $\Z^2$, isometric to the vertices of $\widehat P$. The
isometry is given by the translation $(x_1,x_2) \to (x_1+2n,x_2)$.
Taking care to get the sign right, we see that 
$Q=P+P'$, where there is a bijection between the terms of
$P$ and the terms of $P'$, having the form
$$E(ax_1+bx_2) \to -E((a+2n)x_1+x_2).$$
When $x_1=\pi/2n$ we see that corresponding terms take on
the same value.  Hence $Q=2P$ when $x=\pi/2n$.  
But this means that $F(\frac{\pi}{2n},x_2)=0$.
\endproof

Equation \ref{veech1} follows from
Lemma \ref{elim}, Lemma \ref{final}, and Lemma \ref{yline}.
It only remains to prove Lemma \ref{final}.  The rest of
the chapter is devoted to this.

\subsection{Variation of Edgelength}

In proving Lemma \ref{final} we will connect
various vertices of the unfolding together by the
same sorts of paths we used in the proof of
Lemma \ref{elim}.  These paths mainly involve
the long edges, which all have unit length, but
sometimes they involve a short edge as well.
Even though we are evaluating the derivatives
of our defining functions at $V_n$, we still
need to understand how these short edges vary
in length for points nearby $V_n$.  In this
section, we will deal with this issue.

Suppose $T$ is a triangle with small angles
$x_1$ and $x_2$, normalized so that the long
side of $T$ has unit length.  Let
$l_j$ denote the length of the side of $T$
which is opposite the $x_j$ angle.  Of
course, $l_j$ depends on the parameters
$x_1$ and $x_2$.   When $x_1=x_2$ we
have $l_1=l_2$. When
$x_1=x_2=\pi/2n$ we have

\begin{equation}
\lambda:=l_1=l_2=\frac{\sin(x_1)}{\sin(x_1+x_2)}=
\frac{\sin(\pi/2n)}{\sin(\pi/n)}=
\frac{1}{2\cos(\pi/2n)}=
\frac{1}{\omega+\omega^{-1}}.
\end{equation}
As usual, $\omega=E(\pi/2n)$.
Here $\lambda$ is as in Equation \ref{elimX}.
Our calculations below require the quantities:

$$\lambda_1:=\frac{dl_2}{dx_2}\bigg|_{V_n}=
\frac{dl_1}{dx_1}\bigg|_{V_n}=\frac{\sin(x_2)}{\sin^2(x_1+x_2)}\bigg|_{V_n}=
\frac{2i\omega}{(\omega-\omega^{-1})(\omega+\omega^{-1})^2}$$
\begin{equation}
\label{diff}
\lambda_2:=\frac{dl_2}{dx_1}\bigg|_{V_n}=
\frac{dl_1}{dx_2}\bigg|_{V_n}=
\frac{-\cos(x_1+x_2)\sin(x_1)}{\sin^2(x_1+x_2)}\bigg|_{V_n}=
\frac{-i\omega(\omega^2+\omega^{-2})}{(\omega-\omega^{-1})(\omega+\omega^{-1})^2}.
\end{equation}

\subsection{Proof of Lemma \ref{final}}

We define
\begin{equation}
F_{\alpha,j}={\rm height\/}(\alpha_j)-{\rm height\/}(a_1); \hskip 30 pt
F_{\beta,j}={\rm height\/}(\beta_j)-{\rm height\/}(a_1).
\end{equation}
here $a_1=(1,L,O)^*$.
Again, we measure these heights when the $U_n$ is normalized
so that the long edges are unit length.  We have the
obvious equation
\begin{equation}
H_{ij}=F_{\alpha,i}-F_{\beta_j}.
\end{equation}
We will deduce Lemma \ref{final} from our computations
of $F_{\alpha,i}$ and $F_{\beta,j}$.

In our proof of Lemma \ref{elim} we constructed a path
from $b_1=(1,L,I)$ to and given point $p$.  The
first edge of this path joined $b_1$ to $a_1$. So,
the path we use to connect $a_1$ to $p$ is just
the same one we used above, except with the first
edge chopped off.  In describing our paths, we
let $Y_k$ denote the path made from the first
$k$ edges of the $3$-spine.  We let
$e_k^{\pm}$ denote the short edge such that
$$\widehat {e_k^{\pm}}=\widehat e_k+(\pm 1,0).$$
Here $e_k$ is the $k$th edge of the $3$-spine.
(The correspondence $e \to \widehat e$ is
discussed in detail in \S \ref{sect:background}.)
The three
edges $e_k^-,e_k,e_k^+$ correspond to $3$ consecutive
horizontal dots in Figure 4.8

With this notation, we have:
\begin{itemize}
\item The path connecting $a_1$ to $\alpha_1$ is $Y_{13}$.
\item The path connecting $a_1$ to $\alpha_2$ is 
$Y_5 \cup e_6^+$.  The short edge has type $2$.
\item The path connecting $a_1$ to $\alpha_3$ is
$Y_9 \cup e_{10}^+$. The short edge has type $2$.
\item The path connecting $a_1$ to $\beta_1$ is $Y_3$.
\item The path connecting $a_1$ to $\beta_2$ is
$e_1^-$.  This (short) edge has type $1$.
\item The path connecting $a_1$ to $\beta_3$ is
$Y_{15} \cup e_{16}^-$.  This edge has type $2$.
\end{itemize}

We will leave the details of our calculation to Mathematica,
but here we outline the main points.  In each case we
$\widetilde F=P \overline Q$, so that
$F={\rm Im\/}(\widetilde F)$.  By the product rule we have
\begin{equation}
\partial_j \widetilde F(V_n)=
P(V_n) \partial_j \overline Q(V_n)+\partial_j P(V_n) \overline Q(V_n).
\end{equation}
We evaluate $P$ and $Q$ using Equation \ref{elimX} (without the first term).
We evaluate $\partial_j P$ and $\partial_j Q$ essentially by
differentiating Equation \ref{elimX} (without the first term.)
We now explain how the differentiation works.   Our
calculations use the lists from Equation \ref{master}.

Let $R$ be one of the expressions we want to differentiate.  If
$Y_k$ appears in the definition of the path associated to $R$
then we see a contribution of
$$\sum_{i=1}^k L_j(i) \omega^{L(i)}$$
in the expression for $\partial_j Y(V_n)$.
In conjugating (for the case $R=Q$) we simply reverse the
signs of the list of numbers in $L$.  For instance
$$
\partial_1 \overline Q(V_n)=\sum_{i=1}^{20} L_1(i) \omega^{-L(i)}=$$
\begin{equation}
\frac{2i}{\omega^3}[(-4-4\omega^2-7\omega^4-3\omega^6)n+(10+14\omega^2+16\omega^4+8\omega^6)].
\end{equation}
If we see $e_k^{\pm}$ in our expression, and this edge has type $1$, then
we see a contribution of
\begin{equation}
\lambda (L_1(k)\pm 1)\omega^{L(k)\pm 1} + \lambda_1 \omega^{L(k) \pm 1}
\end{equation}
in the expression for $\partial_1 R(V_n)$ and a contribution of
\begin{equation}
\lambda (L_2(k)+0)\omega^{L(k)+1} + \lambda_2 \omega^{L(k) \pm 1}
\end{equation}
in the expressions for $\partial_2 R(V_n)$.
If $e_k^{\pm 1}$ has type $2$, then we see the
same contributions, but with $\lambda_1$ and $\lambda_2$ switched.

These are the ingredients for our calculations. 
We let $\widetilde H_{ij}=\widetilde F_{\alpha,i}-\widetilde F_{\beta,j}$.
When we compute these quantities using the expressions above,
we find that the result always has the form
\begin{equation}
\label{starform}
\partial_i\widetilde H_{jk}= \frac{f(n,\omega,\omega^{-1})}{\omega^a} \hskip 20 pt
{\rm or\/} \hskip 20 pt
\frac{f(n,\omega,\omega^{-1})}{\omega^a(\omega^2-\omega^{-2})}.
\end{equation}
Here $f$ is some polynomial in $\omega$, $\omega^{-1}$ and $n$ which
is linear in $n$. (This polynomial, and the exponent $a$, both depend
on the indices $i,j,k$.)  The first case occurs $5$ times and the second
case occurs $13$ times.

Since we only care about the sign of the imaginary part of
$\partial_i \widetilde H_{jk}$ we clear denominators by
multiplying the second form by the positive real expression
\begin{equation}
I(\omega^2-\omega^{-2}).
\end{equation}
Call the resulting expression $\widetilde L_{ijk}$.
When $\partial_i \widetilde H_{jk}$ has the first form
we simply set $\widetilde L_{ijk}=\partial_i \widetilde H_{jk}$.

In all cases we find that
\begin{equation}
\widetilde L_{ijk}=
\left[\begin{matrix}c_{-8} \cr c_{-6} \cr c_{-4} \cr c_{-2} \cr c_0 \cr c_2 \cr c_4 \cr c_6 \cr c_8 \end{matrix} \right] \cdot 
\left[\begin{matrix}\omega^{-8} \cr \omega^{-6} \cr \omega^{-4} \cr \omega^{-2} \cr \omega^0 \cr 
\omega^2 \cr \omega^4 \cr \omega^6 \cr \omega^8\end{matrix}\right]\to
\left[\begin{matrix}c_0-\overline c_0 \cr c_2-\overline c_{-2} \cr c_4-\overline c_{-4} \cr c_6-\overline c_{-6} \cr
c_8 -\overline c_{-8}\end{matrix}\right] \cdot
\left[\begin{matrix}1 \cr \omega^2 \cr \omega^4 \cr \omega^6 \cr \omega^8\end{matrix}\right]=
\left[\begin{matrix}a_0 \cr a_1 \cr a_2 \cr a_3 \cr a_4\end{matrix}\right] \cdot
\left[\begin{matrix}1 \cr \omega^2 \cr \omega^4 \cr \omega^6 \cr \omega^8\end{matrix}\right].
\end{equation}
This time $c_j$ has the form
\begin{equation}
c_j=c_{j0}+c_{j1}n; \hskip 30 pt
c_{0j},c_{j1} \in \Z[i].
\end{equation}
Again, the expression $x \to y$ means that the imaginary parts of these
two expressions are positive multiples of each other.

In listing the results of our calculations, we just write out
the coefficient vectors $(a_0,...,a_5)$.  Here are the
$9$ expressions for $-\widetilde L_{1jk}$:

\begin{equation}
\begin{matrix}
* & 120 i n&92 i n&40 i n&8 i n&0 \cr
 & 0&68 n\!-\!60&88 n\!-\!72&48 n\!-\!40&8 n\cr
 & 0&124 n\!-\!220&92 n\!-\!168&24 n\!-\!40&0\cr
 & 0&44 n\!-\!60&52 n\!-\!72&24 n\!-\!40&0\cr
 &  0&32 n\!-\!120&56 n\!-\!144&32 n\!-\!80&0\cr
 &  0&88 n\!-\!280&60 n\!-\!240&8 n\!-\!80&\!-\!8 n\cr
 &  0&140 n\!-\!176&128 n\!-\!168&56 n\!-\!84&6 n\!-\!12\cr
 &  0&128 n\!-\!236&132 n\!-\!240&64 n\!-\!124&6 n\!-\!12\cr
* &  i (224 n\!-\!520)&i (134 n\!-\!348)&i (40 n\!-\!124)&i (\!-\!2 n\!-\!12)&0\end{matrix}
\end{equation}
(We have deliberately listed $-\widetilde L_{ijk}$.)
Here are the $9$ expressions for $\widetilde L_{2jk}$.
\begin{equation}
\label{mat2}
\begin{matrix}
*&  0&0&0&0&0\cr
&  0&52 n\!-\!60&56 n\!-\!72&32 n\!-\!40&0\cr
&  0&88 n\!-\!176&92 n\!-\!184&44 n\!-\!84&0\cr
&  0&52 n\!-\!60&56 n\!-\!72&32 n\!-\!40&0\cr
&  0&104 n\!-\!120&112 n\!-\!144&64 n\!-\!80&0\cr
&  0&140 n\!-\!236&148 n\!-\!256&76 n\!-\!124&0\cr
&  0&116 n\!-\!220&88 n\!-\!168&16 n\!-\!40&0\cr
*&  i (216 n\!-\!360)&i (144 n\!-\!240)&i (48 n\!-\!80)&0&0\cr
*&  i (264 n\!-\!520)&i (180 n\!-\!352)&i (60 n\!-\!124)&0&0\end{matrix}
\end{equation}

The starred lines come from
the first form in Equation \ref{starform} and
the unstarred lines come from the second form.
The first line in Equation \ref{mat2} is
$\widetilde L_{211}$, the quantity we expect to vanish
in light of Lemma \ref{yline}.
For $n \geq 4$ the positive imaginary terms in the other 
starred lines
dominate the negative imaginary terms, and the
positive terms in the unstarred lines dominate the
negative terms.  Hence, with the exception of the
one line which represents the quantity $0$, all
the other lines represent quantities with positive
imaginary part.  (To be sure, we checked this
numerically for $n=4,...,20$.) 

 This
completes the proof of Lemma \ref{final},
and we are done.

\newpage

%% file: 5word.tex
\section{The Words for Theorem \ref{limit1}}
\label{s3}

Theorem \ref{limit1} is the most delicate of our
existence results.  Here we introduce the necessary
words.  In later chapters, we will analyze these
words, as we did for the proofs of Theorems
\ref{nonveech} and \ref{easy}.

\subsection{The Squarepaths}
\label{s3.1}

Theorem \ref{limit1} involves the words
$\{W_{nk}\}$ for $n=3,4,5...$
and $k=0,1,2...$.   In this chapter we introduce these
words and consider the corresponding unfoldings. It turns out that
$W_{nk}$ has length $24n+30k^2-68k-20$.  The shortest
word, $W_{30}$, has length $52$.  Rather than present 
$W_{nk}$ as a long string of digits, we will draw the
square path $\widehat Q_{nk}:=\widehat Q_3(W_{nk})$.
The path $\widehat Q_{nk}$ is not embedded,
but is the union of two embedded halves.
Reflection about a diagonal line swaps these
two halves. We will draw one half
of $\widehat Q_{nk}$, we well as the diagonal line.

$\widehat Q_{nk}$ is based an in $(n-1) \times n$ grid
of squares, which we call an $n$-{\it stamp\/}.  Each
square in the stamp has edge-length $2$, as in
Figure 2.3. 

Figure 5.1 shows $3$ representations of the word
$W_{30}$.  The leftmost figure shows the squarepath
$\widehat Q_{30}$.  This closed path is composed of $2$ halves
that are swapped by reflection in a certain diagonal
line of symmetry.  The middle figure shows one half
$\widehat Q_{30}$.  This is the representation we will
use in the other figures.  We prefer this representation
because the corresponding path is always embedded, and the full squarepath
can be recovered in a straightforward way.  Just reflect
and concatenate.  The right figure shows the corresponding
half of the hexpath $H_{30}$.

\begin{center}
\resizebox{!}{2.6in}{\includegraphics{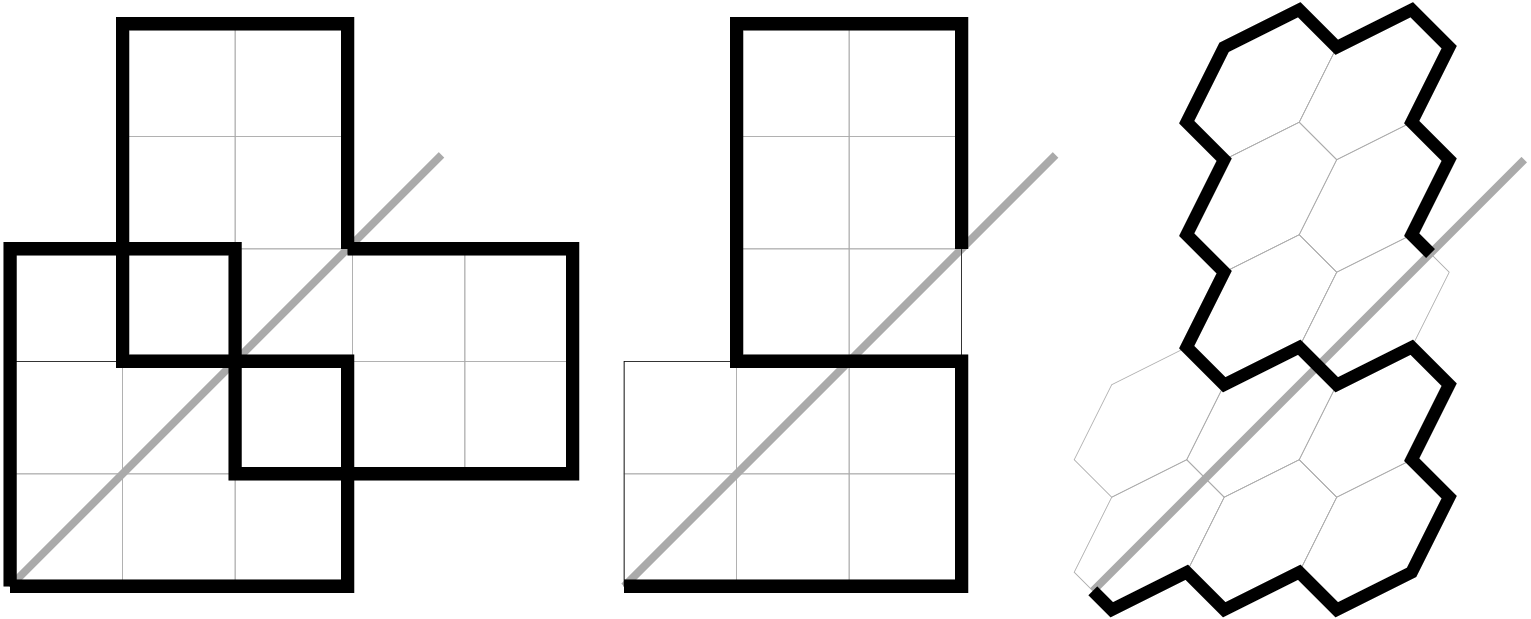}}
\newline
Figure 5.1: Half of $\widehat Q_{31}$.
\end{center}

Figure 5.1 is the beginning of an infinite pattern of
paths.  
Figure 5.2 shows the corresponding halves of
$\widehat Q_{3k}$ for $k=1,2,3$.   The small grid
of squares has been erased in Figure 5.2.

\begin{center}
\resizebox{!}{3in}{\includegraphics{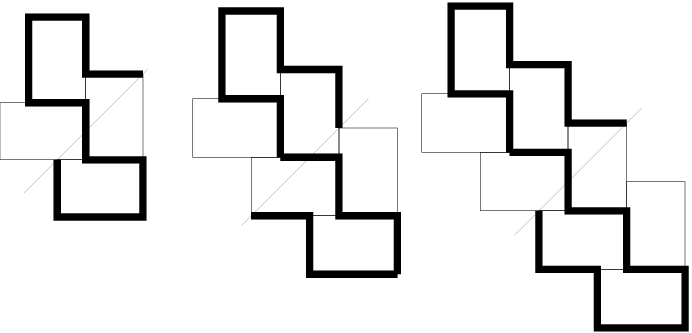}}
\newline
Figure 5.2: Half of $\widehat Q_{3k}$ for $k=1,2,3$.
\end{center}

The path $\widehat Q_{n+1,k}$ is obtained by increasing the
length of each edge of $\widehat Q_{nk}$ by $2$ units.
Figure 5.3 shows the left halves of
$\widehat Q_{31}$, $\widehat Q_{41}$, and $\widehat Q_{51}$.

\begin{center}
\resizebox{!}{2.6in}{\includegraphics{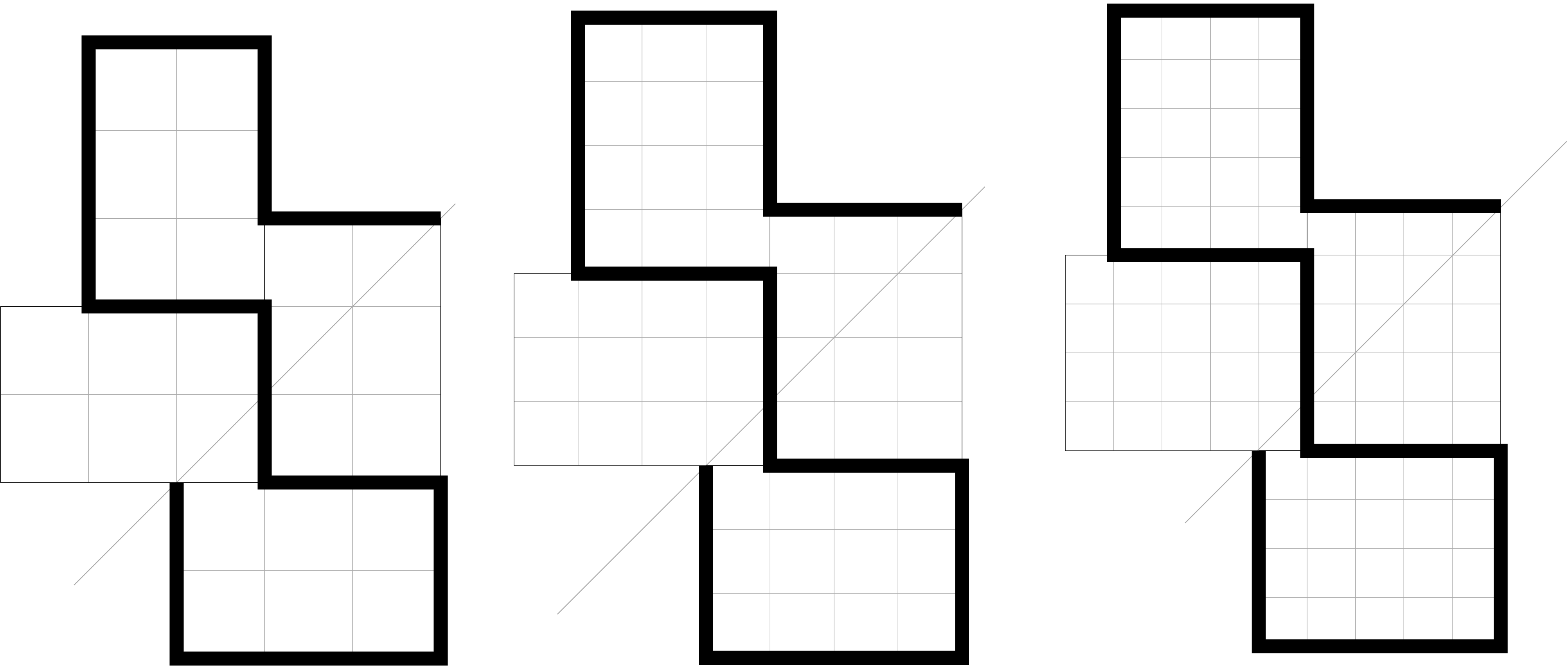}}
\newline
Figure 5.3: Half of the rectilinear paths for $\widehat Q_{30}$ and $\widehat Q_{40}$ and $\widehat Q_{50}$
\end{center}

\subsection{The Unfoldings}
\label{s3.2}

We will see that $U_{nk}$ consists of $4$ ``strips'',
attached along $4$ ``hinges''.  Figure 5.4 shows
this structure.
  
\begin{center}
\includegraphics[height=2.8in]{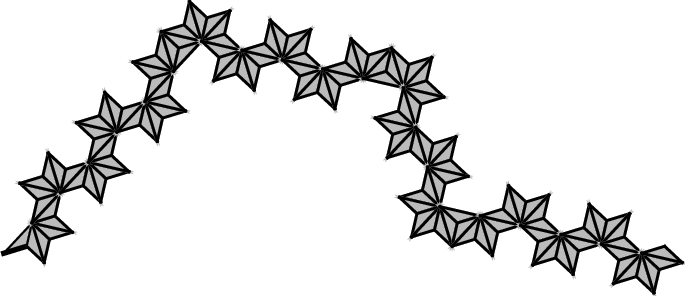}
Figure 5.4: $U(W_{42},V_3)$.
\end{center}

When we change the point relative to which we
unfold, the strips do not change much and
the hinges open and close, so to speak.
For points in the orbit tiles, the hinges
adjust so that the whole unfolding is
practically a straight line.  In Figure 4.4,
we have chosen a point that is far from the
relevant tile.  Our remaining pictures show
the unfoldings for points actually in the
relevant orbit tile.

The strips are essentially composed of units that we
call {\it blocks\/}.  The
left hand side of Figure 5.5 shows what we call a {\it block\/}.
In general, a $k$-block is defined to be $k$ blocks lined
up in sequence.  The right hand side of Figure 5.5 shows a $2$-block.
The triangles in a $k$-block all have the same shape, and this
shape depends on the point in parameter space of interest to us.
If we glue the opposite sides of a block together we get a
space which is naturally the union of two $2$-darts.

\begin{center}
\includegraphics{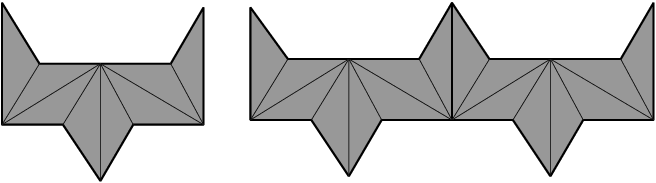}
\newline
Figure 5.5: A $1$-block and a $2$-block
\end{center}

The strips are essentially
composed of blocks.  Once we describe $U_{30}$, we
will describe $U_{nk}$ as a modification which amounts to
changing the combinatorial structure of each strip.
Figure 5.6 shows $U(W_{30},V_3)$.  

\begin{center}
\includegraphics[width=6in]{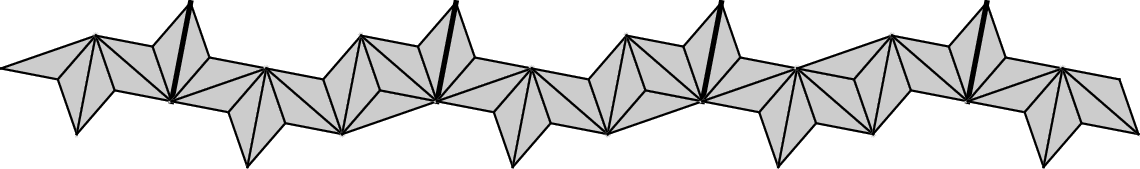}
\newline
Figure 5.6: $U(W_{30},V_3)$.
\end{center}

Note that $U(W_{30},V_3)$ has $12$ long edges which are
all parallel and nearly vertical.   We have highlighted
$4$ of these edges.  The unfolding $U_{3k}$ is obtained
by cutting $U_{30}$ open along each of the $4$ highlighted
edges and inserting a $k$-block.   Figure 5.7 shows
$U(W_{31},V_3)$.  The pattern continues in the obvious way.
In describing our surgery, we have used the geometry of
$U(W_{30},V_3)$ to highlight $4$ particular edges along
which we cut.  However, this surgery has a combinatorial
meaning for any parameter.

\begin{center}
\includegraphics[width=6in]{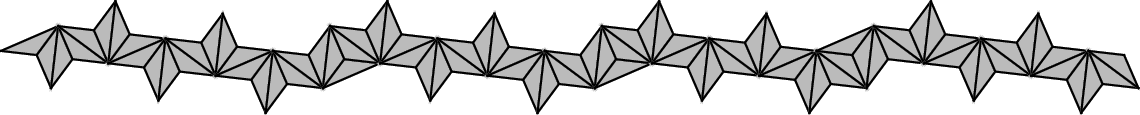}
Figure 5.7: $U(W_{31},V_3)$.
\end{center}

We obtain $U_{n,k}$ from $U_{3k}$ by replacing each
maximal $m$-dart with a maximal $m'$-dart,
where
$$m'=m+(n-3).$$
This fits exactly with our description of
$\widehat Q_{nk}$ as being obtained from
$\widehat Q_{3k}$ by lengthening each edge
by $2(n-3)$ units.
  
\begin{center}
\includegraphics[width=6.25in]{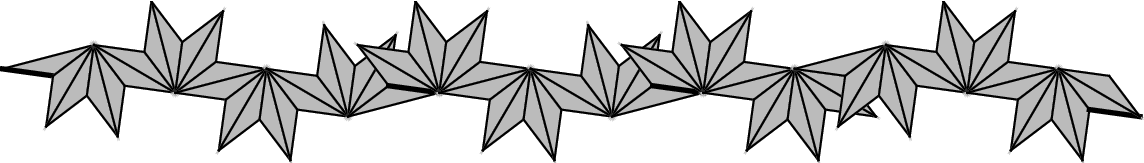}
\newline
Figure 5.8: $U(W_{40},V_4)$.
\end{center}

\begin{center}
\includegraphics[width=6.25in]{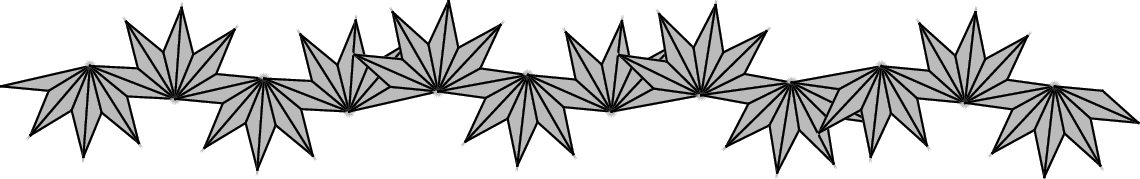}
\newline
Figure 5.9: $U(W_{50},V_5)$.
\end{center}

We end this section with a computation.  The
formulas in the next result will be useful
when we make explicit computations in \S \ref{sect:rescaling} and \S \ref{sect:pivot}.

\begin{lemma}
\label{slope1}
Let $n$ be fixed and let
$e_k$ be the leftmost edge of $U(W_{nk},V_n)$.
As $k \to \infty$, the slope of $e_k$ converges to $0$.
\end{lemma}

\startproof
We normalize so that $e_k$ has unit length.
If we trim off the portions of the darts 
from the set $U:=U(W_{nk},V_n)$ 
we see that the resulting set is the union of
$4$ parallel annuli attached along $4$ edges.
We compute by elementary trigonometry
and induction that the $4$ strips have total
length
\begin{equation}
\Psi_1+\Psi^{\#}k; \hskip 30 pt
\Psi_1=12(1+\cos(\frac{\pi}{n}));
\hskip 15 pt
\Psi^{\#}=8(1+\cos(\frac{\pi}{n})).
\end{equation}
Each annulus has width
\begin{equation}
\Psi_2=\sin(\frac{\pi}{n})
\end{equation}
If we rotate $U$ so that the first edge is horizontal
then the holonomy has coordinates
\begin{equation}
\Psi_1+i4\Psi_2+\Psi^{\#}k.
\end{equation}
The line determined by this complex number converges
to a horizontal line as $k \to \infty$.  Thus,
if we rotate $U$ so that the holonomy is
horizontal, then the slope of the first edge
tends to $0$ as $k \to \infty$.
\endproof

\noindent
{\bf Remark:\/}
In \S \ref{sect:pivot} we will give a more explicit and
combinatorial derivation of the formulas in
Lemma \ref{slope1}.

\subsection{The Pivot Region}
\label{s3.3}

Here we isolate $4$ basic features of the unfolding
$U(W_{nk},V_n)$.
\begin{itemize}
\item There is a family of $12+8k$ parallel and nearly vertical edges.
We call these edges {\it quasi-vertical\/}, or QV for short.
\item There is a family of $24+16k$ parallel and nearly horizontal
edges.  We call these edges {\it quasi-horizontal\/}, or QH for short.
\item Each QV edge is {\it flanked\/} by two QH edges, in the sense
that reflection in this QV edge swaps the two QH edges flanking it.
\item There are exactly $4$ QH edges which connect top to bottom
vertices.  We call these edges the {\it hinges\/}.
\end{itemize}
These facts are all established inductively.  They hold true for
the parameter $(3,0)$, and then we check easily that they remain
true when we perform one of the surgeries described above.

We have distinguished the above edges just for the unfoldings
attached to specific parameters. However, we extend our definitions
of QH and QV, using continuity, to the unfoldings attached to
any point of parameter space.  Of course, for points remote to
the regions of interest to us, the QH edges need not be close
to horizontal and the QV edges need not be close to vertical.
Moreover, these edges need not be parallel to each other
at other parameters.  To talk precisely about the situation
we make the following definition:
\newline
\newline
{\bf Definition:\/} 
Let $A_n$ denote the region $(x_1,x_2) \subset \Delta$
such that $$x_j \in [\pi/2n,\pi/(2n-2)]; \hskip 30 pt j=1,2.$$  The points
$V_n$ and $V_{n-1}$ are two opposite corners of the
little square defined by these conditions.
Let $R'_{nk} \subset A_n$ denote the set of
points such that all the QH edges have negative slope.
Let $R_{nk}$ denote the path 
connected component of $R'_{nk}$ which
contains $V_n$.  We call $R_{nk}$ the
{\it pivot region\/}.

\begin{lemma}
For any point in $R_{nk}$ the QV edges all have positive slope.
\end{lemma}

\startproof
As we pointed out above,
each QV edge $V$ is flanked by two QH edges
$H_1$ and $H_2$.  That is, reflection in $V$ swaps
$H_1$ and $H_2$.  This is a property that holds
for all parameters:  It is a combinatorial symmetry.
Now, let $X \in R_{nk}$ be some
point.  We consider what happens as we vary the
parameter continuously from $V_n$ to $X$,
staying inside $R_{nk}$.  If the slope of
$V$ changes from positive to negative then
$V$ must be either vertical or horizontal at some point. 
But then it is impossible for $H_1$ and $H_2$
to both have negative slope at this point
and still flank $V$.
This is a contradiction. 
\endproof

Referring to \S \ref{sect:dart}, we note that a vertex of $U_{nk}$
is superior if and only if it is incident to a QH edge.
This fact is seen by inspection for $U_{30}$, and then
is unchanged by any of the surgeries we perform.
We call a superior vertex a {\it pivot\/} if it is
incident to one of the pivot edges.
There are $4$ top pivots and $4$ botton pivots.

\begin{lemma}
\label{ruleout}
Let $X \in R_{nk}$ be any point. If all the top
pivots lie above all the bottom pivots of
$U(W_{nk},X)$ then $X \in O(W_{nk})$.
\end{lemma}

\startproof
Let $v$ be a top superior vertex of $U_{nk}$ which is
not a pivot.  Inspecting our unfoldings, we see
that $v$ has one of two properties:
\begin{itemize}
\item $v$ is the left vertex of a QH edge which is
not a hinge.
\item $v$ lies to the left of another superior vertex $v'$,
and reflection in a QV edge swaps $v$ and $v'$.
\end{itemize}
This property is easily seen, by inspection, for
$U_{30}$, and our surgery operations do not destroy
this property.
In either of the above cases, our conditions on
the slopes of the QV and QH edges forces $v$ to lie
above $v'$.   Since this works for all superior vertices
which are not pivots, we see that only a pivot can
be the lowest top superior vertex.  Likewise,
only a pivot can be a bottom superior vertices.
Hence, by hypothesis, all the top superior
vertices lie above all the bottom superior vertices.
Hence $X \in O(W_{nk})$ by the Dart Lemma.
\endproof

\subsection{The Structure of the Quasi-Horizontal Edges}
\label{QH}
\label{s3.4}

Say that a {\it QH point\/} in $\Z^2$ is a point $\widehat e$
which corresponds to a QH edge.   Here we
describe the pattern of QH points associated to
our unfoldings.

We use the McBilliards labeling convention that the
leftmost edge $e_0$ of $U_{nk}$, which happens to be a QH edge,
corresponds to $(0,0) \in \Z^2$.  If $e$ is any other
QH edge, then $e$ is parallel to $e_0$ at the point
$V_n=(\pi/2n,\pi/2n)$.   But the angle between $e$ and $e_0$
is given by
$\widehat e \cdot (\pi/2n,\pi/2n).$
This angle must be an integer multiple of $\pi$.   Hence
\begin{equation}
\label{modwidth}
\widehat e_1+\widehat e_2 \equiv 0 \hskip 5 pt
{\rm mod\/} \hskip 5 pt 2n.
\end{equation}
The map $(x,y) \to x+y$ maps the hexpath $H_{nk}$ to a subset
of $\Z$ having diameter less than $(4+\frac{1}{2})n$.  Hence
Equation \ref{modwidth} forces the QH points to lie
along at most $4$ lines of slope $-1$ in $\Z^2$.  Hence
the QH edges fall into at most $4$ pseudo-parallel
families, in the sense of \S \ref{pp}.  

After some trial and error we figured out how to draw the
QH points.  We find that these points fall into
exactly 4 pseudo-parallel families, and we compute
the extreme points of these families as follows:
Setting
\begin{equation}
Z_{nk}=(2n-2)(k+2)-1,
\end{equation}
the coordinates for the extreme points, in each family, are given by
\begin{enumerate}
\item $(-3,-4n+3)$ and $(-4n+3,-3)+(Z_{nk},-Z_{nk})$.
\item $(-2,-2n+2)$ and $(-2n+2,-2)+(Z_{nk},-Z_{nk})$.
\item $(0,0)$ and $(0,0)+(Z_{nk},-Z_{nk})$.
\item $(2n-2,2)$ and $(2,2n-2)+(Z_{nk},-Z_{nk})$.
\end{enumerate}
Actually, we don't care so much about these formulas.
The main feature of the QH points we use is that
the coordinates of the northwest extreme endpoints$-$the
first ones listed in each line above$-$are independent of
$k$.  This property, together with symmetry, will tell us
everything we want to know about these edges.

Figure 4,10 shows the first few examples of
the QH points in $\Z^2$.
The $8$ larger dots are the extreme points.
$4$ of these dots are black and the other $4$ are grey.
The $4$ grey dots in these figures correspond to the hinges.
The northwest grey dot is labeled $(0,0)$.
The southeast grey dot is $(A,-A)$.
The black path is the path $\widehat Q$, discussed
in \S \ref{s2.4}, which corresponds to the $3$-spine of
the unfolding.  This path is obtained by
doubling the rectilinear paths discussed in
\S \ref{s3.1}.

\begin{center}
\includegraphics[height=6.5in]{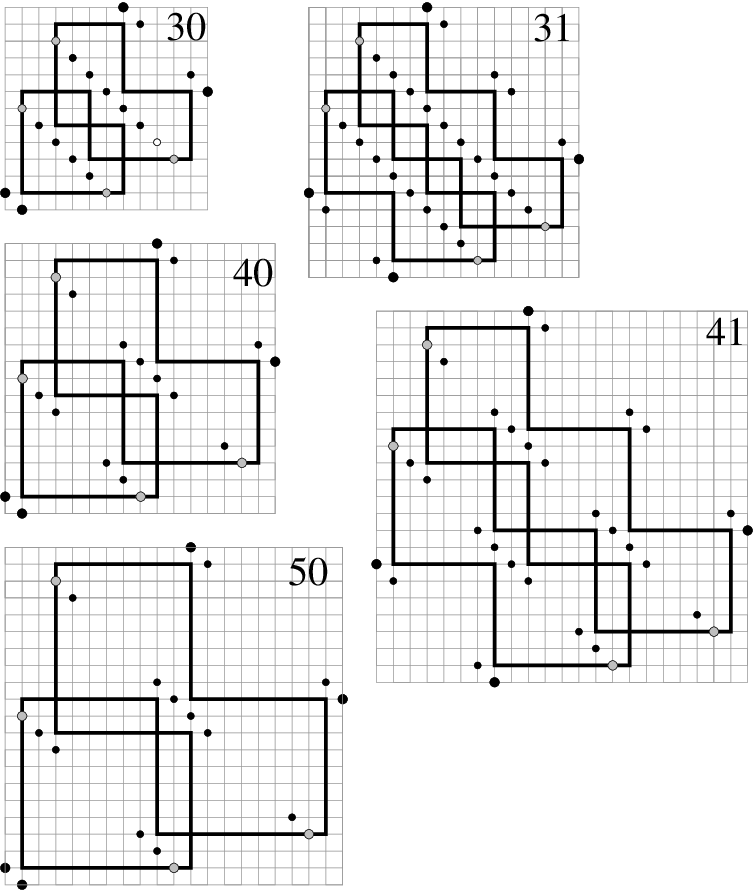}
\newline
Figure 5.10: Points for the QH families.
\end{center}

\newpage

%% file: 6qrt.tex
\section{The Quadratic Rescaling Theorem}
\label{sect:rescaling}

\subsection{Overview}

We are interested in studying infinite sequences $\{W_{nk}\}$
of words introduced in the last chapter. 
We will hold $n$ fixed and set $W_k=W_{nk}$. 
We wish to understand the asymptotic shape of the orbit tiles
$O(W_k)$ as $k \to \infty$.
Each $O(W_k)$ is a piecewise analytic
polygon, whose sides are given as the $0$-level sets
of analytic functions.   It turns out that
$O(W_k)$ has a uniformly bounded
number of sides, independent of both $n$ and $k$,
and we will be able to
make sense of the notion of a {\it side\/} of
$O(W_k)$ which is independent of $k$. This
allows us to group the various defining functions
involved into {\it families\/}.
Our analysis is done one function-family at a time.

As we saw in \S \ref{sect:background}, the function $F_k$ has
the special form
\begin{equation}
\label{pairing}
F_k(X)= {\rm Im\/}(P_k(X-X_0)\overline Q_k(X-X_0)),
\end{equation}
(The only difference between the set-up here and in
\S \ref{sect:background} is that we take special care to translate $F$ so
that $(0,0)$, rather than $X_0$, is the main point
of interest to us in the domain.)
Here $X_0$ is the point in parameter space to which
the orbit tiles converge.  Summarizing
the discussion in \S \ref{sect:background}, 
$P_k$ is the development image of
a certain saddle connection associated to $W_k$,
and $Q_k$ is the holonomy of the unfolding.

We want to place conditions on
$\{P_k\}$ and $\{Q_k\}$ so that the rescaled
functions $\{G_k\}$ converges to a linear
map (whose formula we can compute explicitly).
Here
\begin{equation}
\label{rescale}
G_k(X)=F_k(Xk^{-2}).
\end{equation}
The conditions we place on $\{P_k\}$ and $\{Q_k\}$
have to do with the growth patterns of the
supports of their Fourier transforms.
A few glances at the figures in \S \ref{s3} should
be enough to convince the reader that the conditions
we discuss are satisfied, in particular, by
the functions associated to the words introduced
in \S \ref{s3}.
We will state the Quadratic Rescaling Theorem
in the next section and then spend the rest of
the chapter proving it.

As in \S \ref{sect:background} we frequently let $E(x)=\exp(ix)$.

\subsection{The Main Result}

All our constructions are based on a translation
$T: \Z^2 \to \Z^2$ and a point $X_0 \in \pi \Q^2$.
More specifically, we have
\begin{equation}
\label{basics}
T(x_1,x_2)=(x_1+M_1,x_2+M_2); \hskip 30 pt
X_0=2 \pi \bigg(\frac{p_1}{q_1},\frac{p_2}{q_2}\bigg).
\end{equation}
There is a natural homomorphism associated to $X_0$:
\begin{equation}
\label{hom}
\phi(x_1,x_2)=\bigg[\frac{p_1q_2x_1+p_2q_1x_2}{{\rm G.C.D.\/}(q_1,q_2)}\bigg] \in \Z/N; \hskip 30 pt
N=q_1q_2/{\rm G.C.D.\/}(q_1,q_2).
\end{equation}
Here ${\rm G.C.D.\/}$ stands for {\it greatest common divisor\/}.
We require that $T$ is compatible with $X_0$ in the sense that
$\phi(M_1,M_2)=0$.
\newline
\newline
\noindent
{\bf Remark:\/}
In this paper we will have $n=3,4,5...$ and
\begin{equation}
\label{examples}
T(x_1,x_2)=(x_1+(2n-2),x_2-(2n-2)); \hskip 30 pt
X_0=2 \pi \bigg(\frac{1}{4n},\frac{1}{4n}\bigg).
\end{equation}
In this case, we have
$\phi(x_1,x_2)=[x_1+x_2] \in \Z/4n$, and $\phi(M_1,M_2)=0$.
\newline

Going back to the general case,
we let $\{R_k\}$ stand for either the
sequence $\{P_k\}$ or $\{Q_k\}$.   According to the
theory developed in \S \ref{sect:background}, we can write
\begin{equation}
\label{fourier}
R_k(X)=\sum_{V \in \Z^2} \widehat R_k(V) E(X \cdot V).
\end{equation}
Here $\widehat R_k: \Z^2 \to \Z$ is the Fourier transform of $\R_k$.
We say that $\widehat R_k$ has {\it linear growth\/} if there
is some map $\widehat R^{\#}: \Z^2 \to \Z$ such that
\begin{equation}
\widehat R_{k+1}=\widehat R_k + \widehat R^{\#} \circ T^k; \hskip 30 pt 
k=0,1,2...
\end{equation}
Here $T$ is the translation above.  
We call $\widehat R^{\#}$ the {\it growth generator\/} for
$\{R_k\}$.
We require that both $\widehat R^{\#}$ and $\widehat R_0$ are supported
on finitely many points of $\Z^2$.  Intuitively, the support of
$\widehat R_k$ grows linearly along the fibers of the homomorphism
$\phi$.   As the notation suggestions, define
\begin{equation}
R^{\#}_k(X)=\sum_{V \in \Z^2} \widehat R^{\#}_k(V) E(X \cdot V).
\end{equation}

Supposing that both $\{\widehat P_k\}$ and $\{\widehat Q_k\}$ have
linear growth, we define
$$
\delta={\rm det\/}\left[\begin{matrix}
 P^{\#}(X_0) & P_0(X_0) \cr
 Q^{\#}(X_0) & Q_0(X_0)\end{matrix}\right]; 
$$
$$
\delta_j={\rm det\/}\left[\begin{matrix}
 P^{\#}(X_0) & \partial_jP^{\#}(X_0) \cr
Q^{\#}(X_0) & \partial_jQ^{\#}(X_0)\end{matrix}\right];
$$
\begin{equation}
\label{dets}
\Delta_j=\frac{M_j \delta}{2} + {\rm Im\/}(\delta_j).
\end{equation}
Here $\partial_j$ is the partial derivative with respect to $x_j$.

The rest of this chapter is devoted to proving:

\begin{theorem}[Quadratic Rescaling]
Suppose that $\{\widehat P_k\}$ and $\{\widehat Q_k\}$ have
linear growth with respect to $T$, and the quantities
$P^{\#}(X_0)$ and $Q^{\#}(X_0)$ and $\delta$ are all real.
Then $\{G_k\}$ converges in the $C^{\infty}$-topology to $G$, whose
equation is given by
$$G(x_1.x_2)=F_0(0,0)-\Delta_1x_1-\Delta_2x_2.$$
\end{theorem}

\noindent
{\bf Remarks:\/} \newline
The $C^{\infty}$ convergence means that each partial derivative
of $G_k$ converges, uniformly on compact subsets, to the corresponding
partial derivative of $G$. 

\subsection{Quadratic Growth Conditions}

Let $\cal A$ denote the set of all globally defined
and analytic complex-valued functions on $\R^2$.
Given a {\it multi-index\/} $I=(i_1,i_2)$ we define
\begin{equation}
X^I=x_1^{i_1}x_2^{i_2}; \hskip 20 pt
|I|=i_1+i_2.
\end{equation}
Given an infinite sequence $\{F_k\} \in \cal A$ we can 
write out the power series expansions
\begin{equation}
F_k(X)=\sum_{I} C_{k,I}X^I
\end{equation}
We say that $\{F_k\}$ forms a {\it quadratic growth family\/} if
$\{F_k(0,0)\}$ is a constant sequence we have the following
finite limits for some $\epsilon>0$.
\begin{equation}
\label{cg}
\lim_{k \to \infty} C_{k,(1,0)}k^{-2}=C_1;
\hskip 20 pt
\lim_{k \to \infty} C_{k,(0,1)}k^{-2}=C_2;
\hskip 20 pt
\lim_{k \to \infty} 
\sum_{|I| \geq 2} |C_{k,I}|\ k^{(-2+\epsilon)|I|}=0.
\end{equation}
Note that $\epsilon$ only enters into the third equation.

Recall that $\{G_k\}$ is the rescaled version of $\{F_k\}$, as in
Equation \ref{rescale}. 

\begin{lemma}[Convergence]
Suppose that $\{F_k\}$ is a quadratic growth family. Then
$\{G_k\}$ converges in the $C^{\infty}$ topology
to the linear function $G$, whose formula is given by
$G(x_1,x_2)=F_0(0,0)+C_1x_1+C_2x_2$.
\end{lemma}

\startproof
From the chain rule, we get the following series expansion:
\begin{equation}
G_k(X)=\sum_I C_{k,I}\ k^{-2|I|}\ X^I.
\end{equation}
Consider the difference
$$\widehat G(X)-G_k(X)=L_k(X)+R_k(X); \hskip 30 pt
R_k(X)=\sum_{|I| \geq 2} C_{k,I}\ k^{-2|I|}\ X^I.
$$
Here $L_k(x)$ is a linear function whose coefficients vanish as
$k \to \infty$, and $R_k$ is everything else.  It suffices to
to show that $R_k$ and all its derivatives tend to $0$ uniformly
on compact subsets.

Let $\partial$ stand for some partial derivative and let
$\Omega$ be some big constant.  Suppose $X=(x_1,x_2)$ is 
such that $|x_j| \leq \Omega$ for
$j=1,2$.  There is some constant $N$, depending on $\partial$, such that
$$|\partial R_k(X)| \leq \sum_{|I| \geq 2} |I|^N\ |C_{k,I}|\ k^{-2|I|}\ \Omega^{|I|}=
\sum_{|I| \geq 2} \{|I|^N (\Omega/k^{\epsilon})^{|I|}\}\ |C_{k,I}|k^{(-2+\epsilon)|I|}.$$
For
$$k>\bigg(\Omega N^N\bigg)^{1/\epsilon}$$ 
the term in braces is less than $1$.  Hence
$$|\partial R_k(X)| \leq \sum_{|I| \geq 2}  |C_{k,I}|\ k^{(-2+\epsilon)|I|}.$$
By hypothesis, this last sum tends to $0$ as $k \to \infty$.
\endproof

\subsection{A Fact about the Fourier Transform}

Now we begin to use the information about
the sequences $\{P_k\}$ and $\{Q_k\}$ to establish
the conditions on $\{F_k\}$ discussed in the
Convergence Lemma above.
Let $X_0$ be as in Equation \ref{basics} and let
$\phi$ be the associated homomorphism, given
in Equation \ref{hom}.  In particular, the value
$N$ is given by Equation \ref{hom}.
Consider a function of the form
\begin{equation}
\label{fourier2}
R(X)=\sum_{V \in \Z^2} \widehat R(V) E(X \cdot V).
\end{equation}

Choosing any residue class $k \in \Z/N$ we
define the {\it modular transform\/}:
\begin{equation}
R_{\phi}(k)=\sum_{\phi^{-1}(k)} \widehat R(V).
\end{equation}
In all cases of interest to us, the sum in
Equation \ref{fourier2} is a finite sum.
This sum defines a function
$R_{\phi}: \Z/N \to \C$.

\begin{lemma}[Modular Transform]
With the notation as above, we have
$$R(X_0)=\sum_{j=1}^N R_{\phi}(j)E(2\pi j/N).$$
\end{lemma}

\startproof
Let 
$$N=\frac{q_1q_2}{D}; \hskip 30 pt D={\rm G.C.D.\/}(q_1,q_2).$$
We can write
$R(X_0)=\sum_{j=1}^N R_j$, where
$$R_j=\sum_{(x_1,x_2) \in \phi^{-1}(j)} 
\widehat R(V) E(\frac{2 \pi p_1x_1}{q_1}+\frac{2 \pi p_2x_2}{q_2})=$$
$$\sum_{(x_1,x_2) \in \phi^{-1}(j)}
\widehat R(V)E\bigg(\frac{2 \pi}{N}\bigg(\frac{x_1p_1q_2}{D}+\frac{x_2p_2q_1}{D}\bigg)\bigg)=$$
$$\sum_{(x_1,x_2) \in \phi^{-1}(j)}\widehat R(V) E \bigg(\frac{2\pi j}{N}\bigg)=
R_{\phi}(j) E(2\pi j/N).$$
Summing over $j$ we get the result.
\endproof

\subsection{Growth Formulas}

Let $X_0$ and $T$ be as in Equation \ref{basics}.
In particular, recall that $T$ represents translation
by the vector $(M_1,M_2) \in \Z^2$.
Let $\cal V$ denote the set of sequences of the
form $\{R_k\}$ which have $(M_1,M_2)$-linear growth.
We want to be clear that each individual {\it element\/}
of $\cal V$ is a sequence of functions.  $\cal V$ is
a vector space.  The vector
space laws on $\cal V$ are given by componentwise
scaling and addition.  That is
\begin{equation}
a \cdot \{R_k\}=\{aR_k\}; \hskip 30 pt
\{R_k\}+\{R'_k\}=\{R_k+R'_k\}.
\end{equation}

Let ${\cal V\/}_0$ denote the subspace consisting
of elements $\{R_k\}$ with $R_0=0$.  There is
a natural projection ${\cal V\/} \to {\cal V\/}_0$.
The sequence $\{R_k\}$ is mapped to
$\{S_k\}$ where $S_k=R_k-R_0$.  We call
$\{S_k\}$ the {\it pure projection\/} of
$\{R_k\}$.
We say that an element $\{R_k\}$ of ${\cal V\/}_0$ is {\it simple\/}
if its growth generator $\widehat R^{\#}$ is the indicator function for
a single lattice point.  That is, there is some integral
point $A$ such that $\widehat R^{\#}(X)=1$ iff $X=A$ and 
$\widehat R^{\#}(X)=0$
otherwise.  The simple elements of ${\cal V\/}_0$ form
a basis for ${\cal V\/}_0$. 

\begin{lemma}
\label{leading}
Let $\{R_k\}$ be any element
of ${\cal V\/}$ and let $I=(i_1,i_2)$ be any multi-index.
Then
\begin{equation}
\label{growth}
D_IR_k(X_0)=R^{\#}(X_0) \times \frac{i^{|I|} M_1^{i_1}M_2^{i_2}}{|I|+1} \times k^{|I|+1}+
O(k^{|I|}).
\end{equation}
\end{lemma}

\startproof
If $\{S_k\}$ is the pure projection of $\{R_k\}$ then
$$|D_IS_k(X_0)|-|D_IR_k(X_0)| = O(1)$$
Thus, it suffices to prove this lemma for
elements of ${\cal V\/}_0$.  
Given the scaling and additivity properties of
Equation \ref{growth}, it suffices
to establish Equation \ref{growth} for
the simple elements of ${\cal V\/}_0$.

Suppose $\widehat R^{\#}$ is
the indicator function for $(a_1,a_2) \in \Z^2$.
Then
\begin{equation}
\label{indicator}
\widehat{D_IR_k}(x_1,x_2)=i^{|I|} x_1^{i_1}x_2^{i_2} 
 \hskip 15 pt
\Longleftrightarrow \hskip 15 pt (x_1,x_2) \in 
\bigcup_{j=0}^{k-1} (a_1+jM_1,b_1+jM_2),
\end{equation}
and otherwise this function vanishes.
Let $\beta=\phi(a_1,a_2) \in \Z/N$. Note
that $\beta=R^{\#}(X_0)$ for simple elements.
All the points in Equation \ref{indicator}
lie in the same fiber of $\phi$, namely
$\phi^{-1}(\beta)$.  From the Modular
Transform Lemma we have
$$
D_IR_k(X_0)=i^{|I|}\beta \times \sum_{j=0}^{k-1}
(a_1+jM_1)^{i_1}(a_2+jM_2)^{i_2}=$$
$$i^{|I|}\beta M_1^{i_1}M_2^{i_2} \times \sum_{j=0}^{k-1} j^{|I|}+O(k^{|I|})=$$
\begin{equation}
\label{derive}
R^{\#}(X_0) \times \frac{i^{|I|} M_1^{i_1}M_2^{i_2}}{|I|+1} k^{|I|+1} + O(k^{|I|}).
\end{equation}
This completes the proof.
\endproof

One useful special case of Lemma \ref{leading}, stated more
precisely, is:
\begin{equation}
\label{growth0}
R_k(X_0)=R^{\#}(X_0)k + R_0(X_0).
\end{equation}

\noindent
{\bf Remark:\/} We can relate Equation \ref{growth0}
to Lemma \ref{slope1} as follows.  In our examples
we have $X_0=V_n$, the Veech Point.
If $R$ is the holonomy
function (the $Q$-function)
associated to either the $1$-spine or
$2$-spine of the unfolding $U_{nk}$, then
$R_0(X_0)=\Psi_1+i4\Psi_2$ and $R^{\#}(X_0)=\Psi^{\#}$.
\newline

\begin{lemma}
\label{order}
Let $\{R_k\}$ be any element
of ${\cal V\/}$. Then
for $j=1,2$,
\begin{equation}
\label{growth1}
\partial_j R_k(X_0)=
\bigg(\frac{iM_j}{2}R^{\#}(X_0)\bigg)k^2+
\bigg(\partial_j R^{\#}(X_0)-\frac{iM_j}{2}R^{\#}(X_0)\bigg) k
 + {\rm const.\/}
\end{equation}
\end{lemma}

\startproof
Let $\{S_k\}$ be the pure projection of $\{R_k\}$.
Assuming for the moment that $\{S_k\}$ is a simple element,
Equation \ref{derive}, applied to $I=(1,0)$, gives us
$$\partial_j S_k(X_0)-\frac{iM_1}{2} R^{\#}(X_0)k^2 =
\partial_j S_k(X_0)+\frac{iM_1}{2} \beta k^2 =$$
$$\frac{i\beta M_1}{2}k+a_1 k=
\bigg(\frac{iM_j}{2}S^{\#}(X_0)k+\partial_j S^{\#}(X_0)\bigg) k.$$
Hence Equation \ref{growth1} holds, with zero constant term, for
the simple elements.   Both sides of Equation \ref{growth} (with
zero constant term) can be interpreted as homomorphisms from
$\cal V$ into $\C$.  Hence, Equation \ref{growth1} holds,
with zero constant term, for all
elements of ${\cal V\/}_0$.
Finally we note that $\partial_j S_k(X_0)$ and
$\partial_j R_k(X_0)$ differ by a constant.
\endproof

\subsection{Consequences of the Growth Formulas}

Now we assume that $\{P_k\}$ and $\{Q_k\}$ and $\{F_k\}$ are
all as in the Quadratic Rescaling Theorem.

\begin{lemma}
\label{constantseq}
$\{F_k(0,0)\}$ is a constant sequence.
\end{lemma}

\startproof
Using Equation \ref{growth0} we compute
$F_k(0,0)=Xk^2+Yk+F_0(0,0),$
where
$$X={\rm Im\/}(P^{\#}(X_0)\overline{Q^{\#}(X_0)})=0;$$
$$Y={\rm Im\/}[P_0(X_0)\overline{Q^{\#}(X_0)}+ P^{\#}(X_0)\overline{Q_0(X_0)}]={\rm Im\/}(\delta)=0.$$
This completes the proof.
\endproof

\begin{lemma}
\label{gradseq}
$\partial_j F_k(0,0)=-\Delta_jk^2+O(k).$
\end{lemma}

\startproof
For easy reference we repeat Equations \ref{growth0} and \ref{growth1}.
For ease of notation we write $\rho=\rho(X_0)$, understanding that all
our functions $\rho$ are evaluated at $X_0$ unless we explicitly
indicate otherwise.  We will also set $M=M_j$ and $\partial=\partial_j$.
$$
R_k=R^{\#}k + R_0.
$$
$$
\partial R_k=
\bigg(\frac{iM}{2}R^{\#}\bigg)k^2+
\bigg(\partial R^{\#}-\frac{iM}{2}R^{\#}\bigg) k
 + {\rm const.\/}
$$
Using these equations and the identity
$$\partial F_k(0,0)=
{\rm Im\/}\big[P_k \overline{\partial Q_k}+\partial P_k \overline{Q_k}\big],$$
we just expand everything out and cancel many terms in pairs.
We find that
$$\partial F_k(0,0)=Xk^3+k^2(Y_0+Y_1+Y_2)+O(k).$$
Here
$$X=\frac{M}{2}{\rm Im\/}\big[-iP^{\#}\overline{Q^{\#}}+iP^{\#}\overline{Q^{\#}}\big]=0;$$
$$Y_0={\rm Im\/}\big[P^{\#} \frac{iM}{2}\overline{Q^{\#}}-
\frac{iM}{2}P^{\#}\overline{Q^{\#}}\big]=0;$$ 

Now we get to the nontrivial quantities.  In our calculations
we use the fact that $P^{\#}$ and $Q^{\#}$ are both real at $X_0$.

$$Y_1={\rm Im\/}\big[P^{\#} \overline{\partial Q^{\#}}+
\partial P^{\#} \overline{Q^{\#})}\big]=
{\rm Im\/}\big[-P^{\#}\partial Q^{\#} + (\partial P^{\#})Q^{\#}]=
-{\rm Im\/}(\delta);$$

$$Y_2={\rm Im\/}\big[-\frac{iM}{2}P_0\overline{Q^{\#}}+\frac{iM}{2}P^{\#}\overline{Q_0}\big]=$$
$${\rm Re\/}\big[-\frac{M}{2}P_0\overline{Q^{\#}}+\frac{M}{2}P^{\#}\overline{Q_0}\big]=$$
$${\rm Re\/}\big[-\frac{M}{2}P_0Q^{\#}+\frac{M}{2}P^{\#}Q_0\big]=
\frac{-M\delta}{2}.$$
This completes the proof.
\endproof

Now we turn to the task of establishing Equation \ref{cg}.

\begin{lemma}
There is some constant $M$ such that
\begin{equation}
 \label{qbound2}
|C_{k,I}| < (Mk)^{|I|+2}.
\end{equation}
\end{lemma}

\startproof
If follows from the linear growth of $\{\widehat P_k\}$ and
$\{\widehat Q_k\}$  and Equation \ref{pairing} that
\begin{equation}
\label{formm1}
F_k(x_1,x_2)=\sum_{j=1}^{N_k} J_{kj} \sin(A_{kj} x_1+B_{kj} x_2),
\end{equation}
Where, for some fixed constant $M$, we have
\begin{equation}
\label{formm2}
N_k<Mk^2; \hskip 15 pt
\max_j(|J_{kj}|)<M \hskip 15 pt
\max_j(|A_{kj}|)<Mk; \hskip 15 pt
\max_j(|B_{kj}|)<Mk.
\end{equation}
Let $D$ be a differential operator of order $\alpha$.  We have
$$DF_k(x_1,x_2)=\sum_{j=1}^{N_k} J_{D,kj} {\rm trig\/}(A_{kj} x_1+B_{kj} x_2),$$
where ${\rm trig\/}$ stands for either the sine or the cosine function,
depending on the parity of $\alpha$.   Equation \ref{formm2} gives us 
$\max_j(|J_{D,kj}|) < M^{\alpha+1} k^{\alpha}.$
Given that the sine and cosine functions lie between $-1$ and $1$, and that
there are at most $Mk^2$ terms in the sum for $DF_k$, we have
\begin{equation}
\sup |DF_k| \leq (Mk)^{\alpha+2}.
\end{equation}

Let $(n,k)$ stand for ''$n$ choose $k$''.
If $I=(i_1,i_2)$ is a multi-index, with $\alpha=|I|$, then we have
\begin{equation}
|C_{k,I}| \leq \frac{(\alpha,i_1) \sup |D_{I}F_k|}{\alpha!} =
\frac{(\alpha,i_1)(Mk)^{\alpha+2}}{\alpha!}.
\end{equation}
Summing over all multi-indices of weight $\alpha$ we get
$$
\sum_{I:\ |I|=\alpha} |C_{k,I}| \leq
\frac{2^{\alpha}(Mk)^{\alpha+2}}{\alpha!} \leq (Mk)^{\alpha+2}.
$$
This is Equation \ref{qbound2}.
\endproof

Here is an improvement on Equation \ref{qbound2}, at least
when $|I|=2$.

\begin{lemma}
There is some constant $M$ such that
\begin{equation}
 \label{qbound}
|I| = 2 \hskip 20 pt \Longrightarrow \hskip 20 pt |C_I| < (Mk)^3
\end{equation}
\end{lemma}

\startproof
To simplify our notation, an expression like $D_IP_k$ shall stand
for $D_IP_k(X_0)$.  The functions $P_k$ and $Q_k$ are always
evaluated at $X_0$.
We will consider $D_{2,0}F_k(0)$, the other
partial derivatives of interest having a similar analysis.
By the chain rule we have
$D_{2,0}F_k(0)=A+B+C,$ where
$$A\!=\!{\rm Im\/}(D_{2,0}P_k\overline{D_{0,0}Q_k}); \hskip 10 pt
  B\!=\!{\rm Im\/}(D_{1,0}P_k\overline{D_{1,0}Q_k}); \hskip 10 pt
  C\!=\!{\rm Im\/}(D_{0,0}P_k\overline{D_{2,0}Q_k}).$$
From Lemma \ref{growth} and our assumptions, there are
constants $a,b \in \R$ such that
$$D_{2,0}P_k = ak^2 D_{0,0}Q_k+O(k^2); \hskip 30 pt
  D_{1,0}P_k=b D_{1,0}Q_k+O(k).$$
Hence
$$A={\rm Im\/}((ak^2D_{0,0}Q_k \times \overline{D_{0,0}Q_k})+O(k^2) \times O(k) =
O(k^3),$$ and
$$B={\rm Im\/}((bD_{1,0}Q_k \times \overline{D_{1,0}Q_k})+O(k) \times O(k) =
O(k^2).$$
The term $C$ has the same treatment as $A$.  All in all,
$|D_{2,0}F|=O(k^3)$.
\endproof

\begin{lemma}
\label{contg}
Equation \ref{cg} holds for $\epsilon=1/4$.
\end{lemma}

\startproof
Without loss of generality, we may take $k>2+M^{100}$.
Given Equations 
\ref{qbound2} and \ref{qbound} we have
$$
\sum_{|I| \geq 2} |C_{k,I}| k^{(-2+\epsilon)|I|}  \leq 
\sum_{|I|=2} |C_{k,I}| k^{-4+1/2} +
\sum_{\alpha=3}^{\infty} \sum_{|I|=\alpha} |C_{I,k}| k^{-7\alpha/4} \leq $$
$$3M k^{-1/2} + \sum_{\alpha=3}^{\infty} M^{\alpha+2} k^{-(3\alpha/4)+2}
\leq
3M k^{-1/2} + \sum_{\alpha=3}^{\infty} k^{-(3\alpha/2)+2+(\alpha/200)+(1/100)} \leq $$
$$3Mk^{-1/2} + \sum_{k=3}^{10} k^{-1/8} + \sum_{\alpha=10}^{\infty} k^{-\alpha} \leq
3Mk^{-1/2}+ 7 k^{-1/8} + 2 k^{-1}.$$
This last expression tends to $0$ as $k \to \infty$.
\endproof

The Quadratic Rescaling Theorem is an immediate consequence of
the Convergence Lemma, Lemma \ref{constantseq}, Lemma \ref{gradseq},
and Lemma \ref{contg}.

\newpage

%% file: 7pivot.tex
\section{Calculating the Pivot Region}
\label{sect:pivot}

\subsection{The Main Result}

We defined the pivot region $R_{nk}$ in \S \ref{s3.3}.
(We will recall the definition in the next section.)  
In this section we compute the asymptotic shape of
this region when $n$ is fixed and $k \to \infty$.
Let $T_{nk}$ denote the dilation which maps
$V_n$ to $(0,0)$ and
dilates by a factor of $k^2$. 
Define
\begin{equation}
C_n=\frac{s}{(2n-2)c}; \hskip 30 pt
c=\cos(\frac{\pi}{2n}) \hskip 20 pt
s=\sin(\frac{\pi}{2n}).
\end{equation}

\begin{lemma}[Pivot]
For any $n$, the set $T_{nk}(R_{nk})$ converges to the
infinite strip $\Sigma_n$ defined by the inequalities
$|x-y| < C_n$.
\end{lemma}

\noindent
{\bf Remarks:\/} 
(i)
When we restrict to any bounded region of the
plane, the convergence we have in mind is the same
discussed in 
Theorem \ref{limit1}. \newline
(ii) 
The Pivot Lemma is sharp.
As we will 
see in the next chapter, the limit
$$\lim_{k \to \infty} T_{nk}(O(W_{nk}))$$
turns out to have vertices on
both components of $\partial \Sigma_n$.,
\newline

The rest of the chapter is devoted to proving the
Pivot Lemma.

\subsection{Reducing to Defining Functions}

Suppose for the moment that $e(t)$ is a continuously varying
segment in $\R^2$ for $t \in [0,1]$.
Suppose that $e(0)$ has negative slope.
Let $f(t)$ denote the height of the left
endpoint of $e(t)$ minus the height of the right endpoint
of $e(t)$.
Note that $f(0)>0$.
We would like to make the following statements, which we
call the {\it slope statements:\/}
\begin{enumerate}
\item
$e(t)$ has negative slope $\forall t \in [0,1]$ if
$f(t)>0$ $\forall t \in [0,1]$.
\item If $f(t)<0$ for some parameter $t$ then $e(t)$ has
positive slope at $t$.
\end{enumerate}
We would like the slope statements to be true because we
would like to define the pivot region in terms of
the defining functions associated to the endpoints
of the QH edges.
Unfortunately, the slope statements are not
necessarily true.  The problem is that $e(t)$
could become vertical at some point.  However,
the slope statements are true if $e(t)$ has
finite slope for all $t \in [0,1]$.
Most of this section is devoted to dealing with
this irritating hitch in the slope statements.
Once we have the kink worked out, we will proceed
to define the pivot region in terms of defining
functions.

Recall that $R_{nk}$ is the path component of
$R'_{nk}$ which contains $V_n$; and
$R'_{nk}$ is the subset of $A_n$ consisting
of points where all the QH edges of $U_{nk}$ have
negative slope.  (We defined $A_n$ in \S \ref{s3.3}; this
set is about to be replaced so we will not bother
to recall the definition here.)
Let $A_{nk} \subset \Delta$ denote the subset consisting
of points which are within $k^{-3/2}$ of $P_n$.
For $k$ sufficiently large, $A_{nk}$ is a subset
of $A_n$, the set defined in \S \ref{s3.3}.
Notice that $T_{nk}(A_{nk})$ is a disk of
radius $k^{1/2}$. Hence
$\lim_{k \to \infty}
T_{nk}(R'_{nk})=\lim_{k \to \infty}
T_{nk}(R'_{nk} \cap A_{nk})$.
For this reason, we will always work within $A_{nk}$
when we analyze $R'_{nk}$ and $R_{nk}$.

For the next several results we choose two points
$X_0,X_1 \in A_{nk}$. We might as well take
$X_0=V_n$.  Let $U_j=U(W_{nk},X_j)$ for $j=0,1$.

\begin{lemma}
Let $e_0$ and $e_1$ be corresponding edges of
$U_0$ and $U_1$, which are edges of the $j$th
triangle  from the left.  If $U_0$ and $U_1$ are
both rotated so that
the leftmost edge is horizontal then
the angle between $e_1$ and $e_2$ is
at most $O(jk^{-3/2})$.
\end{lemma}

\startproof
The point $\widehat e \in \Z^2$ corresponding to
$e_0$ and $e_1$ has norm $O(j)$.
Also $|X_0-X_1|=O(k^{-3/2})$ by hypothesis.
But the angle between our two edges is
$|\widehat e \cdot (X_0-X_1)|$, a quantity which
is $O(jk^{-3/2})$.
\endproof

\begin{corollary}
If $U_0$ and $U_1$ are both rotated so that the
leftmost edges are horizontal, then the angle
between the holonomy of $U_0$ and the holonomy
o $U_1$ is $O(k^{-1/2})$.
\end{corollary}

\startproof
Here we will use the fact that $X_0=V_n$.
Let $L_0$ denote the line which joins up
the endpoints of the $3$-spine $S_0$ of $U_0$.
Given the structure of $U_0$ discussed in
\S \ref{s3}, we see that $L$ has length which is
linear in $k$.  The holonomy of $U_0$ maps
the left endpoint of $L_0$ to the right endpoint
of $L_0$.
Let $S_1$ be the $3$-spine of $U_1$.  If the
left endpoints of corresponding $j$th edges of $L_0$
and $L_1$ are matched up, then the right endpoints
differ by at most $O(jk^{-3/2})$.  Hence, by vector
addition, we see that the right endpoints of
$S_0$ and $S_1$ differ by at most
$$\sum_{j=1}^{C_nk} O(jk^{-3/2})=O(k^{1/2}),$$
assuming that the left endpoints have been
matched up.  Our notation in the last
estimate is a bit informal. The constant $C_n$
is present in the sum to indicate
that there are at most $C_nk$ edges in $U_{nk}$.
  Hence $L_1$ also has length
which is linear in $k$.  It now follows
from basic trigonometry that the angle
between $L_0$ and $L_1$ is $O(k^{-1/2})$.
\endproof

\begin{corollary}
\label{novert}
If $X \in A_{nk}$ and $k$ is sufficiently large then
none of the QH edges in $U(W_{nk},X)$ is vertical.
\end{corollary}

\startproof
If $U_0$ and $U_1$ are both rotated so as to have
horizontal holonomy then the angle between
corresponding edges of $U_0$ and $U_1$
is at most $O(k^{-1/2})$.  This is an immediate
consequence of the previous
two results, and the fact that there are $O(k)$ triangles in
$U_0$ and $U_1$.
This lemma now follows from Lemma \ref{slope1}.  The
idea is that the QH edges are nearly horizontal
for one point in $A_{nk}$, and then that cannot
rotate much as we move around in $A_{nk}$.
\endproof

Now we
can proceed with the analysis of the region $R_{nk}$
by means of defining functions.  For each QH edge
$e$, let $F_{nk,e}$ denote the defining function
which measures the height of the right endpoint of
$e$ minus the height of the left endpoint of $e$.
(We normalize so that $e$ has length $1$, as in
\S \ref{sect:background}.)   In particular, let $e_1,...,e_8$
be the Qh edges corresponding to the extreme
QH points, as discussed in \S \ref{s3.4}.  Compare
Figure 5.10.
Let $\widetilde R_{nk} \subset A_{nk}$ denote those
points $X$ such that
$F_{nk,a}(X)<0$ for $a=1,...,8$.  Here is the
main result of this section.

\begin{lemma}
If $\lim_{k \to \infty} T_{nk}(\widetilde R_{nk})=\Sigma_n$, then
the Pivot Lemma is true.
\end{lemma}

\startproof
As above, we use the convention that our defining functions
measure the height of the left vertex minus the height
of the right vertex.
Let $X \in \widetilde R_{nk}$.  The Convex Hull Lemma says
that $F_{nk,e}(X)>0$ for all QH edges $e$.  Given Corollary
\ref{novert} we now know that all the QH edges have negative
slope.  Hence $\widetilde R_{nk} \subset R'_{nk}$ for
$k$ sufficiently large.  If
$T_{nk}(\widetilde R_{nk})$ converges to $\Sigma_n$ then
the connected component $U_k$ of
$T_{nk}(\widetilde R_{nk})$ containing $(0,0)$ also
converges to $\Sigma_n$.  Since $\widetilde R_{nk} \subset R'_{nk}$,
we see that $U_k$ is a connected subset of $R'_{nk}$ which contains
$V_n$.  Hence $U_k \subset T_{nk}(R_{nk})$.  In summary,
some subset $U_k$ of $T_{nk}(R_{nk})$ converges to
$\Sigma_n$.  

Suppose we could find a sequence $\{X_k\}$ of points
such that $X_k \in R_{nk}$ but $T_{nk}(X_k)$ converges
to some point of $\R^2-\Sigma_n$.  Then some defining
function $F_{nk,a}$ would be negative at $X_k$.  But then,
by the second slope statement and Corollary \ref{novert},
some QH edge would have positive slope at $X_k$. This
is a contradiction.  Hence
$T_{nk}(R_{nk})$ itself converges to
$\Sigma_n$.
\endproof

\subsection{Bilateral Symmetry}
\label{bilateral}

Now we will focus our attention on the $8$ extreme QH points, the
ones listed in \S \ref{s3.4}.
The important
fact for us is that the coordinates of the $4$
{\it northwest\/} extreme points are
independent of $k$ and the $4$ {\it southwest\/}
extreme points are, in a geometric sense,
symmetrically located with respect to the
northwest extreme points.
Each pseudo-parallel family is bounded by a
northwest point and a southeast point.  We call two
such extreme points {\it partners\/}.

\begin{lemma}
\label{bil}
Suppose that $F_1$ and $F_2$ are the defining functions
associated to a pair of partner extreme points.
Then $F_1(x,y)=F_2(y,x)$.
\end{lemma}

\startproof
The hexpath $H_{nk}$ has bilateral
symmetry across a diagonal line.  If we reflect
$H_{nk}$ in its line of symmetry and trace it
backwards we get the same path.   Correspondingly,
if we rotate $U_{nk}$ by $180$ degrees and
trace it backwards we get a cyclic permutation
of $U_{nk}$, except with the edges of
type $1$ and $2$ reversed.  From this symmetry we
see that {\it partner extreme points correspond
to edges of different types\/}.  
We label so that the extreme point corresponding
to $F_j$ has type $j$.

Consider the effect of changing the origin in $\Z^2$.  This amounts to
adding some vector $V_0$ to all the lattice points.  If we
compute our defining functions with the new labeling at $X$ we
simply multiply both $P$ and $Q$ by the same quantity
$E(X \cdot V_0)$.  Hence $F$ is unchanged.
For the duration of the lemma, we change the origin so that
it lies on the line of bilateral symmetry for the
hexpath $H_{nk}$.   In this case we have
$P_1(x,y)=P_2(y,x)$ by symmetry.   Reflection
in the main diagonal interchanges the two sets
$\widehat Q_1$ and $\widehat Q_2$.  Hence
$Q_1(x,y)=Q_2(y,x)$.   From Equation \ref{pairing}
we see that $F_1(x,y)=F_2(y,x)$.
\endproof

\subsection{The Computation}

According to the symmetry above, we just have to analyze the
defining functions associated to the $4$ northwest extreme points.
Let $\{P_k\}$ and $\{Q_k\}$ and $\{F_k\}$ be the functions
associated to one of these points.   We use the convention
that $F_k$ measures the height of the right vertex minus
the height of the left vertex.
Referring to the setup for the Q.R.T., we have the basic
constants
\begin{equation}
X_0=V_n=(\frac{\pi}{2n},\frac{\pi}{2n}) \hskip 30 pt
M_1=2n-2; \hskip 30 pt M_2=-2n+2.
\end{equation}

The fundamental translation $T$ moves points $2n-2$
units south and east.
Inspecting the figures in \S \ref{s3.1}, we see that the hexpath
$H_{nk}$, interpreted as a function from $\Z^2$ into $\Z$
has $T$-linear growth.  The same therefore is true of the
path $\widehat Q$, which is derived from $H_{nk}$ 
as discussed in \S \S \ref{s2.4}.   By Lemma \ref{slope1}
we have  $Q^{\#}(V_n)=\Psi^{\#} \in \R$.  
(We will give a second derivation of this fact
in \S \ref{s6}, based on the combinatorics of the $3$-spine
of the unfolding.)

Note that $\widehat P_k$ is the indicator function for a
single point whose coordinates do not change with $P$.
Hence $P^{\#}$ is the $0$-function.  The value of
$P_k(V_n)$ is independent of both $k$ and the choice of
northwest extreme point.  For the point $(0,0)$ we can
see that $P_k(V_n)= \pm 1$.  Our convention of the
position of the left vertex minus the position of
the right vertex leads to $P_k(V_n)=-1$.
All in all, we have
\begin{equation}
\label{vanish}
P_0(V_n)=-1; \hskip 16 pt
P^{\#}(V_n)=\partial_i P^{\#}(V_n)=0; \hskip 16 pt
\delta=Q^{\#}(V_n); \hskip 16 pt
\delta_j=0.
\end{equation}
From Lemma \ref{slope1}, we have $\delta=\Psi^{\#} \in \R$.
Hence $\{F_k\}$ satisfies the conclusions of the Q.R.T.

Define
\begin{equation}
c'=\cos(\frac{\pi}{n}); \hskip 30 pt
s'=\sin(\frac{\pi}{n}).
\end{equation}
We have $s'=2cs$ and $c'=2c^2-1$.
From Lemma \ref{slope1} we have

\begin{equation}
Q_0(V_n)=\Psi_1+i4\Psi_2=12(1+c')+4is'=24c^2+8ics.
\end{equation}
\begin{equation}
\label{computeQ}
Q^{\#}(V_n)=\Psi^{\#}=8(1+c')=16c^2.
\end{equation}
(We will give another derivation of this in \S \ref{s6}.
See Equation \ref{check}.)
Combining Equation \ref{vanish} with Equation
\ref{computeQ} we see that
$$
F_0(0,0)={\rm Im\/}(P_0(V_0) \overline{Q_0(V_n)})=8cs; \hskip 30 pt
\delta=16c^2;$$
\begin{equation}
-\Delta_2=\Delta_1=-\frac{M_1}{2}\delta+0=-(16n-16)(c^2)=-\frac{8cs}{C_n}.
\end{equation}
According to the Q.R.T, the family of functions
$\{G_k\}$ converges in the smooth topology to the function
$$
G(x_1,x_2)=8cs - \Delta_1 x_1 - \Delta_2 x_2=
$$
$$
8cs + (16n-16)c^2x_1 - (16n-16)x_2=
$$
\begin{equation}
\label{firstfour}
8cs \times (1+\frac{x_1}{C_n}-\frac{x_2}{C_n}).
\end{equation}
This calculation agrees, for small $n$, with the
automatic computations done by McBilliards.
We now see that $G$ is one of the two
linear functions defining $\Sigma_n$.
The other function comes from symmetry,  
as we have discussed above.
This completes the proof of the Pivot Lemma.

\newpage

%% file: 8tile.tex
\section{Rescaling The Orbit Tiles}
\label{s6}

\subsection{Overview}

We continue using the notation from previous chapters.
In particular, $T_{nk}$ is the dilation which maps
$V_n$ to $0$ and expands by $k^2$.   We are interested
in understanding the limits, as $k \to \infty$,
of the sets 
$$T_{nk}(O(W_{nk})).$$
Let $a_1,a_2,a_3,a_4$ be the
top pivot vertices of $U_{nk}$, labeled from left to right.
Likewise let $b_1,b_2,b_3,b_4$ be the bottom
pivot vertices of $U_{nk}$, labeled from left to right.
For any pair $(p,q)$ of vertices amongst these,
let $F_k[p,q]$ be the corresponding defining function.
We are suppressing $n$ from our notation.  If we want
to evaluate our function at a point $X$, we write it
as $F_k[p,q](X)$.
In computing these functions we will use the following sign
conventions
\begin{itemize}
\item $F_k[a_i,b_j]>0$ iff $a_i \uparrow b_j$.
\item $F_k[a_i,a_j]>0$ iff $a_j \uparrow a_i$.
\item $F_k[b_i,b_j]>0$ iff $b_j \uparrow b_i$.
\end{itemize}
Assuming that $\{F_k[p,q]\}$ satisfies the 
$G[p,q]$ denote the rescaled limit of the sequence $\{F_k[p,q]\}$.
Below we will prove

\begin{lemma}
\label{basis}
Each of the $3$ families
$$\{F_k[a_1,b_2]\}; \hskip 20 pt
\{F_k[b_2,b_3]\}; \hskip 20 pt
\{F_k[b_3,a_4]\}$$ 
satisfies the hypotheses of the Q.R.T.
for $j=2,3,4$.
\end{lemma}

We will compute the scaling limits of these functions
explicitly below.   The work in \S \ref{sect:pivot} shows that each of the
function families
$\{F_k[a_i,b_i]\}$ satisfies the hypotheses (and hence the
conclusion) of the Q.R.T.  Any other function
$F_k[p,q]$ can be written as a linear combination of
the functions $F_k[a_i,b_i]$ and the functions
from Lemma \ref{basis}.  (We will find these
linear combinations explicitly below.)
 Hence, all our function families
satisfy the conclusions of the Q.R.T.

Let $\Omega_n \subset \R^2$ denote the convex set on which
all the functions $G[a_i,b_j]$ are positive.  By the Q.R.T.
$\Omega_n$ is a convex polygon$-$possibly empty or infinite$-$with
at most $16$ sides.   By the same symmetry as discussed in \S \ref{sect:pivot},
we know that $\Omega_n$ is symmetric with respect to reflection
in the line $x_1=x_2$.  We also know that $\Omega_n \subset \Sigma_n$,
because the defining functions for $\Sigma_n$, namely
$G[a_i,b_i]$, are by definition positive on $\Omega_n$.

\begin{lemma}
Assume that $\Omega_n$ is bounded. Then
$T_{nk}(O(W_{nk}))$ converges to $\Omega_n$, in the sense
mentioned in the introduction.
\end{lemma}

\startproof
Let $U$ be any open set whose closure is contained in the
interior of $\Omega_n$.   Let
$X$ be any point of
$T_{nk}^{-1}(U)$.  By the Pivot Lemma, 
$X \in R_{nk}$ for $k$ sufficiently large.
Lemma \ref{ruleout} now says that the
lowest top vertex of $U(W_{nk},X)$ is one of the $a_i$ and
the highest bottom vertex is one of the $b_j$.
Once $k$ is sufficiently large, the functions
$G_k[a_i,b_j]$ will all be positive on $U$. Hence
$F_k[a_i,b_j](X)>0$ provided that $k$ is sufficiently
large.  But now we know that all the top vertices
of $U(W_{nk},X)$ lie above all the bottom
vertex.  Hence $X \in O(W)$.  Hence
$U \subset T_{nk}(O(W))$.  This is Property 1
of our definition of convergence.

Suppose that we can find points $X_k \in O(W_{nk})$ such that
$T_{nk}(X_k)$ converges to a point in $\R^2$ which is not
contained in the closure of $\Omega_n$.  Then, for $k$
sufficiently large, at least one of the functions
$F_k[a_i,b_j]$ is negative on $X_k$.  But then
some bottom vertex of $U(W_{nk},X_k)$ lies above some
top vertex.  Hence $X_k \not \in O(W_{nk})$.  This
is a contradiction.  This contradiction, establishes
Property 2 of our convergence.
\endproof

During our proof of Lemma \ref{basis} we will
gather enough information to compute all the
functions defining $\Omega$ exactly.
It is then a simple matter to check that these
functions
cut out precisely the region
advertised in Theorem \ref{limit1}.

\subsection{Asymptotic Limit Calculations}

Here will prove Lemma \ref{basis} and compute
the rescaled limits for the relevant families of
functions.
For each of the $3$ function families of interest to us,
the corresponding vertices can be
connected to each other using part of the $3$-spine.
Figure 8.1 shows by example the general
pattern for the
paths $\widehat P$ and $\widehat Q$.
The example corresponds to $U_{41}$, but every
picture has the same combinatorial structure.
In our figures in this chapter, the grids have
edgelength $1$ rather than $2$, as in \S \ref{s3}.
We make the change because
we want to point out integer coordinates
of various points in the squarepath, and some of
these coordinates might be odd.

\begin{center}
\includegraphics{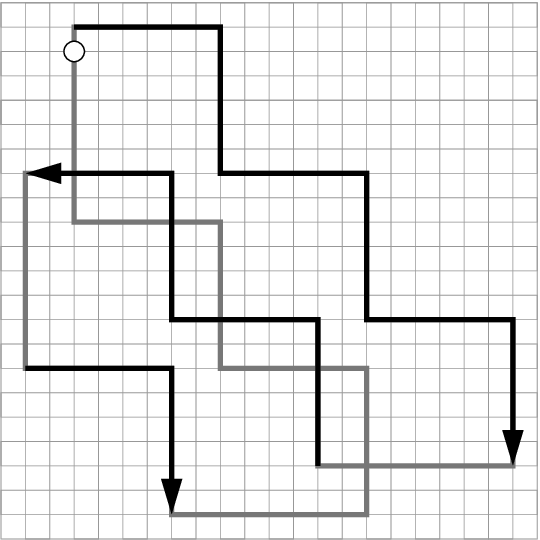}
\newline
Figure 8.1: Paths in $\Z^2$.
\end{center}

The whole path represents
$\widehat Q$.  The white dot denotes the origin.
Each of the three black paths corresponds to
a different one of our function families.
Tracing the path around clockwise, starting
at the origin, we encounter the paths in
the same order that the function families
are listed.
From this picture we see clearly that 
both $\widehat P$ and $\widehat Q$ have $T$ linear growth,
where $T$ is as in \S \ref{sect:pivot}.  
It remains to compute all the quantities relevant to
the Q.R.T., for each of the three families.
We will do this in a step by step fashion.
To keep our pictures concrete, we will draw
the case $n=4$, though the general case is
extremely similar.

\subsubsection{The Holonomy Calculation}

Here we compute the quantities associated to
the family $\{Q_k\}$.  Just by scaling we get
\begin{equation}
\label{check}
Q_0(X_0)=\lambda_n \Psi; \hskip 30 pt
Q^{\#\/}(X_0)=\lambda_n \Psi^{\#}; \hskip 30 pt
\lambda_n=\frac{1}{2c}.
\end{equation}
We explain the constant $\lambda_n$ in
\S \ref{scaled} below.  Unfortunately, 
the geometric method used in the proof
of Lemma \ref{slope1} does not readily
shed light on the derivatives of $Q^{\#}$.
So, here we will use the combinatorial method explained in \S \ref{sect:background}.
At any rate, our calculations here serve as a second proof
of Lemma \ref{slope1}.

\begin{center}
\includegraphics{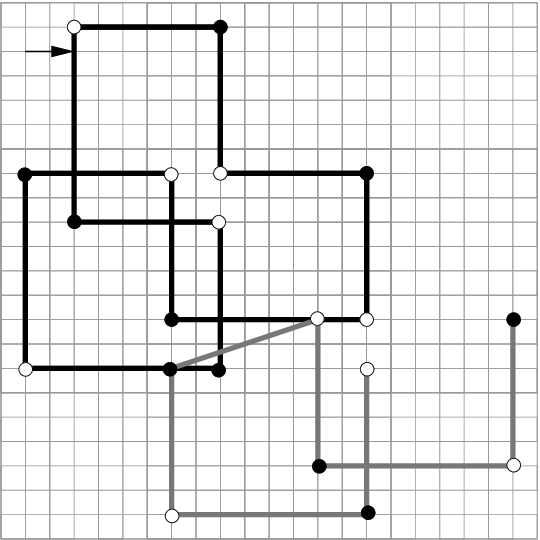}
\newline
Figure 8.2: Fourier Transform of the Holonomy
\end{center}

Figure 8.2 shows $\widehat Q_0$,
using the representation we discussed in \S \ref{sect:background}. The big grey dot
is the origin.  The dots connected by the grey path
are not part of
$\widehat Q_0$.  These dots are the support of
$\widehat Q^{\#}$.  

Let
\begin{equation}
M=2n-2
\end{equation}  
we arrive at the following tableaux for $Q_0$ and $Q^{\#}$.
\begin{equation}
\label{tableau1}
\begin{matrix}
(+)&0&1 \cr
&M&1 \cr
&M&1-M \cr
&2M&1-M \cr
&2M&1-2M \cr
&M-2&1-2M \cr
&M-2&1-M \cr
&-2&1-M \cr
&-2&-1-2M \cr
&M&-1-2M \cr
&M&-1-M \cr
&0&-1-M \end{matrix}
\hskip 50 pt
\begin{matrix}
(+)&& \cr
.&&\cr
.&&\cr
&3M&1-2M \cr
&3M&1-3M \cr
&2M-2&1-3M \cr
&2M-2&1-2M \cr
&-2+M&1-2M \cr
&-2+M&-1-3M \cr
&2M&-1-3M \cr
&2M&-1-2M \cr
.&& \end{matrix}
\end{equation}
We determined the sign for the second tableaux by trial and error:
We expect $Q^{\#}(X_0)$ to be positive real rather than negative
real because our unfoldings grow in the positive direction.
To help show the pattern we have staggered the entries in the
tableau for $\widehat Q^{\#}$ to indicate their corresponding lines
in the tableau for $\widehat Q_0$.

we use the notation
\begin{equation}
c=\cos(\frac{\pi}{2n}); \hskip 30 pt
s=\sin(\frac{\pi}{2n}).
\end{equation}
as in \S \ref{sect:pivot}. 
We could use the Modular Transform Lemma from \S \ref{sect:rescaling} to evaluate
the above functions and their derivatives at $X_0$.  However,
we will do the calculations symbolically in Mathematica.
From Equation \ref{tableau1} we get:
$$
\label{Q0}
Q_0(X_0)=12 c + 4is;
$$
$$
Q^{\#}(X_0)=8c.
$$
$$
\partial_1Q^{\#}(X_0)=8i(4n-5)c;$$
\begin{equation}
\label{Q1}
\partial_2Q^{\#}(X_0)=8(n-1)s -40i(n-1)c.
\end{equation}
In particular, using the double angle formulas we see that
Equation \ref{check} is indeed true.

\subsubsection{The First Function Family}

Here we compute the quantities associated to
$\{P_k\}$, for $P_k[a_1,b_2]$.  Figure 8.3 shows
$\widehat P_0$ and $\widehat P^{\#}$.

The bottom two dots represent $\widehat P^{\#}$.
Setting $M=2n-2$ we read off the following
function tableaux for $P$ and $P^{\#}$.
\begin{equation}
\label{tableau2}
\begin{matrix}
(-)&0&1 \cr
&M&1 \cr
&M&1-M \cr
&2M&1-M \cr
&2M&1-2M \end{matrix}
\hskip 40 pt
\begin{matrix}
(+)&3M&(1-2M) \cr
&3M&(1-3M) \end{matrix}
\end{equation}

\begin{center}
\includegraphics{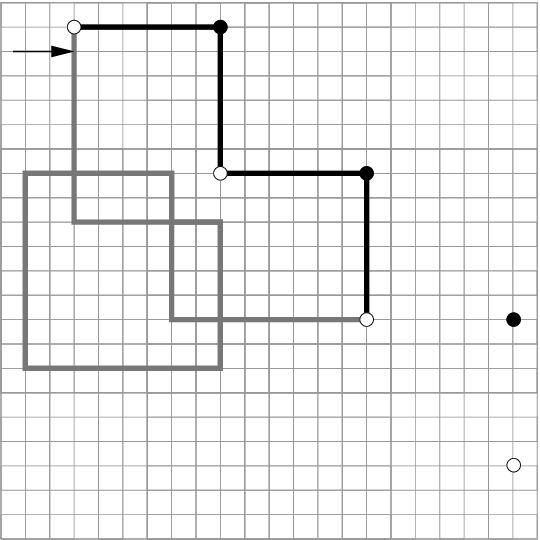}
\newline
Figure 8.3: Fourier Transform for the First Path
\end{center}

\noindent
{\bf Remarks:\/} (i)
The tableau for $\widehat P_0$ consists of the
first $5$ lines of the tableau for $\widehat Q_0$, and
the tableau for $\widehat P^{\#}$ consists of the
first $2$ lines of the tableau for $\widehat Q^{\#}$.
In this sense, the tableaux for $\widehat Q_0$ and
$\widehat Q^{\#}$ form a kind of master list.
\newline
(ii)
The fact that (by inspection) $F_0(0,0)>0$ determines the
global signs for our tableaux.
\newline

We compute symbolically from
Equation \ref{tableau2} that
\begin{equation}
\label{P10}
P_0(X_0)=-5c - i s; \hskip 30 pt
P^{\#}(X_0)=-2 c.
\end{equation}
Next, we compute that:
\begin{equation}
\label{P11}
\partial_1 P^{\#}(X_0)=12 i(n-1) c; \hskip 30 pt
\partial_2 P^{\#}(X_0)=2(n-1)s - 2i(5n-6) c.
\end{equation}
Equations \ref{Q0} and \ref{P10} tells us, in particular,
that $Q^{\#}(X_0)$ and $P^{\#}(X_0)$ are both real.
(We could also deduce this geometrically.
Combining Equations \ref{Q0} and \ref{P10} we compute
that
\begin{equation}
\label{PQ11}
F_0(0,0)={\rm Im\/}(-5c-is)(\overline{12c + 4 i s})=8cs=4 \sin(\frac{\pi}{n}).
\end{equation}
This is the same value we got for the 
family considered in \S \ref{sect:pivot} in connection with the Pivot Lemma.
We don't have an explanation for this agreement.

Next, we compute
\begin{equation}
\delta=\det \left[\begin{matrix}
-5c-is & -2c \cr
12c + 4is & 8c \end{matrix} \right]=-16c^2 \in \R.
\end{equation}
All the relevant quantities are real, and so the associated
family $\{F_k\}$ satisfies the hypotheses of the Quadratic
Rescaling Theorem.  
Using the formulas above we compute
$${\rm Im\/}(\delta_1)={\rm Im\/} \det \left[\begin{matrix}
-2c & -12 i (n-1) c \cr
8c & 8i(4n-5)c \end{matrix} \right] =
16(2n-1)c^2.$$
Similarly
$$\delta_2=16c^2.$$
Finally
\begin{equation}
\Delta_1=\frac{2n-2}{2}\delta+\delta_1=16nc^2; \hskip 30 pt
\Delta_2=\frac{2-2n}{2}\delta+\delta_2=16nc^2.
\end{equation}
By the Q.R.T., the rescaled limit of the sequence $\{F_k\}$ is
\begin{equation}
G(x_1,x_2)=8cs-16nc^2x_1-16nc^2x_2.
\end{equation}
The result agrees with the automatic computations done by
McBilliards for small values of $n$.  Note the similarity
to the middle line of Equation \ref{firstfour}.

\subsubsection{The Second Function Family}

Here we compute the quantities associated to
$\{P_k\}$, for $P_k[b_2,b_3]$.   Comparing
Figures 8.1 and 8.2 we find that the
tableau for $\widehat P_0$ consists of lines
$6,7,8$ of the tableau for $\widehat Q_0$ and
the tableau for $\widehat P^{\#}$ consists of
lines $3,4$ of the tableau for $\widehat Q^{\#}$. 
Here are these tableaux:

\begin{equation}
\label{tableau3}
\begin{matrix}
(-) &M-2&1-2M \cr
&M-2&1-M \cr
&-2&1-M \end{matrix}
\hskip 40 pt
\begin{matrix}
(+)&2M-2&1-3M \cr
&2M-2&1-2M 
\end{matrix}
\end{equation}

From these tableaux we compute that
$$P_0(X_0)=3cis; \hskip 30 pt
P^{\#}(X_0)=2c;$$
\begin{equation}
\partial_1 P^{\#}(X_0)=4 i(2n-3)c; \hskip 30 pt
\partial_2 P^{\#}(X_0)=2(n-1)s- 2i(5n-6)c.
\end{equation}
Using these equations, and the ones above for $Q$, we compute that
\begin{equation}
F(0,0)=0; \hskip 15 pt
\delta=0; \hskip 15 pt
\Delta_1={\rm Im\/}(\delta_1)=16c^2; \hskip 15pt
\delta_2={\rm Im\/}(\delta_2)=-16c^2.
\end{equation}
There $\{F_k\}$ satisfies the hypotheses of the Q.R.T. and
the rescaled limit is
\begin{equation}
G(x_1,x_2)=-16c^2x_1+16c^2x_2.
\end{equation}
\newline
\newline
{\bf Remark:\/}
For small values of $n$ this agrees with the
calculations made by McBilliards.  For this
example, it took many tries before we got the
sign right.  The problem is that the constant
term vanishes, making it much trickier to
deduce the correct sign without making an
error.

\subsubsection{The Third Function Family}

Here we compute the quantities associated to
$\{P_k\}$, for $P_k[b_3,a_4]$.   Comparing
Figures 8.1 and 8.2 we find that the
tableau for $\widehat P_0$ consists of line $9$
of the tableau for $\widehat Q_0$ and
the tableau for $\widehat P^{\#}$ consists of
lines $5,6$ of the tableau for $\widehat Q^{\#}$. 
Here are these tableaux:

\begin{equation}
\label{tableau4}
\begin{matrix}
&-2&-1-2M \end{matrix}
\hskip 50 pt
\begin{matrix}
(+)&-2+M&1-2M \cr
&-2+M&-1-3M \end{matrix}
\end{equation}
From these tableaux we compute that
$$P_0(X_0)=c+is; \hskip 30 pt
P^{\#}(X_0)=2c;$$
\begin{equation}
\partial_1 P^{\#}(X_0)=4 i(n-2)c; \hskip 30 pt
\partial_2 P^{\#}(X_0)=2(n-1)s-2i(5n-4)c.
\end{equation}
Using these equations, and the ones above for $Q$, we compute that
$$
F(0,0)=8cs; \hskip 15 pt
\delta=-16c^2; \hskip 15 pt
\delta_1=(16+32n)c^2; \hskip 15 pt
\delta_2=16c^2; 
$$
\begin{equation}
\Delta_1=16nc^2; \hskip 15 pt
\Delta_2=16nc^2.
\end{equation}
Therefore $\{F_k\}$ satisfies the hypotheses of the Q.R.T, and
the rescaled limit is
\begin{equation}
G(x_1,x_2)=8cs-16nc^2x_1-16nc^2x_2.
\end{equation}
For small values of $n$ this agrees with the calculations
made by McBilliards.

\subsection{The Fudge Factor}
\label{scaled}

Suppose that $p_1,p_2,p_3$ are three
vertices on our unfolding.  Let $F_{ij}$
be the defining function which computes
(up to scale) the height difference
between $p_i$ and $p_j$.
Here we refer to the function defined
in \S \ref{sect:background}.  We would like to say that
$F_{13}=F_{12}+F_{23}$
but there is a catch. The
functions might not all be computed with respect
to the same spine of the unfolding, and thus
the function values might represent differences
in heights measured at different scales.  This explains
the fudge factor $\lambda_n$ in Equation \ref{check}.
In this section we will address this issue systematically.

Let $\theta_d(X)$ denote the sine of the angle
of the triangle $T_X$ which is opposite the $d$th edge.
Supposing that $F$ has been computed in terms of
the $d$-spine, we define
\begin{equation}
\widetilde F(X)=\sin^2(\theta_d)F.
\end{equation}
Then $\widetilde F$ measures the height difference
between the relevant points when the edge of type
$d$ is scaled to have length $\sin(\theta_d)$.
the exponent $2$ appears in the definition
because the functions $P$ and
$Q$ both scale linearly with the edge length, and
$F$ is the imaginary part of their product.

It follows from the Law of Sines that a triangle
may be scaled, with a single scale, so that its
type $d$ edge has length $\sin(\theta_d)$.
Therefore, the functions $\widetilde F$ is
computed at the same scale, independent of which
spine it uses.  Thus we really do have
\begin{equation}
\widetilde F_{13}=\widetilde F_{12}+\widetilde F_{23}.
\end{equation}

Our modification works well with the Q.R.T.  The
set of functions that satisfy the conclusion of
the Q.R.T. with respect to the same
point $X_0$ forms an
algebra.  That is, they can be added, scaled, and
multiplied together.  Letting
$A_k(X)=\sin(\theta_d(X))$, a function which
is actually independent of the parameter $k$.
We see easily that the sequence
$\{A_k\}$ satisfies the conclusions of the Q.R.T.
and has rescaled limit function $\sin_d(X_0)$.
Therefore, if $\{F_k\}$ satisfies
the conclusions of the  Q.R.T. and
has rescaled limit $G$, then
$\{\widetilde F_k\}$ also satisfies the
conclusions of the Q.R.T. and has rescaled
limit
\begin{equation}
\widetilde G=\sin(\theta_d(X_0) G.
\end{equation}
In the examples of interest to us, we have the
following conversions: If $G$ is based on the $3$-spine then
\begin{equation}
\widetilde G=(s')^2G=4c^2s^2G.
\end{equation}
If $G$ is based on either the $1$-spine or the $2$-spine then
\begin{equation}
\widetilde G=s^2G.
\end{equation}
Working with the $\widetilde G$ limits instead of the $G$ limits,
we can add and subtract with impunity.

\subsection{The Shapes of the Tiles}

In this section we compute the region $\Omega_n$, using the
explicit formulas for the functions defining $\Omega_n$.
Our first task is to write down all these functions.
To make our notation simpler, we will understand that
our functions are always evaluated at the point
$(x_1,x_2)$.  We will also use the notation
\begin{equation}
\left[\begin{matrix}A\cr B\cr C \end{matrix} \right]=A+Bx_1+Cx_2.
\end{equation}
In the language of \S \ref{bilateral}, the first and
third hinges of the unfolding $U_{nk}$ correspond to
northwest extreme points, and the second and fourth
hinges correspond to southeast extreme points.
Combining Lemma \ref{bil} and Equation \ref{firstfour} we
see that

\begin{equation}
\widetilde G[a_1,b_1]=\widetilde G[a_3,b_3]=
\left[\begin{matrix}8cs^3 \cr 16(n-1)c^2s^2 \cr
-16(n-1)c^2s^2 \end{matrix} \right]
\end{equation}

\begin{equation}
\widetilde G[a_2,b_2]=\widetilde G[a_4,b_4]=
\left[\begin{matrix}8cs^3 \cr -16(n-1)c^2s^2 \cr
16(n-1)c^2s^2 \end{matrix} \right]
\end{equation}

Again, these are the functions which define the strip
$\Sigma_n$ from the Pivot Lemma.  Summarizing the
calculations we made in this chapter, we have

\begin{equation}
\label{a1b2}
\widetilde G[a_1,b_2]=\widetilde G[a_4,b_3]=
\left[\begin{matrix}32c^3s^3 \cr - 64nc^4s^2 \cr -64nc^4s^2 \end{matrix} \right];
\hskip 30 pt
\widetilde G[b_2,b_3]=\left[\begin{matrix}0 \cr -64c^4s^2 \cr +64c^4s^2 \end{matrix} \right].
\end{equation}

To explain the rules we use to compute the rest of the defining functions,
we simplify our notation.  We let $[a_ib_j]$ stand for $\widetilde G[a_i,b_j]$.
Using our sign conventions above (and checking the signs against the
output from McBilliards) we find that
\begin{itemize}
\item $[a_1b_3]=[a_1b_2]-[b_2b_3]$
\item $[a_4b_2]=[a_4b_3]+[b_2b_3]$
\item $[a_1b_4]=[a_1b_3]-[a_4b_3 ]+[ a_4b_4]$
\item $[a_2b_1]=[a_2b_2]-[a_1b_2]+[a_1b_1]$
\item $[a_2b_3]=[a_2b_2]-[b_2b_3]$
\item $[a_2b_4]=[a_2b_3]-[a_4b_3]+[a_4b_4]$
\item $[a_3b_1]=-[a_1b_3]+[a_3b_3]+[a_1b_1]$
\item $[a_3b_2]=[a_3b_3]+[b_2b_3]$
\item $[a_3b_4]=-[a_4b_3]+[a_3b_3]+[a_4b_4]$
\item $[a_4b_1]=[a_4b_4]+[a_1b_1]-[a_1b_4]$
\end{itemize}
Since we are only interested in pairs of the
form $[a_ib_j]$ we further compress our notation
and write $[ij]=[a_ib_j]=\widetilde G[a_i,b_j]$.
We computed all the above quantities
symbolically and noticed a lot of
symmetry.  To help express this symmetry we
write $[i_1j_1] \sim [i_2j_2]$ if the map
$(x_1,x_2) \to (x_2,x_1)$ 
conjugates the one function to the other.
We compute the following relations:
\begin{enumerate}
\item $[13] \sim [42]$.
\item $[24]\sim [31]$.
\item $\frac{1}{2} \times [13]+\frac{1}{2} \times [42]=[12]=[43]$.
\item $\frac{1}{2} \times [24]+\frac{1}{2} \times [31]=[21]=[34]$.
\item $\frac{1}{2} \times [13]+\frac{1}{2} \times [31]=[11]=[33] \sim
[22] = [44]$.
\item $t \times [11]+(1-t)\times [44]=[14]=[32] \sim [23]=[41]$, where
$t=2c^2/(n-1)$.
\end{enumerate}
(We recall that $c=\cos(\pi/2n)$.)

The above relations easily imply that all the defining
functions are convex combinations of
$$[13] \sim [42]; \hskip 15 pt
[24] \sim [31]$$
Hence $\Omega_n$ is defined by these $4$ alone.
Now that we are done adding these functions together, we
can replace them by positive multiples and still define
the same region in the plane.  Setting
\begin{equation}
\sigma=\sec(\pi/n)=\frac{1}{2c^2-1}
\end{equation}
we compute

\begin{equation}
[13] \propto \left[\begin{matrix}1 \cr
-2(n-1) (c/s) \cr
-2(n+1) (c/s)\end{matrix} \right]; 
\hskip 20 pt
[24] \propto \left[\begin{matrix}-1 \cr
2(c/s)(n+1+\sigma) \cr
2(n-1)(c/s)(1+2 \sigma) \end{matrix} \right].
\end{equation}

Now we dilate the plane by 
$$\zeta_n:=2(n-1)(c/s).$$
The region $\zeta_n\Omega_n$
is the region cut out by the defining functions obtained from
the ones above (and their symmetric conjugates)
by dividing all the second and third coordinate entries by $\zeta_n$.
That is, $\zeta_n \Omega_n$ is defined by the functions:
\begin{equation}
\label{uuu}
\left[\begin{matrix}1\cr -\frac{n+1}{n-1} \cr -1 \end{matrix} \right]; \hskip 20 pt
\left[\begin{matrix}-1\cr \frac{1+n+2 \sigma}{n-1} \cr 1+2 \sigma \end{matrix} \right]; \hskip 20 pt
\left[\begin{matrix}-1\cr 1+2 \sigma \cr \frac{1+n+2 \sigma}{n-1} \end{matrix} \right]; \hskip 20 pt
\left[\begin{matrix}1\cr -1 \cr -\frac{n+1}{n-1} \end{matrix} \right].
\end{equation}

Setting
\begin{equation}
\mu_n=\frac{1}{2}-\frac{\tan^2(\pi/2n)}{2}
\end{equation}
we now list the vertices from Theorem \ref{limit1}, modified so that their
first coordinate is padded with a $1$.
\begin{equation}
\label{vvv}
\left[\begin{matrix}1 \cr -\frac{1}{n} \cr 1-\frac{1}{n} \end{matrix} \right]; \hskip 20 pt
\left[\begin{matrix}1\cr\frac{1}{2}-\frac{1}{2n} \cr\frac{1}{2}-\frac{1}{2n} \end{matrix} \right]; \hskip 20 pt
\left[\begin{matrix}1 \cr1 -\frac{1}{n} \cr -\frac{1}{n} \end{matrix} \right]; \hskip 20 pt
\hskip 20 pt
\left[\begin{matrix}1\cr\mu_n(\frac{1}{2}-\frac{1}{2n}) \cr\mu_n(\frac{1}{2}-\frac{1}{2n}) \end{matrix} \right];
\end{equation}

We claim that $\zeta_n\Omega$ is the convex hull of the vertices from
Theorem \ref{limit1}.  To see this, it suffices to show that
the matrix of dot products between the vectors
in Equation \ref{uuu} and the vectors in Equation \ref{vvv}
is non-negative, and has two zeros in each row and column.
Here is the matrix

\begin{equation}
\left[\begin{matrix}
\frac{2}{n-1}&\beta&0&0 \cr
0&\sigma &\sigma & 0 \cr
0 & 0 &\beta &\frac{2}{n-1} \cr
\frac{1}{2c^2}&0&0&\frac{1}{2c^2} \end{matrix}\right]; 
\hskip 30 pt
\beta=\frac{2 \sigma(n-2-\cos(\pi/n))}{n-1}
\end{equation}
This completes our verification that $\zeta_n \Omega_n$ is the
convex hull of the above mentioned vertices. 

We got 
$\Omega_n$ as the limit of rescaling by a quadratic family
$\{T_{nk}\}$ of dilations.  If we had used the
quadratic family $\{\zeta_n T_{nk}\}$ instead, we
would get $\zeta_n\Omega$ right on the nose.  This
is the family $\{S_{nk}\}$ we use for the proof of Theorem
\ref{limit1}.  This completes our proof of
Theorem \ref{limit1}. 

As a final remark, we calculate that both
functions $[13]$ and $[31]$ vanish at a common
vertex of $\Omega_n$.  Hence the function $[11]$,
one of the defining functions for $\Sigma_n$, also vanishes
at this vertex. Hence $\Omega_n$ has a vertex on
$\partial \Sigma_n$.  By symmetry, $\Omega_n$ intersects
both components of $\partial \Sigma_n$ in a vertex.
This justifies our comment, in \S \ref{sect:pivot}, that the
Pivot Lemma is sharp.

\newpage

%% file: 9stability.tex
\section{Stability questions}
\label{sect:stability}
\def\D{{\mathcal D}}
\def\SL{{\textit{SL}}}

In this section, we prove Theorem \ref{thm:stability}. Recall that $V_n$ for $n \geq 3$ is the obtuse isosceles triangle with two angles of measure $\pi/(2n)$. We will prove that
if $n$ is a power of two, then $V_n$ has no stable periodic billiard paths. For $n$ not a power of two, we will construct a stable periodic billiard path.

\begin{remarks}[Errors Corrected]
\label{rem:corrections}
This version of the paper corrects errors in the published version of the document, \cite{HS}. These errors occurred in sections 9.3 and 9.5.
Namely, Lemma \ref{lem:stability}
is corrected and clarified from the corresponding statement in the published document. More fundamentally, the homology classes
corresponding to stable periodic billiard paths in the published version of section 9.5 were incorrect. This version contains corrected
formulas for these classes, and more detailed proofs. The authors would like to thank Yilong Yang for pointing out the errors in \cite{HS}.
\end{remarks}

\subsection{A Homological Condition for Stability}

Given a triangle $T$, let ${\mathcal D}T$ denote its double with the vertices removed. (The double
of a polygon can be thought of a pillowcase for a pillow made in the shape of the polygon.)
See Figure 9.1. There is a natural folding map $f:{\mathcal D}T \to T$ which sends each of the two triangles making up ${\mathcal D}T$
isometrically to $T$. (This map folds in the sense that is $2$-$1$ except on the edges of the triangle, where it is $1$-$1$.) If $\tilde p$ is a closed geodesic on $\D T$, then $f(\tilde p)$ is a periodic billiard path. Conversely,
if $p$ is a periodic billiard path in $T$ (which hits an even number of sides in a period), then there is a closed geodesic lift $\tilde p$ 
with $f(\tilde p)=p$. The lift $\tilde p$ is unique up to the single non-trivial automorphism of the folding map $f$, which preserves the labeling of edges and swaps the two triangles.

\begin{lemma}
\label{lem:homological_criterion}
A periodic billiard path $p$ in a triangle $T$ is stable if and only if its lift $\tilde p$ to the double $\D T$ is null homologous (equivalent to zero in $H_1(\D T, \Z)$).
\end{lemma}

\noindent {\bf Remark:} Because we removed the vertices of the triangles from $\D T$, $\D T$ is topologically a $3$-punctured sphere. Thus, $H_1(\D T, \Z) \cong \Z^2$

\startproof
This is formally a restatement of Lemma \ref{lem:stability}. The total sign counted for the edge labeled $1$ in Lemma \ref{lem:stability} 
is equivalent to computing the algebraic sign of the intersection of $\tilde p$ with the lift of the edge labeled $1$ to $\D T$. Having zero algebraic intersection
number with each of the lifts of an edge to $\D T$ is equivalent to being null homologous.
\endproof


\begin{center}
\includegraphics[width=3.5in]{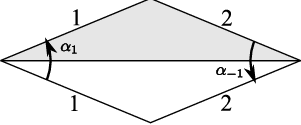}
\newline
Figure 9.1 The double $\D T$ of the triangle $T=V_4$. The numerals indicate pairs of edges which are glued together by an orientation preserving isometry of the plane.
\end{center}

\subsection{Translation surfaces and Veech's lattice property}

We will need to understand some of the implications of work of Veech \cite{V}. For this, we introduce some of the ideas appearing in the study of translation surfaces.
See \cite{MT} for a more thorough introduction.

A {\em translation surface} is a union of polygonal subsets of the plane with edges glued together pairwise by translation. There is a natural translation surface associated to every triangle $T$. 
Let $G=\langle r_1, r_2, r_3 \rangle$ denote the subgroup of $\textit{Isom}(\R^2)$ generated by the reflections in the sides of the triangle $T$. 
The translation surface $S(T)$ is the disjoint union of the triangles $g(T)$ with $g \in G$ with some identifications. Two triangles $g_1(T)$ and $g_2(T)$ are identified by translation if
$g_1 \circ g_2^{-1}$ is a translation. Also, we identify two triangles $g_1(T)$ and $g_2(T)$ along the $i$-th edge by translation if $g_1 \circ g_2^{-1}$ 
can be written as a composition of $r_i$ and a translation.

The resulting surface $S(T)$ is the smallest translation surface cover of $\D T$. The covering map is the map which sends each triangle in
$S(T)$ isometrically to the triangle in $\D T$ with the same orientation. In this section, we follow the convention that the vertices of the triangles making up 
$S(T)$ are removed. (These points are really cone points, we only remove them to make our topological notation simpler.)

Since a translation surface is built from polygonal subsets of the plane glued together by translations, the surface inherits a notion of the {\em direction} of a unit tangent vector. This notion of direction is just the measure of angle compared to a horizontal vector. This notion of direction is a fibration from the unit tangent bundle $T_1 S(T)$ to the circle $\R/2 \pi \Z$.

There is a natural action of the affine group $\SL(2,\R)$ on the space of translation surfaces. Suppose $A \in \SL(2,\R)$ and $S$ is a translation surface, then we will define
$A(S)$. Suppose $S$ is the disjoint union of the polygonal subsets of the plane $P_i$ with $i \in \Lambda$ with edges identified pairwise by translation. Let $A(S)$ be the 
disjoint union of the polygons $A(P_i)$ with $A$ acting linearly on the plane with the same edge identifications. The new edge identifications are also by translation. This
is possible because $A$ sends parallel lines to parallel lines and preserves the ratio's of lengths of pairs of parallel segments. 

The {\em Veech group} $\Gamma(S)$ is the subgroup of elements $A \in \SL(2,\R)$ such that there is a direction preserving isometry $\varphi_A: A(S) \to S$. We abuse notation by
using $A$ to denote the natural map from $S \to A(S)$ given by the restriction of the action of $A$ on the plane to the polygonal subsets of the plane making up $S$. 
The map $\varphi_A \circ A: S \to S$ is known as an affine automorphism of $S$. The set of such maps forms a group known as the {\em affine automorphism group} of $S$.

The Veech group $\Gamma(S)$ is always discrete. We say $S$ has the {\em lattice property} if $\Gamma(S)$ is a lattice in $\SL(2,\R)$. We will utilize the following consequence of Veech's work. If a translation surface has the lattice property, then the collection of closed geodesics on $S$ can be well understood.

A direction $\theta \in \R/2 \pi \Z$ is called a {\em completely periodic direction} for a translation surface $S$ if every bi-infinite geodesic in this direction is closed.

\begin{theorem}[Veech Dichotomy]
Suppose $S$ has the lattice property.
Let $\theta \in \R/2 \pi \Z$ be a direction. Then either $\theta$ is a completely periodic direction for $S$ or the geodesic flow in the direction $\theta$ is uniquely ergodic. Moreover, $\theta$ is completely periodic if and only if there is a parabolic $A \in \Gamma(S)$ for which $\theta$ is an eigendirection.
\end{theorem}


\begin{center}
\resizebox{!}{3.5in}{\includegraphics{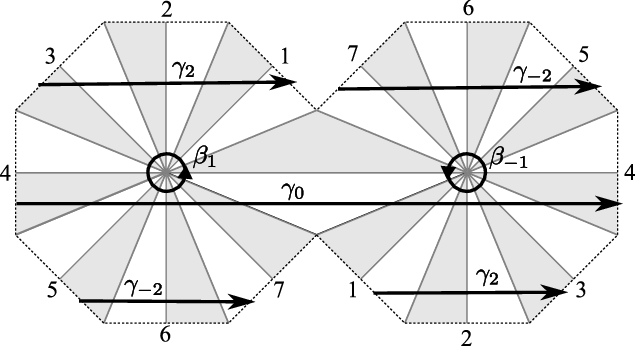}}
\end{center}
Figure 9.2: The translation surface $S(V_4)$. All but two of the obtuse isosceles triangles have been split along their axes of symmetries.
Numbers indicate edge identifications by translations. Curves in the homology classes $\beta_{-1}$, $\beta_1$, $\gamma_{-2}$, $\gamma_0$ and $\gamma_2$ are shown.
\newline

Veech showed that $S(V_n)$ has the lattice property.

We will now carefully describe $S(V_n)$, so that we can explicitly use Veech's property. See Figure 9.2 for visual guidance. 
To build $S(V_n)$, start with a copy of $V_n$ oriented in the plane so that longest side lies on the $x$-axis. Cut this triangle in two along the axis of symmetry. Reflecting each half repeated along the edges it shares with $V_n$  yields 
two regular $2n$-gons. The halves can then be reassembled by gluing according to appropriate translations. This amounts to gluing each edge of the left regular $2n$-gon to the opposite side of the right 
$2n$-gon by translation. We remove the center point of each $2n$-gon and also the vertices of the $2n$-gon, since these points correspond to vertices of lifts of our triangle $V_n$. (This removal of vertices will make our homological notation simpler.)

\begin{theorem}[Veech]
Let $n \geq 3$.
The Veech group $\Gamma\big(S(V_n)\big)$ is a hyperbolic $(n, \infty, \infty)$ triangle group. The group 
$\Gamma\big(S(V_n)\big)= \langle A, B~:~(A \circ B^{-1})^n=e\rangle$ is generated by parabolics which fix the directions of angle $0$ or $\frac{\pi}{2n}$. The corresponding affine automorphisms act as single right Dehn twists in each maximal cylinder in the respective eigendirection.
\end{theorem}

We enjoy the following consequence.

\begin{lemma}[Enumeration]
If $p_1$ is a closed geodesic on $S(V_n)$ then there is a closed geodesic $p_0$ in the directions $0$ or $\frac{\pi}{2n}$ and an affine automorphism
$\varphi_D \circ D$ which maps $p_1$ onto $p_0$.
\end{lemma}
\startproof
The direction of $p_1$ is not uniquely ergodic. Thus by Veech dichotomy, this direction must be an eigendirection for a parabolic $C \in \Gamma\big(S(V_n)\big)$.
Because of the structure of the group, there must be an $D \in \Gamma\big(S(V_n)\big)$ and an integer $k \neq 0$ such that either
$D\circ C \circ D^{-1}=\pm A^k$ or $D\circ C \circ D^{-1}=\pm B^k$. Then $D$ maps the direction of $p_1$ onto either the direction $0$ or $\frac{\pi}{2n}$. It follows
that the affine automorphism corresponding to $D$ sends the geodesic $p_1$ to a geodesic which travels in either the direction $0$ or $\frac{\pi}{2n}$. This image is our
$p_0$.
\endproof

\subsection{Generators for the Affine Automorphism Group}

In this subsection, we combine the idea of Lemma \ref{lem:homological_criterion} with Veech's lattice property to yield necessary and sufficient conditions for stability of a periodic billiard path
in $V_n$ constructed via the affine automorphism group of $S(V_n)$.

Topologically, ${\mathcal D}V_n$ is a 3-punctured sphere. We choose generators $\alpha_{1}$ and $\alpha_{-1}$ for $H_1({\mathcal D}V_n, \Z) \cong \Z^2$. See Figure 9.1.


\begin{center}
\resizebox{!}{3.5in}{\includegraphics{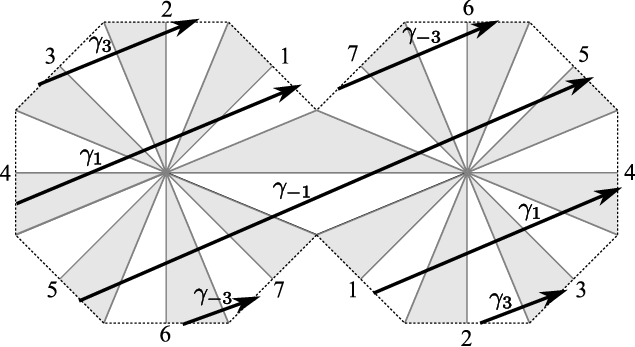}}
\end{center}
Figure 9.3: The translation surface $S(V_4)$. Curves in the homology classes $\gamma_{-3}$, $\gamma_{-1}$, $\gamma_1$ and $\gamma_3$ are shown
\newline

We choose a basis for the homology group $H_1\big(S(V_n), \Z\big) \cong \Z^{2n+1}$. 
Our basis is the homology classes of the collection of curves $${\mathcal B}=\{ \beta_1, \beta_{-1}, \gamma_{1-n}, \ldots, \gamma_{n-1}\}$$ depicted in Figures 9.2 and 9.3.
 The homology classes $\gamma_k$ with $k$ even all contain horizontal geodesics. We order them so that 
the portions of these geodesics in the left polygon increase in $y$ coordinate as the index $k$ increases. The homology class $\gamma_0$ is chosen so that a geodesic in this class travels below the centers of the two $2n$-gons.
The homology classes $\gamma_k$ with $k$ odd all contain geodesics which travel in the direction of angle $\frac{\pi}{2n}$. With the correct choice of indices, we have algebraic intersection numbers
$\gamma_{2i} \cap \gamma_{2i+1}=1$ and $\gamma_{2i} \cap \gamma_{2i-1}=1$, but no other pairs of curves intersect.

We think of first cohomology, $H^1\big(S(V_n), \Z\big)$, as the dual space to homology. 
Let $${\mathcal B}^\ast=\{\beta_1^\ast, \beta_{-1}^\ast, \gamma_{1-n}^\ast, \ldots, \gamma_{n-1}^\ast \}$$
denote the dual basis for $H^1\big(S(V_n), \Z\big)$.

Consider the natural covering map $\phi:S(V_n) \to {\mathcal D}T$. Let $\phi: H_1\big(S(V_n), \Z\big) \to H_1({\mathcal D}T, \Z)$ be the induced map on homology. This satisfies the following properties:
\begin{itemize}
\item $\phi(\beta_i)=2n \alpha_i$ for $i \in \{1, -1\}$. 
\item $\displaystyle \phi(\gamma_k)=\begin{cases}
(n+k) \alpha_1-(n+k) \alpha_{-1} & \textrm{if $k<0$} \\
n \alpha_1+n \alpha_{-1} & \textrm{if $k=0$} \\
(k-n) \alpha_1+(n-k) \alpha_{-1} & \textrm{if $k>0$}
\end{cases}$
\end{itemize}
Using $\phi$, we get elements of the cohomology group $H^1\big(S(V_n), \Z\big)$, given by the coefficients of $\alpha_1$ and $\alpha_{-1}$.
Denote these elements by $\phi_1^\ast$ and $\phi_{-1}^\ast$. We think of these as maps of the form $H_1\big(S(V_n), \Z\big) \to \Z$.
They satisfy $\phi(x)=\phi_1^\ast(x) \alpha_1+\phi_{-1}^\ast(x) \alpha_{-1}$. We have
\begin{equation}
\label{eq:coverdual}
\begin{array}{c}
\phi_1^\ast=2n \beta_1^\ast+\gamma_{1-n}^\ast+2 \gamma_{2-n}^\ast+ \ldots + n \gamma_0^\ast
-(n-1) \gamma_{1}^\ast-(n-2) \gamma_{2}^\ast- \ldots - \gamma_{n-1}^\ast \\
\phi_{-1}^\ast=2n \beta_{-1}^\ast-\gamma_{1-n}^\ast-2 \gamma_{2-n}^\ast- \ldots - (n-1) \gamma_{-1}^\ast + n \gamma_0^\ast
+(n-1) \gamma_{1}^\ast+ \ldots + \gamma_{n-1}^\ast.
\end{array}
\end{equation}

We now describe the action of the affine automorphism group on $S(V_n)$. 
If $g$ is an element of this affine automorphism group, we also use $g$ to denote
the induced action on homology, and use $g^\ast$ to denote the pullback action. That is,
if $\eta^\ast \in H^1\big(S(V_n), \Z\big)$ and $g$ is an affine automorphism of $S(V_n)$ acting on homology, 
then
$$\big(g^\ast(\eta^\ast)\big)(x)=\eta^\ast\big(g^{-1}(x)\big) \quad \text{for all $x \in H_1\big(S(V_n), \Z\big)$}.$$
(So, viewed as matrices acting on our bases, $g^\ast$ is the inverse transpose of $g$.)
For the formulas that follow,
we follow the convention $\gamma_{-n}=\gamma_{n}=0$ and $\gamma^\ast_{-n}=\gamma^\ast_{n}=0$.
The affine automorphism group of $S(V_n)$ is generated by the following elements:
\begin{itemize}
\item The involution which swaps the two regular $2n$-gons. The actions on homology and cohomology are given by 
\begin{eqnarray}
\sigma:~  \beta_i \mapsto \beta_{-i} ~:~
\gamma_k \mapsto \gamma_{-k} 
\\
\label{eq:sigmaast}
\sigma^\ast:~ \beta_i^\ast \mapsto \beta_{-i}^\ast~:~
\gamma_k^\ast \mapsto \gamma_{-k}^\ast.
\end{eqnarray}

\item The right Dehn twists in the odd cylinders. The actions on homology and cohomology are given by 
\begin{eqnarray}
\tau_o:~ \beta_i \mapsto \beta_{i} ~:~
\gamma_k \mapsto \begin{cases}
\gamma_k & \textrm{if $k$ is odd} \\
\gamma_{k-1}+\gamma_k+\gamma_{k+1} & \textrm{if $k$ is even} 
\end{cases}
\\
\label{eq:tauoast}
\tau_o^\ast:~ \beta_i^\ast \mapsto \beta_{i}^\ast ~:~
\gamma_k^\ast \mapsto \begin{cases}
-\gamma_{k-1}^\ast+\gamma_k^\ast-\gamma_{k+1}^\ast & \textrm{if $k$ is odd} \\
\gamma_k^\ast & \textrm{if $k$ is even} 
\end{cases}
\end{eqnarray}

\item The right Dehn twists in the even cylinders. The actions on homology and cohomology are given by 
\begin{eqnarray}
\tau_e:~ \beta_i \mapsto \beta_{-i} ~:~
\gamma_k \mapsto \begin{cases}
\gamma_k & \textrm{if $k$ is even} \\
-\beta_{-k}+\gamma_{k-1}+\gamma_k+\gamma_{k+1} & \textrm{if $k\in \{-1,1\}$} \\
\gamma_{k-1}+\gamma_k+\gamma_{k+1} & \textrm{if $k$ is odd and $k \notin \{-1,1\}$} 
\end{cases}
\\
\label{eq:taueast}
\tau_e^\ast:~ \beta_i^\ast \mapsto \beta_{-i}^\ast+\gamma_{i}^\ast ~:~
\gamma_k^\ast \mapsto \begin{cases}
\gamma_k^\ast & \textrm{if $k$ is odd} \\
-\gamma_{k-1}^\ast+\gamma_k^\ast-\gamma_{k+1}^\ast & \textrm{if $k$ is even}
\end{cases}
\end{eqnarray}
We also record the action of $(\tau_e^\ast)^{-1}$:
$$(\tau_e^\ast)^{-1}:~ \beta_{i}^\ast \mapsto \beta_{-i}^\ast-\gamma_{-i}^\ast ~:~
\gamma_k^\ast \mapsto \begin{cases}
\gamma_k^\ast & \textrm{if $k$ is odd} \\
\gamma_{k-1}^\ast+\gamma_k^\ast+\gamma_{k+1}^\ast & \textrm{if $k$ is even}
\end{cases}$$
\end{itemize}

The work above gives us an method to prove the existence or non-existence of a stable periodic billiard path.
By the Enumeration Lemma, every closed geodesic $S(V_n)$ is the image of one of the geodesics in one of the homology classes 
$\gamma_{1-n}, \ldots, \gamma_{n-1}$ under an affine automorphism in $\langle \sigma, \tau_o, \tau_e \rangle$.
Let $x \in H_1\big(S(V_n), \Z\big)$ be the homology class of a closed geodesic on $S(V_n)$. Then there
is a $w \in \langle \sigma, \tau_o, \tau_e \rangle$ such that $x= w(\gamma_k)$. For $x$ to be stable, we must have
$\phi(x)=0 \in H_1({\mathcal D}T, \Z)$. This is equivalent to saying that
$\phi^\ast_i(x)=0$ for each $i \in \{-1,1\}$.
Since $w^\ast$ denotes the pullback action on cohomology, observe that for $i \in \{-1,1\}$, we have
$$\phi^\ast_i(x)= \phi^\ast_i \big( w(\gamma_k) \big) = \big((w^\ast)^{-1}(\phi^\ast_i)\big)(\gamma_k).$$
In particular, the geodesic in the homology class $x$ is stable if and only if 
for each $i \in \{-1,1\}$, when $(w^\ast)^{-1}(\phi^\ast_i)$ is written in the basis ${\mathcal B}^\ast$, the coefficient of $\gamma_k^\ast$ is zero.
We summarize these conclusions in the proposition below.

\begin{proposition}
\label{prop:stability}
Let $x=w(\gamma_k)\in H_1\big(S(V_n), \Z\big)$  be the homology class of a geodesic in $S(V_n)$.
The corresponding billiard path in $V_n$ is stable if and only if for each $i \in \{-1,1\}$, when $(w^\ast)^{-1}(\phi^\ast_i)$ is written in the basis ${\mathcal B}^\ast$, the coefficient of $\gamma_k^\ast$ is zero.
\end{proposition}

We simplify this proposition by noting the following:
\begin{proposition}
The cohomology class $\phi^\ast_1+\phi^\ast_{-1}=2n \beta_1+2n \beta_{-1}+2n \gamma_0^\ast$ is
invariant under the group action $\langle \sigma^\ast, \tau_o^\ast, \tau_e^\ast \rangle$.
\end{proposition}
\noindent
{\bf Proof:}
We check the least trivial case of the action of $\tau_e^\ast$, and leave the other cases to the reader:
$$\begin{array}{rcl}
\tau_e^\ast(2n \beta_1^\ast+2n \beta_{-1}^\ast+2n \gamma_0^\ast) & = &
2n(\beta_{-1}^\ast+\gamma_1^\ast)+2n(\beta_1^\ast+\gamma_{-1}^\ast)+2n(-\gamma_{-1}^\ast+\gamma_0^\ast-\gamma_1^\ast) \\
& = &
2n \beta_1^\ast+2n \beta_{-1}^\ast+2n \gamma_0^\ast.\end{array}$$
\endproof

We give a more careful statement of our final conclusions below:

\begin{lemma}[Stability]
Suppose $p$ is a periodic billiard path in $V_n$. Then, there is a lift to a closed geodesic $\tilde p$ in $S(V_n)$. Let $x$ be the homology class
of $\tilde p$. By the Enumeration Lemma, $x= w(\gamma_k)$ for some affine automorphism $w$ in the group $\langle \sigma, \tau_o, \tau_e \rangle$ and some $k \in \{1-n, \ldots, n-1\}$. 
Then, $p$ is stable if and only if $k \neq 0$ and when $(w^\ast)^{-1}(\phi^\ast_1)$ is written in the basis ${\mathcal B}^\ast$, the coefficient of $\gamma_k^\ast$ is zero.
\end{lemma}
\noindent
{\bf Proof:}
This Lemma is essentially the same as Proposition \ref{prop:stability}. We observe that stability is equivalent
to the statement that the coefficients of $\gamma_k^\ast$ is zero for the cohomology class
$(w^{-1})^\ast(\phi^\ast_1)$ and the invariant cohomology class
$(w^{-1})^\ast(\phi^\ast_1+\phi^\ast_{-1})=\phi^\ast_1+\phi^\ast_{-1}$.
The coefficient of $\gamma_k^\ast$ of $\phi^\ast_1+\phi^\ast_{-1}$ is zero 
if and only if $k \neq 0$. 
\endproof

\subsection{Instability}
Suppose that $n=2^m$ for an integer $m \geq 2$. We will use the Stability Lemma to show that $V_{n}$ has no stable periodic billiard paths.
We will show that the condition for stability given in the Stability Lemma can not hold modulo $2n=2^{m+1}$;
i.e. for the corresponding cohomology classes in $H^1\big(S(V_n), \Z/2n\Z\big)$.

The following holds for all $n$, though we will only use it for powers of two.

\begin{proposition}
\label{prop:instability}
Let $w^\ast \in \langle \sigma^\ast, \tau_o^\ast, \tau_e^\ast \rangle$. There exist odd integers $r,s$ such that
$$w^\ast(\phi^\ast_1) \equiv \sum_{i=1-n}^{n-1} c(i) \gamma_{i}^\ast \pmod{2n}$$
with coefficients given by
$$c(i)=\begin{cases}
r (i+n) & \textrm{if $i$ is odd} \\
s (i+n) & \textrm{if $i$ is even.} 
\end{cases}$$
\end{proposition}

\noindent
{\bf Proof of part 1 of Theorem \ref{thm:stability}:}
We assume $n=2^m$. Our numbers $i$ lie within the set $\{1-n,\ldots,n-1\}$ in particular
$i+n$ is never equivalent to zero modulo $2n$. Multiplication by an odd number permutes
the classes of $\Z/2n\Z=\Z/2^{m+1}\Z$ and preserves zero. Thus, $r(i+n)$ is never equivalent to
zero modulo $2^{m+1}$. Therefore, the Stability Lemma implies that there can be no stable periodic billiard paths in $V_n$.
\endproof

\noindent
{\bf Proof of Proposition \ref{prop:instability}:}
The proof is by induction in the group $\langle \sigma, \tau_o, \tau_e \rangle$. The statement is true for the identity element with $r=s=1$. See Equation \ref{eq:coverdual}.

We set up some notation.
Let ${\bf E}$ and ${\bf O}$ denote the even and odd integers in $\{1-n,\ldots,n-1\}$. 
We find it convenient (as before) to let $\gamma_{n}^\ast$ and $\gamma_{-n}^\ast$ represent the 
zero cohomology class.

Now suppose that the statement is true for the group element $w_0^\ast$ for the odd numbers $r$ and $s$. That is, suppose that in $H^1(S(V_n), \Z/2n\Z)$,
$$w_0^\ast(\phi^\ast_1)=\sum_{i \in {\bf O}} r(i+n) \gamma_i^\ast+\sum_{i \in {\bf E}} s(i+n) \gamma_i^\ast.$$
By Equation \ref{eq:sigmaast}, applying the involution $\sigma^\ast$ yields:
$$
\begin{array}{rcl}
\sigma^\ast \circ w_0^\ast(\phi^\ast_1) & = & \sum_{i \in {\bf O}} r(i+n) \gamma_{-i}^\ast+\sum_{i \in {\bf E}} s(i+n) \gamma_{-i}^\ast \\
& = & \sum_{i \in {\bf O}} r(n-i) \gamma_{i}^\ast+\sum_{i \in {\bf E}} s(n-i) \gamma_{i}^\ast \\
& = & \sum_{i \in {\bf O}} -r(n+i) \gamma_{i}^\ast+\sum_{i \in {\bf E}} -s(n+i) \gamma_{i}^\ast.
\end{array}$$
So, the equation holds for $w^\ast=\sigma^\ast \circ w_0^\ast$ for the choice of odds $-r$ and $-s$. Now consider
the action of $\tau^\ast_o$. We have 
$$
\begin{array}{rcl}
\tau_o^\ast \circ w_0^\ast(\phi^\ast_1) & = & \sum_{i \in {\bf O}} r(i+n) 
(-\gamma_{i-1}^\ast+\gamma_{i}^\ast-\gamma_{i+1}^\ast)+\sum_{i \in {\bf E}} s(i+n) \gamma_{i}^\ast \\
& = & \sum_{i \in {\bf O}} r(i+n) \gamma_{i}^\ast+\sum_{i \in {\bf E}} 
\big(-r(i-1+n)-r(i+1+n)+s(i+n)\big) \gamma_{i}^\ast\\
& = & \sum_{i \in {\bf O}} r(i+n) \gamma_i^\ast+\sum_{i \in {\bf E}} (s-2r)(i+n) \gamma_i^\ast.
\end{array}$$
So, the equation holds for $w^\ast=\tau_o^\ast  \circ w_0^\ast$
with the choice of odds $r$ and $s-2r$. Similarly, the equation holds when
$w^\ast=(\tau_o^\ast)^{-1}  \circ w_0^\ast$
with the choice of odds $r$ and $s+2r$. We need to be slightly more careful with the even twist
because of the role of the zero cohomology classes $\gamma^\ast_{-n}$ and $\gamma^\ast_{n}$. 
Working modulo $2n$, we have:
$$
\begin{array}{rcl}
\tau_e^\ast \circ w_0^\ast(\phi^\ast_1) & = & \sum_{i \in {\bf O}} r(i+n) \gamma_{i}^\ast+
\sum_{i \in {\bf E}} s(i+n) 
(-\gamma_{i-1}^\ast+\gamma_{i}^\ast-\gamma_{i+1}^\ast) \\
& = & \big(r-2s\big)\gamma_{1-n}^\ast+\big(-r-(2n-2)s\big)\gamma_{n-1}^\ast + \\
& & \sum_{i \in {\bf O} \smallsetminus \{1-n,n-1\}} \big(-s(i-1+n)-s(i+1+n)+r(i+n)\big) \gamma_{i}^\ast + \\
& & \sum_{i \in {\bf E}} s(i+n) \gamma_i^\ast \\
& = & \sum_{i \in {\bf O}} (r-2s)(i+n) \gamma_i^\ast+ \sum_{i \in {\bf E}} s(i+n) \gamma_i^\ast.
\end{array}$$
So, the equation holds when $w^\ast= \tau_e^\ast  \circ w_0^\ast$ for the choice of odds
$r-2s$ and $s$. Similarly, when $w^\ast=(\tau_e^\ast)^{-1}  \circ w_0^\ast$,
we use the odds $r+2s$ and $s$. 
\endproof

\subsection{Existence of Stable Trajectories}
Suppose that $n \geq 3$ is not a power of two. We will show the second part of Theorem \ref{thm:stability};
$V_n$ has a stable periodic billiard path. We prove this by establishing homology classes of geodesics
on $S(V_n)$ which project to stable geodesics. These homology classes are provided by the following two theorems.

\begin{theorem}[Odd case]
\label{thm:odd case}
Suppose $n \geq 3$ is odd. Then a closed geodesic in the homology class
$(\tau_e \circ \tau_o^{-1})^{\frac{n-1}{2}} \circ \tau_e^{\frac{3-n}{2}}(\gamma_{n-2})$ in $S(V_n)$ 
projects to a stable periodic billiard path in $V_n$ via the folding map $S(V_n) \to V_n$.
\end{theorem}

\begin{theorem}[Even case]
\label{thm:even case}
Suppose $n$ is even and not a power of two. Then $n=2^a b$ for an odd integer $b \geq 3$ and an integer $a \geq 1$. Let 
$w=(\tau_e \circ \tau_o^{-1})^{\frac{n}{2}} \circ \tau_o^{\frac{b-1}{2}}$. 
Then, a closed geodesic in the homology class $w(\gamma_{n-2^{a+1}})$ projects to a stable periodic billiard path in $V_n$ via the folding map
$S(V_n) \to V_n$.
\end{theorem}

\begin{remarks}
The above formulas correct the incorrect formulas in \cite{HS}. The proofs were wrong due in large part to calculation errors. For this reason, we make the calculations very explicit below. Also, work in this section was checked with Mathematica.
\end{remarks}

To prove these theorems, we investigate the action of $(w^\ast)^{-1}$ on $\phi^\ast_1$ in order to apply the Stability Lemma.
The $w$ provided by each theorem begins with a power of $\tau_e \circ \tau_o^{-1}$, which is conjugate but not equal to
a rotation of $S(V_n)$. The following is the main formula necessary in the proofs of the above results.
\begin{lemma}[Elliptic orbit]
\label{lem:elliptic}
Let $n \geq 3$ be an integer and let $k$ be an integer with $0 \leq k \leq \frac{n}{2}$. Then,
$$
\begin{array}{rcl}
\big(\tau_{o}^\ast \circ (\tau_{e}^\ast)^{-1}\big)^k(\phi^\ast_1)
& = &
n\big(1+(-1)^k\big) \beta^\ast_1+n\big(1-(-1)^k\big) \beta^\ast_{-1}+n \gamma_0^\ast + \\
& & (-1)^k \sum_{j=1}^{\min \{2k,n-1\}} \big(j  + (-1)^j (4jk-2jn+n)\big)(\gamma_j^\ast-\gamma_{-j}^\ast)+\\
& & (-1)^k \sum_{j=2k+1}^{n-1} (j - n) \big(1 + (-1)^j 4 k\big)(\gamma_j^\ast-\gamma_{-j}^\ast).
\end{array}
$$
\end{lemma}

To simplify the computations required to prove the lemma, we find it useful to break $\phi_1^\ast$ into pieces. Namely, $\phi_1^\ast$ is the sum of 
$$\Phi_A^\ast=2n \beta_1^\ast+n \gamma_0^\ast 
\quad \text{and} \quad 
\Phi_B^\ast=\sum_{j=1}^{n-1} (j-n)(\gamma_{j}^\ast-\gamma_{-j}^\ast).$$
The following two propositions give
formulas for the image of each under $\big(\tau_{o}^\ast \circ (\tau_{e}^\ast)^{-1}\big)^k$. Combining them gives the proof of the lemma above.

\begin{proposition}
Let $n \geq 3$ be an integer and suppose $0 \leq k \leq \frac{n}{2}$. Then, we have
$$
\begin{array}{rcl}
\big(\tau_{o}^\ast \circ (\tau_{e}^\ast)^{-1}\big)^k(\Phi^\ast_A)
& = &
n\big(1+(-1)^k\big) \beta^\ast_1+n\big(1-(-1)^k\big) \beta^\ast_{-1}+n \gamma_0^\ast + \\
& & n (-1)^k \sum_{j=1}^{2k} (-1)^j (\gamma_j^\ast-\gamma_{-j}^\ast).
\end{array}
$$
To make sense of the case when $k=\frac{n}{2}$, recall that we assigned
$\gamma_n^\ast$ and $\gamma_{-n}^\ast$ to be zero.
\end{proposition}
\noindent
{\bf Proof:}
We prove this by induction in $k$. The statement is true when $k=0$ by definition of $\Phi^\ast_A$. 
Now we will prove it when $k=1$. We have
$$
\begin{array}{rcl}
(\tau_{e}^\ast)^{-1}(\Phi^\ast_A)  & = & (\tau_{e}^\ast)^{-1}(2n \beta_1^\ast+n \gamma_0^\ast) \\
& = & 2n(\beta_{-1}^\ast- \gamma_{-1}^\ast)+n (\gamma_{-1}^\ast+\gamma_0^\ast+\gamma_{1}^\ast) \\
& = & 2n \beta_{-1}^\ast-n\gamma_{-1}^\ast+n\gamma_0^\ast+n\gamma_{1}^\ast.
\end{array}$$
Then applying $\tau_o^\ast$ to the above yields:
$$
\begin{array}{rcl}
\tau_o^\ast \circ (\tau_{e}^\ast)^{-1}(\Phi^\ast_A)  & = & 
\tau_o^\ast(2n \beta_{-1}^\ast-n\gamma_{-1}^\ast+n\gamma_0^\ast+n\gamma_{1}^\ast) \\
& = & 2n\beta_{1}^\ast-n(-\gamma_{-2}^\ast+\gamma_{-1}^\ast-\gamma_0^\ast)+n\gamma_0^\ast+
n(-\gamma_{0}^\ast+\gamma_{1}^\ast-\gamma_2^\ast)\\
& = & 2n \beta_1^\ast+n\gamma_{-2}^\ast-n\gamma_{-1}^\ast+n \gamma_0^\ast+n \gamma_1^\ast-n \gamma_2^\ast
\end{array}$$
This is equivalent to the formula given in the proposition when $k=1$. 

We continue by induction. Assume that the statement is true for $k$ with $1 \leq k<\frac{n}{2}$. We will prove the statement for $k+1$. For conciseness let $\Psi_k^\ast=\big(\tau_{o}^\ast \circ (\tau_{e}^\ast)^{-1}\big)^k(\Phi^\ast_A)$.
Let ${\bf E}$ and ${\bf O}$ denote the even and odd integers non-strictly between $1$ and $2k$.
By inductive hypothesis, we have 
$$\begin{array}{rcl}
(\tau_{e}^\ast)^{-1}(\Psi_k^\ast)
& = &
(\tau_{e}^\ast)^{-1}\Big(n\big(1+(-1)^k\big) \beta^\ast_1+n\big(1-(-1)^k\big) \beta^\ast_{-1}+n \gamma_0^\ast + \\
& & n (-1)^k \sum_{j=1}^{2k} (-1)^j (\gamma_j^\ast-\gamma_{-j}^\ast)\Big) \\
& = & n\big(1+(-1)^k\big) (\beta^\ast_{-1}-\gamma_{-1}^\ast)+
n\big(1-(-1)^k\big) (\beta^\ast_{1}-\gamma_{1}^\ast)+n (\gamma_{-1}^\ast+\gamma_0^\ast+\gamma_{1}^\ast)+ \\
& & n(-1)^k \sum_{j \in {\bf E}} (\gamma_{j-1}^\ast+\gamma_j^\ast+\gamma_{j+1}^\ast-\gamma_{-j-1}^\ast-\gamma_{-j}^\ast-\gamma_{1-j}^\ast)-\\
& & n (-1)^k \sum_{j \in {\bf O}}(\gamma_j^\ast-\gamma_{-j}^\ast)\\
& = &  n\big(1+(-1)^{k+1}\big) \beta^\ast_{1}+n\big(1-(-1)^{k+1}\big) \beta^\ast_{-1}+
n\gamma_0^\ast+n(-1)^k \sum_{j=1}^{2k+1} (\gamma_j^\ast-\gamma_{-j}^\ast).
\end{array}
$$
We now let ${\bf O}'$ be the odds non-strictly between $1$ and $2k+1$. We apply $\tau_o^\ast$ to the result from the above:
$$\begin{array}{rcl}
\Psi_{k+1}^\ast
& = &
\tau_{o}^\ast\Big(n\big(1+(-1)^{k+1}\big) \beta^\ast_{1}+n\big(1-(-1)^{k+1}\big) \beta^\ast_{-1}+n\gamma_0^\ast+\\
& & n(-1)^k \sum_{j \in {\bf O}'} \tau_{o}^\ast(\gamma_j^\ast-\gamma_{-j}^\ast)+
n(-1)^k \sum_{j \in {\bf E}} \tau_{o}^\ast(\gamma_j^\ast-\gamma_{-j}^\ast)\Big)\\
& = &
n\big(1+(-1)^{k+1}\big) \beta^\ast_{1}+n\big(1-(-1)^{k+1}\big) \beta^\ast_{-1}+n\gamma_0^\ast+\\
& & n(-1)^k \sum_{j \in {\bf O}'} (-\gamma_{j-1}^\ast+\gamma_j^\ast-\gamma_{j+1}^\ast+\gamma_{-j-1}^\ast-\gamma_{-j}^\ast+\gamma_{1-j}^\ast)+\\
& & n(-1)^k \sum_{j \in {\bf E}} (\gamma_j-\gamma_{-j})\\
& = & n\big(1+(-1)^{k+1}\big) \beta^\ast_{1}+n\big(1-(-1)^{k+1}\big) \beta^\ast_{-1}+n\gamma_0^\ast+ \\
& & n (-1)^{k+1} \sum_{j=1}^{2k+2} (-1)^j (\gamma_j^\ast-\gamma_{-j}^\ast).
\end{array}
$$
This proves the proposition.
\endproof

\begin{proposition}
For an integer $n \geq 3$ and an integer $k$ satisfying $0 \leq k \leq \frac{n}{2}$, we have
$$
\begin{array}{rcl}
\big(\tau_{o}^\ast \circ (\tau_{e}^\ast)^{-1}\big)^k(\Phi^\ast_B)
& = &
(-1)^k \sum_{j=1}^{\min \{2k,n-1\}} \big(j + (-1)^j (4jk-2jn)\big)(\gamma_j^\ast-\gamma_{-j}^\ast)+\\
& & (-1)^k \sum_{j=2k+1}^{n-1} (j - n) \big(1 + (-1)^j 4 k\big)(\gamma_j^\ast-\gamma_{-j}^\ast).
\end{array}
$$
\end{proposition}
\noindent
{\bf Proof:}
First we setup some notation. We let $\eta_{j}^\ast=\gamma_j^\ast-\gamma_{-j}^\ast$. 
If $j=0$ or $j=n$, then $\eta_{j}$ is zero. 

We will prove this proposition by induction in $k$. We observe the statement is true when $k=0$. Now suppose it holds for $k$
with $0 \leq k < \frac{n}{2}$.  
We define ${\mathcal A}$ and ${\mathcal B}$ so that they sum to $\big(\tau_{o}^\ast \circ (\tau_{e}^\ast)^{-1}\big)^k(\Phi^\ast_B)$:
$${\mathcal A}=(-1)^k \sum_{j=1}^{2k} j \big(1 + (-1)^j (4k-2n)\big)\eta_{j}^\ast.$$
$${\mathcal B}=(-1)^k \sum_{j=2k+1}^{n-1} (j - n) \big(1 + (-1)^j 4 k\big)\eta_{j}^\ast.$$
We first consider the action of $(\tau_{e}^\ast)^{-1}$ on ${\mathcal A}$. For each integer $i>0$ define the collections
of even and odd integers:
$${\bf E}(2i)=\{2,4,\ldots,2i\} \quad \text{and} \quad {\bf O}(2i-1)=\{1,3,\ldots,2i-1\}.$$
We have:
$$
\begin{array}{rcl}
(\tau_{e}^\ast)^{-1}({\mathcal A}) & = & 
(-1)^k (\tau_{e}^\ast)^{-1} \Big(\sum_{j \in {\bf O}(2k-1)} j (1-4k+2n)\eta_{j}^\ast+\sum_{j \in {\bf E}(2k)} j(1+4k-2n)\eta_{j}^\ast\Big) \\
& = & 
(-1)^k \Big(\sum_{j \in {\bf O}(2k-1)} j(1-4k+2n)\eta_{j}^\ast+ \\
& & \sum_{j \in {\bf E}(2k)} j(1+4k-2n)
(\eta_{j-1}^\ast+\eta_{j}^\ast+\eta_{j+1}^\ast)\Big) \\
& = & (-1)^k \Big(2k(1+4k-2n)\eta_{2k+1}^\ast+\sum_{j \in {\bf O}(2k-1)} j(3+4k-2n)\eta_{j}^\ast+\\
& &  \sum_{j \in {\bf E}(2k)} j(1+4k-2n)\eta_{j}^\ast\Big). \\
\end{array}
$$
Now we apply $\tau_o^\ast$ to the end of the above formula:
$$
\begin{array}{rcl}
\tau_o^\ast \circ (\tau_{e}^\ast)^{-1}({\mathcal A}) & = & 
(-1)^k \Big(2k(1+4k-2n)(-\eta_{2k}^\ast+\eta_{2k+1}^\ast-\eta_{2k+2}^\ast)+\\
& & \sum_{j \in {\bf O}(2k-1)} j(3+4k-2n)(-\eta_{j-1}^\ast+\eta_{j}^\ast-\eta_{j+1}^\ast)+\\
& & \sum_{j \in {\bf E}(2k)} j(1+4k-2n)\eta_{j}^\ast\Big) \\
& = & 
(-1)^k \Big(
-(2k-1)(3+4k-2n)\eta_{2k}^\ast+
2k(1+4k-2n)(\eta_{2k+1}^\ast-\eta_{2k+2}^\ast)+\\
&  & 
\sum_{j \in {\bf O}(2k-1)} j(3+4k-2n)\eta_{j}^\ast+
\sum_{j \in {\bf E}(2k-2)} j(-5-4k+2n)\eta_{j}^\ast\Big).
\end{array}
$$
One can observe that the coefficients for $\eta_j$ with $j \in \{1,\ldots,2k-1\}$
agree with our formula for $\big(\tau_{o}^\ast \circ (\tau_{e}^\ast)^{-1}\big)^{k+1}(\Phi^\ast_B)$ given in the proposition. (Note that the sign difference
comes from the fact that in the proposition, we begin with $(-1)^{k+1}$.) 

Now we will do a similar calculation for ${\mathcal B}$. Let ${\bf O}$ and ${\bf E}$
denote the odds and evens satisfying $1 \leq i \leq n-1$. 
We have:
$$
\begin{array}{rcl}
(\tau_{e}^\ast)^{-1}({\mathcal B}) & = & 
(-1)^k (\tau_{e}^\ast)^{-1} \Big(\sum_{j \in {\bf O} \smallsetminus {\bf O}(2k-1)} (j - n) (1 -4 k)\eta_{j}^\ast+ \\
& & 
\sum_{j \in {\bf E} \smallsetminus {\bf E}(2k)} (j - n) (1 + 4k)\eta_{j}^\ast\Big) \\
& = & 
(-1)^k \Big(\sum_{j \in {\bf O} \smallsetminus {\bf O}(2k-1)} (j - n) (1 -4 k)\eta_{j}^\ast+ \\
& & 
\sum_{j \in {\bf E} \smallsetminus {\bf E}(2k)} (j - n) (1 + 4k)(\eta_{j-1}^\ast+\eta_{j}^\ast+\eta_{j+1}^\ast)\Big) \\
& = & 
(-1)^k \Big((3 + 8 k - 2 n)\eta_{2k+1}^\ast+
\sum_{j \in {\bf O} \smallsetminus {\bf O}(2k+1)} (j - n) (3 +4 k)\eta_{j}^\ast \\
& & 
\sum_{j \in {\bf E} \smallsetminus {\bf E}(2k)} (j - n) (1 + 4k)\eta_{j}^\ast\Big)
\end{array}
$$
We apply $\tau_o^\ast$ to the end of the above formula:
$$
\begin{array}{rcl}
\tau_o^\ast \circ (\tau_{e}^\ast)^{-1}({\mathcal B}) & = &
(-1)^k \Big((3 + 8 k - 2 n)(-\eta_{2k}^\ast+\eta_{2k+1}^\ast-\eta_{2k+2}^\ast)+\\
& & 
\sum_{j \in {\bf O} \smallsetminus {\bf O}(2k+1)} (j - n) (3 +4 k)(-\eta_{j-1}^\ast+\eta_{j}^\ast-\eta_{j+1}^\ast) \\
& & 
\sum_{j \in {\bf E} \smallsetminus {\bf E}(2k)} (j - n) (1 + 4k)\eta_{j}^\ast\Big)
\\ & = & (-1)^k \Big((3 + 8 k - 2 n)(-\eta_{2k}^\ast+\eta_{2k+1}^\ast)+ \\
& & (-10 - 16 k + 4 n)\eta_{2k+2}^\ast+ \\
& & 
\sum_{j \in {\bf O} \smallsetminus {\bf O}(2k+1)} (j - n) (3 +4 k)\eta_{j}^\ast +
\sum_{j \in {\bf E} \smallsetminus {\bf E}(2k+2)} (j - n) (-5 - 4k)\eta_{j}^\ast\Big).
\end{array}
$$
The coefficients match with our formula for $\big(\tau_{o}^\ast \circ (\tau_{e}^\ast)^{-1}\big)^{k+1}(\Phi^\ast_B)$ so long as $j \geq 2k+3$. 

We now check the coefficients of $\eta_{j}^\ast$ for $2k \leq j \leq 2k+2$. The
coefficient of $\eta_{2k}^\ast$ of the sum $(\tau_{e}^\ast)^{-1}({\mathcal A}+{\mathcal B})$ is given by 
$$(-1)^k \big(-(2k-1)(3+4k-2n)-(3 + 8 k - 2 n)\big)=
(-1)^{k+1} (2k+2-n)\big(1+4(k+1)\big),$$
as required by the proposition. The
coefficient of $\eta_{2k+1}^\ast$ is given by:
$$(-1)^k\big(2 k (1 + 4 k - 2 n) + (3 + 8 k - 2 n)\big)=(-1)^{k+1}(2k+1)(1-4(k+1)+2n)$$
as required.
The coefficient of $\eta_{2k+2}^\ast$ is given by:
$$(-1)^k \big(-2 k (1 + 4 k - 2 n) - 2 (5 + 8 k - 2 n))\big)=
(-1)^{k+1} (2 k + 2) \big(1 + 4 (k + 1) - 2 n\big).$$
This also coincides with the formula in the proposition. We have now checked all coefficients.
\endproof

\noindent
{\bf Proof of Theorem \ref{thm:odd case} (The odd case):}
By the Stability Lemma, it is equivalent to show that the coefficient
of $\gamma_{n-2}^\ast$ is zero in 
$(\tau_e^\ast)^{\frac{n-3}{2}}\circ \big(\tau_o^\ast \circ (\tau_e^\ast)^{-1}\big)^{\frac{n-1}{2}}(\phi_1^\ast)$.
Let ${\mathcal A}=\big(\tau_o^\ast \circ (\tau_e^\ast)^{-1}\big)^{\frac{n-1}{2}}(\phi_1^\ast)$. Using the Elliptic Orbit Lemma, we compute that:
$${\mathcal A}=3 (-1)^{\frac{n-1}{2}} \gamma_{n-3}^\ast+2(n-3)(-1)^{\frac{n-1}{2}}\gamma_{n-2}^\ast+(-1)^{\frac{n-1}{2}}\gamma_{n-1}^\ast+\ldots.$$
The coefficients of the rest of the expression are irrelevant.
Observe that if $n=3$, the coefficient of $\gamma_{n-2}^\ast$ of
${\mathcal A}$ is zero. This proves this statement in this case. Otherwise, we need to apply $(\tau_e^\ast)^{\frac{n-3}{2}}$ to this expression. We find: 
$$\textstyle (\tau_e^\ast)^{\frac{n-3}{2}}({\mathcal A})=(-1)^{\frac{n-1}{2}}
\Big(
3  \big(\frac{3-n}{2}\gamma_{n-4}^\ast+\gamma_{n-3}^\ast+\frac{3-n}{2}\gamma_{n-2}^\ast\big)+
2(n-3)\gamma_{n-2}^\ast+
(\frac{3-n}{2}\gamma_{n-2}^\ast+\gamma_{n-1}^\ast)\Big)+\ldots.$$
These are the only terms which contribute to the coefficient
of $\gamma_{n-2}^\ast$, so we see this coefficient is zero.
\endproof

\noindent
{\bf Proof of Theorem \ref{thm:even case} (The even case):}
We will again use the Stability Lemma. We will show that the coefficient
of $\gamma_{n-2^{a+1}}^\ast$ is zero in 
$(\tau_o^\ast)^{\frac{1-b}{2}}\circ \big(\tau_o^\ast \circ (\tau_e^\ast)^{-1}\big)^{\frac{n}{2}}(\phi_1^\ast)$.
Let ${\mathcal A}=\big(\tau_o^\ast \circ (\tau_e^\ast)^{-1}\big)^{\frac{n}{2}}(\phi_1^\ast)$. 
Using the Elliptic Orbit Lemma, we compute that:
$${\mathcal A}=
(-1)^{\frac{n}{2}}\Big(
(-2^{a+1}-1)\gamma_{n-2^{a+1}-1}^\ast+
2(n-2^a) \gamma_{n-2^{a+1}}^\ast+
(-2^{a+1}+1)\gamma_{n-2^{a+1}+1}^\ast
\Big)+\ldots.
$$
The coefficients of the rest of the expression are irrelevant.
Then, 
$$\begin{array}{rcl}
(\tau_o^\ast)^{\frac{1-b}{2}}({\mathcal A}) & = & 
(-1)^{\frac{n}{2}}
\Big(
(-2^{a+1}-1)(\frac{b-1}{2}\gamma_{n-2^{a+1}-2}^\ast+\gamma_{n-2^{a+1}-1}^\ast+\frac{b-1}{2}\gamma_{n-2^{a+1}}^\ast)+\\
 & & 2(n-2^a) \gamma_{n-2^{a+1}}^\ast+\\
 & & (-2^{a+1}+1)(
\frac{b-1}{2}\gamma_{n-2^{a+1}}^\ast+
\gamma_{n-2^{a+1}+1}^\ast+
\frac{b-1}{2}\gamma_{n-2^{a+1}+2}^\ast
)
\Big)+\ldots.\end{array}$$
Only these terms contribute to the coefficient
of $\gamma_{n-2^{a+1}}^\ast$. We compute this coefficient:
$$(-1)^{\frac{n}{2}}\big((-n+2^a+\frac{1-b}{2})+(2n-2^{a+1})+(-n+2^a+\frac{b-1}{2})\big)=0.$$
\endproof

\newpage